\documentclass{ArticleModern2023}

\usepackage{SymbolsModern2023}
\usepackage{StylesModern2023}

\newcommand{\TreeT}{\TermT}
\newcommand{\TreeS}{\TermS}
\newcommand{\TreeR}{\TermR}

\newcommand{\SetU}{\mathfrak{U}}

\newcommand{\PWord}{\mathfrak{p}}
\newcommand{\QWord}{\mathfrak{q}}
\newcommand{\RWord}{\mathfrak{r}}

\newcommand{\VariableCount}{\mathrm{vc}}
\newcommand{\OperationCount}{\mathrm{oc}}
\newcommand{\Op}{\mathfrak{op}}
\newcommand{\Action}{\cdot}
\newcommand{\MonoidProduct}{\operatorname{\cdot}}
\newcommand{\MonoidProductExtension}{\operatorname{\BBar{\MonoidProduct}}}
\newcommand{\MonoidUnit}{\mathrm{e}}
\newcommand{\TrivialMonoid}{\mathcal{E}}
\newcommand{\TrivialClone}{\mathcal{T}}
\newcommand{\EndClone}{\mathcal{M}}
\newcommand{\SetLetters}{\mathcal{L}}
\renewcommand{\P}{\mathbf{P}}
\newcommand{\T}{\mathbf{T}}
\newcommand{\PLeq}{\bm{\Leq}}
\newcommand{\Identity}{\mathrm{u}}
\newcommand{\Pigmentation}{\mathrm{p}}
\newcommand{\PWordEmpty}{\bm{\epsilon}}
\newcommand{\PSymbol}{\mathbb{P}}
\newcommand{\EquivCovering}{\sim}
\newcommand{\Covering}{\leadsto}
\newcommand{\CoveringRT}{\ll}

\newcommand{\Interpretation}{\mathsf{int}}
\newcommand{\Frontier}{\mathsf{fr}}
\newcommand{\RightComb}{\mathsf{rc}}
\newcommand{\Sort}{\mathsf{sort}}
\newcommand{\First}{\mathsf{first}}
\newcommand{\Max}{\operatorname{\downarrow}}
\newcommand{\Superposition}[1]{\bm{\left[}#1\bm{\right]}}
\newcommand{\RemoveFirst}{\mathsf{r}\First}
\newcommand{\RelationSetQ}{\mathrel{\mathfrak{Q}}}

\newcommand{\WInc}{\mathsf{WInc}}
\newcommand{\Inc}{\mathsf{Inc}}
\newcommand{\Magn}{\mathsf{Magn}}
\newcommand{\Arra}{\mathsf{Arra}}
\newcommand{\Stal}{\mathsf{Stal}}
\newcommand{\Pill}{\mathsf{Pill}}

\newcommand{\Witness}[3]{
    \tikz[baseline=(line.base)]{
        \node[inner sep=2pt,outer sep=0pt](line){$#3$};
        \draw[shorten <=1mm,shorten >=1mm](line.south west)--(line.south east);
        \ifthenelse{\equal{#1}{y}}{
            \node[circle,draw=ColA,fill=ColA!60,minimum size=1mm,inner sep=0pt,xshift=1mm]
                at(line.south west){};
        }
        {}
        \ifthenelse{\equal{#1}{n}}{
            \node[xshift=1mm,font=\scriptsize,inner sep=0pt]
                at(line.south west){$\ColA{\times}$};
        }
        {}
        \ifthenelse{\equal{#2}{y}}{
            \node[circle,draw=ColB,fill=ColB!60,minimum size=1mm,inner sep=0pt,xshift=-1mm]
                at(line.south east){};
        }
        {}
        \ifthenelse{\equal{#2}{n}}{
            \node[xshift=-1mm,font=\scriptsize,inner sep=0pt]
                at(line.south east){$\ColB{\times}$};
        }
        {}
    }
}

\newcommand{\LeftCombTwo}[5]{
    \begin{tikzpicture}[Centering,xscale=.2,yscale=.29]
        \node[LeafST](0)at(0.00,-3.33){$#3$};
        \node[LeafST](2)at(2.00,-3.33){$#4$};
        \node[LeafST](4)at(4.00,-1.67){$#5$};
        \node[NodeST](1)at(1.00,-1.67){$#1$};
        \node[NodeST](3)at(3.00,0.00){$#2$};
        \draw[Edge](0)--(1);
        \draw[Edge](1)--(3);
        \draw[Edge](2)--(1);
        \draw[Edge](4)--(3);
        \node(r)at(3.00,1.5){};
        \draw[Edge](r)--(3);
    \end{tikzpicture}}

\newcommand{\BinaryCorolla}[3]{
    \begin{tikzpicture}[Centering,xscale=.3,yscale=.35]
        \node[LeafST](0)at(0.00,-1.50){$#2$};
        \node[LeafST](2)at(2.00,-1.50){$#3$};
        \node[NodeST](1)at(1.00,0.00){$#1$};
        \draw[Edge](0)--(1);
        \draw[Edge](2)--(1);
        \node(r)at(1.00,1.25){};
        \draw[Edge](r)--(1);
    \end{tikzpicture}}

\Title{%
    Clones of pigmented words and realizations of special classes of monoids%
    \vspace{-1ex}%
}
\TitleShort{Clones of pigmented words}
\Author{Samuele Giraudo}
\AuthorShort{S. Giraudo}
\Subjects{%
    05E16, 
    08B15, 
    16S15, 
    18M80, 
    68R15. 
}
\Date{\today}
\Address{%
    Université du Québec à Montréal, LACIM, \\
    Pavillon Président-Kennedy, 201 Avenue du Président-Kennedy, Montréal, H2X~3Y7, Canada.
}
\Email{giraudo.samuele@uqam.ca}
\Funding{%
    This research has been partially supported by the projects CARPLO (ANR-20-CE40-0007),
    LambdaComb (ANR-21-CE48-0017) of the Agence nationale de la recherche, and
    by the Natural Sciences and Engineering Research Council of Canada (RGPIN-2024-04465).
}
\Abstract{%
    Clones are specializations of operads forming powerful instruments to describe varieties
    of algebras wherein repeating variables are allowed in their equations. They allow us in
    this way to realize and study a large range of algebraic structures. A functorial
    construction from the category of monoids to the category of clones is introduced. The
    obtained clones involve words on positive integers where letters are accompanied by
    elements of a monoid. By considering quotients of these structures, we construct a
    complete hierarchy of clones involving some families of combinatorial objects. This
    provides clone realizations of some known and some new special classes of monoids as
    among others the variety of left-regular bands, bounded semilattices, and regular band
    monoids.
}

\begin{document}

\MakeFirstPage

\section{Introduction}
Given a variety of algebras specified by a set of fundamental operations together with
equations between the operations, an important question consists in deciding if two terms
describing compound operations are equivalent. For instance, in the variety of groups, the
two operations $\Par{x_1, x_2} \mapsto \Par{x_1 \cdot x_2}^{-1}$ and $\Par{x_1, x_2} \mapsto
x_2^{-1} \cdot x_1^{-1}$ compute always both the same value, where $\Par{x_1, x_2} \mapsto
x_1 \cdot x_2$ is the multiplication operation and $x_1 \mapsto x_1^{-1}$ is the inverse
operation of groups. This general question is known as the word problem and in some cases,
term rewrite systems~\cite{BN98,BKV03} offer solutions by orienting in a suitable way the
equations which define the variety in order to form a terminating and confluent rewrite
system.

While the word problem is in general undecidable, this inherent undecidability does not
obstruct the development of tools capable of resolving specific instances. Rather than
focusing on finding the optimal orientation or completion of the equations within a variety,
an alternative approach involves encoding compound operations by using combinatorial
objects. In this context, the functional composition can be interpreted as a relevant
operation on these objects. Within this framework, operads~\cite{LV12,Men15,Gir18} emerge as
valuable instruments to facilitate these abstractions, called operad realizations of a
variety. An illustrating example can be found in the realization of the variety of pre-Lie
algebras in terms of rooted trees~\cite{CL01} and grafting operations on such trees.
Besides, operads are also great tools for tackling problems originating from combinatorics.
Indeed, by endowing a set of combinatorial objects with an operad structure, we obtain a
framework for enumerating~\cite{Gir20b} and generating~\cite{Gir19} their elements. This is
based on presentations by generators and equations of the operads to study and more
precisely on their orientations in order to form, here again, terminating and confluent term
rewrite systems.

Despite their broad utility, operads have limitations, particularly when dealing with
varieties that are defined through equations with repeating variables. This issue arises,
for instance, in the variety of groups, lattices, or flexible algebras, where natural
descriptions of these varieties require equations involving repeated inputs. Although it is
feasible to capture a certain part of such varieties by working with operads in the category
of vector spaces on a field of zero characteristic and by considering some tricks to encode
equations with repeating variables by linear combinations of linear terms (like in the case
of the variety of flexible algebras~\cite{May72}), operads are not the ideal instrument in
this context. Some other devices have been developed for these purposes. Examples include
abstract clones~\cite{Coh65,Tay93,MMT18}, Lawvere theories~\cite{Law63,ARVL10}, and monads
with arities~\cite{EM65,HP07,BMW12}. The aim of this work is to create bridges between the
theory of abstract clones ---called simply ``clones'' here henceforth--- and combinatorics.
To our knowledge, contrary to what operad theory has experienced since its rebirth in the
1990s~\cite{Lod96}, not many such connections have been established in the existing
literature. We have opted to work with clones rather than with Lawvere theories or monads
with arities because clones can be perceived as generalized operads with minor distinctions.
Since as presented above, the connections between operads and combinatorics are now very
clear and well-established (see also~\cite{CL01,Gir15,Gir18,Gir20}), we anticipate that new
significant connections between clones and combinatorics could be unearthed.

In an initial, humble, and modest first step in this direction, we introduce a new
combinatorial recipe to build clones of combinatorial objects. More precisely, given a
monoid $\Monoid$, we construct a clone $\P(\Monoid)$ involving $\Monoid$-pigmented words,
which are some words of integers whose each letter is accompanied by an element of
$\Monoid$. The variety of algebras described by $\P(\Monoid)$, called variety of
$\Monoid$-pigmented monoids, bears similarities to the variety of algebras described by the
operad $\T(\Monoid)$, where $\T$ is a construction from monoids to operads introduced
in~\cite{Gir15}. More specifically, the variety of $\Monoid$-pigmented algebras has an extra
nullary fundamental operation (playing the role of a unit) and some equations involving it
compared to the variety of algebras described by $\T(\Monoid)$. For this reason, the present
work can be seen as a continuation and a generalization of~\cite{Gir15}, but in the context
of clones rather than of operads.

The clone $\P(\Monoid)$ is rich enough to contain some notable quotients. In order to
construct quotients of $\P(\Monoid)$, we consider clone congruences $\Equiv$ of
$\P(\Monoid)$ each coming with a so-called $\PSymbol$-symbol to decide whether two
$\Monoid$-pigmented words are $\Equiv$-equivalent. A $\PSymbol$-symbol for a clone
congruence $\Equiv$ is a map sending an $\Monoid$-pigmented word to a representative of its
$\Equiv$-equivalence class. Such maps enable us to obtain concrete realizations and
presentations by generators and equations of quotients of $\P(\Monoid)$. The studied
quotients of $\P(\Monoid)$ fit into a diagram of surjective clone morphisms generalizing
some lattices of varieties of special classes of monoids (see~\cite{GLV22}) and of
semigroups (see~\cite{Eva71,SVV09,KKP11}). In particular, we obtain as main results clone
realizations of commutative monoids, left-regular bands, bounded semilattices, and regular
bands. These clone realizations allow us to solve the word problem in these varieties by
using algorithms akin to those developed in~\cite{SS82,NS00} for idempotent semigroups.

This paper is organized as follows. Section~\ref{sec:clones_varieties} contains preliminary
notions about terms, clones and free clones, presentations of clones, and varieties of
algebras. In particular, we show Proposition~\ref{prop:clone_presentation} which is an
important result to establish presentations of clones. Next, in
Section~\ref{sec:pigmented_monoids}, we introduce the varieties of $\Monoid$-pigmented
monoids and describe the construction $\P$. By
Theorem~\ref{thm:clone_presentation_pigmented_monoids}, the main result of this section, we
show that $\P(\Monoid)$ is a clone realization of the variety of $\Monoid$-pigmented
monoids. In Section~\ref{sec:construction_quotients} we introduce some tools to investigate
quotient clones of $\P(\Monoid)$. In particular, we introduce the concept of
$\PSymbol$-symbol specific to our context and its relationships with clone congruences by
way of Propositions~\ref{prop:p_symbol_fiber_equivalence_relations},
\ref{prop:p_symbol_congruence}, and~\ref{prop:p_symbol_composition}. We show also with
Proposition~\ref{prop:clone_realization_p_symbol} how to obtain a concrete description of a
quotient of $\P(\Monoid)$ by a congruence $\Equiv$ admitting a $\PSymbol$-symbol
$\PSymbol_{\Equiv}$. Continuing this, two clone congruences $\Equiv_\Sort$ and
$\Equiv_{\First_k}$, $k \geq 0$, are introduced. These congruences as well as some of their
compositions are used to build the quotient clones $\WInc(\Monoid)$, $\Arra_k(\Monoid)$, $k
\geq 0$, and $\Inc_k$, $k \geq 0$. By Propositions~\ref{prop:presentation_winc},
\ref{prop:presentation_arra}, and~\ref{prop:presentation_inc}, we describe presentations of
these clones. Finally, Section~\ref{sec:hierarchy} contains the most technical results under
a combinatorial point of view. Here, we construct three quotients of $\P(\Monoid)$ by clone
congruences defined by intersecting some of the congruences $\Equiv_\Sort$ and
$\Equiv_{\First_k}$, $k \geq 0$. The main results are formed by
Theorems~\ref{thm:p_symbol_magn}, \ref{thm:p_symbol_stal}, and~\ref{thm:p_symbol_pill},
describing realizations of these clones, and Theorems~\ref{thm:presentation_magn},
\ref{thm:presentation_stal}, and~\ref{thm:presentation_pill}, giving presentations for these
clones. In particular, we obtain here a clone realization of the variety of regular bands
which seems new at the best of our knowledge. This text ends with a list of open questions
and future research directions.

\subsubsection*{Acknowledgments}
We would like to express our sincere gratitude to the Editor for selecting such a highly
competent and insightful reviewer. We extend our deepest thanks to the Reviewer, whose
detailed and constructive feedback has been invaluable in improving the entire article.

\subsubsection*{General notations and conventions}
For any integers $i$ and $j$, $[i, j]$ denotes the set $\{i, i + 1, \dots, j\}$. For any
integer $i$, $[i]$ denotes the set $[1, i]$ and $\HanL{i}$ denotes the set $[0, i]$. For any
set $A$, $A^*$ is the set of words on $A$. For any $w \in A^*$, $\Length(w)$ is the length
of $w$, and for any $i \in [\Length(w)]$, $w(i)$ is the $i$-th letter of $w$. For any $a \in
A$, $|w|_a$ is the number of occurrences of $a$ in $w$. The only word of length $0$ is the
empty word $\epsilon$. For any $i \leq j \in [\Length(w)]$, $w(i, j)$ is the word $w(i) w(i
+ 1) \dots w(j)$. The word $\Reverse(w)$ is the mirror image $w(\Length(w)) \dots w(1)$ of
$w$. Given two words $w$ and $w'$, the concatenation of $w$ and $w'$ is denoted by $ww'$ or
by $w \Conc w'$.

\section{Clones and realizations of varieties} \label{sec:clones_varieties}
This preliminary section contains the main definitions and notions about abstract clones,
free abstract clones, presentations of abstract clones by generators and equations,
varieties of algebras, and clone realizations of varieties of algebras.

\subsection{Abstract clones}
In this part, we set our notations and main notions about abstract clones. Let us begin
with graded sets.

\subsubsection{Graded sets}
A \Def{graded set} is a disjoint union $G := \bigsqcup_{n \geq 0} G(n)$. For any $x \in G$,
the unique integer $n \geq 0$ such that $x \in G(n)$ is the \Def{arity} of $x$, denoted by
$|x|$. If for any $n \geq 0$, $G(n)$ is finite, then $G$ is \Def{combinatorial}. In this
case, the \Def{sequence of sizes} of $G$ is the sequence $\Par{\# G(n)}_{n \geq 0}$. Let
$G'$ be another graded set. A map $\phi : G \to G'$ is a \Def{graded set morphism} if $\phi$
preserves the arities. Besides, if for any $n \geq 0$, $G'(n) \subseteq G(n)$, then $G'$ is
a \Def{graded subset} of~$G$. A binary relation $\RelationSet$ on $G$ is a \Def{graded set
binary relation} on $G$ if $\RelationSet$ preserves the arities. The \Def{quotient} of $G$
by a graded set equivalence relation $\Equiv$ is the graded set $G /_{\Equiv}$ defined for
any $n \geq 0$ by $G /_{\Equiv}(n) := \Bra{[x]_{\Equiv} : x \in G(n)}$ where $[x]_{\Equiv}$
is the $\Equiv$-equivalence class of $x \in G$.

\subsubsection{Abstract clones}
Abstract clones are devices which can be used to describe algebraic
structures~\cite{Coh65,Neu70,Tay93,MMT18} (see also~\cite{Fuj20} for a point of view from
universal algebra~\cite{BS81}). An \Def{abstract clone} (or \Def{clone} for short) $\Clone$
is a graded set $\Clone$ endowed with maps
\begin{equation}
    -\Superposition{-, \dots, -}_{n, m} : \Clone(n) \times \Clone(m)^n \to \Clone(m),
\end{equation}
where $n, m \geq 0$, called \Def{superposition maps}, and with distinguished elements
$\Unit_{i, n} \in \Clone(n)$, where $n \geq 1$ and $i \in [n]$, called \Def{projections}.
This data has to satisfy, for any $x \in \Clone(n)$, $n \geq 0$, $y_1, \dots, y_n \in
\Clone(m)$, $m \geq 0$, $z _1, \dots, z_m \in \Clone(k)$, $k \geq 0$, and $i \in [n]$, the
relations
\begin{equation} \label{equ:clone_relation_1}
    \Unit_{i, n} \Superposition{y_1, \dots, y_n}_{n, m} = y_i,
\end{equation}
\begin{equation} \label{equ:clone_relation_2}
    x \Superposition{\Unit_{1, n}, \dots, \Unit_{n, n}}_{n, n} = x,
\end{equation}
\begin{equation} \label{equ:clone_relation_3}
    \Par{x \Superposition{y_1, \dots, y_n}_{n, m}} \Superposition{z_1, \dots, z_m}_{m, k}
    = x \Superposition{
        y_1 \Superposition{z_1, \dots, z_m}_{m, k},
        \dots, y_n \Superposition{z_1, \dots, z_m}_{m, k}}_{n, k}.
\end{equation}
To lighten the notation when the context is clear, we shall drop the indices of the
superposition maps in order to write $x \Superposition{y_1, \dots, y_n}$ instead of $x
\Superposition{y_1, \dots, y_n}_{n, m}$ for any $x \in \Clone(n)$, $n \geq 0$ and $y_1,
\dots, y_n \in \Clone(m)$, $m \geq 0$. In the same way, we shall write $\Unit_i$ instead of
$\Unit_{i, n}$ for any $n \geq 1$ and $i \in [n]$ when the value of $n$ is clear or not
significant.

Observe that for any $0 \leq n \leq m$, there is a map $\iota_{n, m} : \Clone(n) \to
\Clone(m)$ such that for any $x \in \Clone(n)$, $\iota_{n, m}(x) := x
\Superposition{\Unit_{1, m}, \dots, \Unit_{n, m}}$. It is easy to check that $\iota$ is an
injection. Therefore, in each set $\Clone(m)$, there is a copy of the elements of
$\Clone(n)$, seen in $\Clone(m)$ as elements of arity~$m$. Observe also that for any $n \geq
0$, $\iota_{n, n}$ is the identity map on $\Clone(n)$, and that for any $0 \leq n \leq m
\leq k$, the relation $\iota_{m, k} \circ \iota_{n, m} = \iota_{n, k}$ holds.

The \Def{trivial clone} is the clone $\TrivialClone$ such that for any $n \geq 0$,
$\TrivialClone(n)$ is a singleton. Observe that there is no choice for the definition of the
superposition maps of $\TrivialClone$. Let $\Clone'$ be another clone. A graded set morphism
$\phi : \Clone \to \Clone'$ is a \Def{clone morphism} if, for any $n \geq 1$ and $i \in
[n]$, $\phi$ sends the projection $\Unit_{i, n}$ of $\Clone$ to the projection $\Unit'_{i,
n}$ of $\Clone'$, and for any $x \in \Clone(n)$, $n \geq 0$, and any $y_1, \dots, y_n \in
\Clone(m)$, $m \geq 0$,
\begin{equation}
    \phi\Par{x \Superposition{y_1, \dots, y_n}}
    = \phi(x) \Superposition{\phi\Par{y_1}, \dots, \phi\Par{y_n}}.
\end{equation}
Besides, if $\Clone'$ is a graded subset of $\Clone$ such that $\Clone'$ contains the
projections of $\Clone$, and $\Clone'$ is closed under the superposition maps of $\Clone$,
then $\Clone'$ is a \Def{subclone} of $\Clone$. Given $S \subseteq \Clone$, the subclone of
$\Clone$ \Def{generated} by $S$ is the smallest subclone $\Clone^S$ of $\Clone$ containing
$S$. When $\Clone^S = \Clone$, $S$ is a \Def{generating set} of $\Clone$. A \Def{clone
congruence} of $\Clone$ is a graded set equivalence relation $\Equiv$ on $\Clone$ such that
for any $x, x' \in \Clone(n)$, $n \geq 0$, and any $y_1, y'_1, \dots, y_n, y'_n \in
\Clone(m)$, $m \geq 0$, if $x \Equiv x'$ and $y_1 \Equiv y'_1$, \dots, $y_n \Equiv y'_n$,
then
\begin{math}
    x \Superposition{y_1, \dots, y_n} \Equiv x' \Superposition{y'_1, \dots, y'_n}.
\end{math}
The \Def{quotient} of $\Clone$ by $\Equiv$ is the clone on the graded set $\Clone
/_{\Equiv}$ such that for any $x \in \Clone(n)$, $n \geq 0$, $y_1, \dots, y_n \in
\Clone(m)$, $m \geq 0$, the superposition maps of $\Clone /_{\Equiv}$ satisfy
\begin{equation}
    [x]_{\Equiv} \Superposition{\Han{y_1}_{\Equiv}, \dots, \Han{y_n}_{\Equiv}}
    = \Han{x \Superposition{y_1, \dots, y_n}}_{\Equiv},
\end{equation}
and for any $n \geq 1$ and $i \in [n]$, the projection $\Unit_{i, n}$ of $\Clone /_{\Equiv}$
is the $\Equiv$-equivalence class of the projection~$\Unit_{i, n}$ of~$\Clone$.

\subsubsection{Algebras over clones}
Let $\Clone$ be a clone. An \Def{algebra} over $\Clone$ (or a \Def{$\Clone$-algebra} for
short) is a structure $\Par{\Algebra, \Op}$ where $\Algebra$ is a set and for any $x \in
\Clone(n), n \geq 0$, $\Op(x)$ is a map from $\Algebra^n$ to $\Algebra$ satisfying the
following relations. For any $a_1, \dots, a_m \in \Algebra$, $m \geq 0$, $i \in [m]$, $x \in
\Clone(n)$, $n \geq 0$, and $y_1, \dots, y_n \in \Clone(m)$,
\begin{equation}
    \Op\Par{\Unit_{i, m}}\Par{a_1, \dots, a_m} = a_i,
\end{equation}
\begin{equation}
    \Op\Par{x \Superposition{y_1, \dots, y_n}}\Par{a_1, \dots, a_m}
    = \Op(x)\Par{
        \Op\Par{y_1}\Par{a_1, \dots, a_m}, \dots, \Op\Par{y_n}\Par{a_1, \dots, a_m}}.
\end{equation}
In other terms, each $x \in \Clone(n)$, $n \geq 0$, gives rise to an operation $\Op(x)$ on
$\Algebra$ with $n$ inputs and one output, and the functional composition of such operations
is coherent with the superposition maps of~$\Clone$.

Algebras over clones admit the following equivalent description. For any set $\Algebra$, let
$\EndClone_\Algebra$ be the graded set such that for any $n \geq 0$, $\EndClone_\Algebra(n)$
is the set of maps from $\Algebra^n$ to $\Algebra$. We endow $\EndClone_\Algebra$ with the
superposition maps satisfying, for any $n \geq 0$, $f \in \EndClone_\Algebra(n)$, $g_1,
\dots, g_n \in \EndClone_\Algebra(m)$, and $a_1, \dots, a_m \in \Algebra$,
\begin{equation}
    \Par{f \Superposition{g_1, \dots, g_n}} \Par{a_1, \dots, a_m}
    = f\Par{g_1\Par{a_1, \dots, a_m}, \dots, g_n\Par{a_1, \dots, a_m}}.
\end{equation}
We also endow $\EndClone_\Algebra$, for any $n \geq 1$ and $i \in [n]$, with the projection
$\Unit_{i, n}$ satisfying, for any $a_1, \dots, a_n \in \Algebra$, $\Unit_{i, n}\Par{a_1,
\dots, a_n} = a_i$. This endows $\EndClone_\Algebra$ with the structure of a clone, called
the \Def{clone of endomorphisms}. Given this, for any clone $\Clone$, a structure
$\Par{\Algebra, \Op}$ is a $\Clone$-algebra if and only if $\Op$ is a clone morphism from
$\Clone$ to~$\EndClone_\Algebra$.

\subsection{Terms and free clones}
In order to describe free clones, we need to introduce some notions and combinatorics about
terms. The reason behind this is that the elements of free clones can be described as terms
and their superposition maps as graftings in terms.

\subsubsection{Terms} \label{subsubsec:terms}
A \Def{signature} is a graded set $\GeneratingSet$. Its elements are called \Def{operation
symbols}. Let $\SetVariables$ be the set $\Bra{\VarX_i : i \in \N \setminus \{0\}}$. Any
element of $\SetVariables$ is a \Def{variable}. For any $n \in \N$, let the subset
$\SetVariables_n$ of $\SetVariables$ consisting of the variables $\VarX_i$ such that $i \in
[n]$. A \Def{$\GeneratingSet$-term} (or simply \Def{term} when the context is clear) is
recursively either a variable or a pair $\Par{\GenG, \Par{\TreeT_1, \dots, \TreeT_k}}$,
where $\GenG \in \GeneratingSet(k)$, $k \geq 0$, and $\TreeT_1$, \dots, $\TreeT_k$ are
$\GeneratingSet$-terms. For convenience, we shall write $\GenG \Superposition{\TreeT_1,
\dots, \TreeT_k}$ instead of $\Par{\GenG, \Par{\TreeT_1, \dots, \TreeT_k}}$. From this
definition, any $\GeneratingSet$-term can be interpreted as a rooted planar tree where
internal nodes are decorated by operation symbols and leaves are decorated by variables. The
graded set of $\GeneratingSet$-terms is denoted by $\SetTerms(\GeneratingSet)$ where, for
any $n \geq 0$, $\SetTerms(\GeneratingSet)(n)$ is a copy of the set of
$\GeneratingSet$-terms having all variables belonging to~$\SetVariables_n$.

Let $\TreeT$ be a $\GeneratingSet$-term. The \Def{operation count} $\OperationCount(\TreeT)$
of $\TreeT$ is the number of internal nodes of $\TreeT$ seen as a tree, that is, the number
of operation symbols of $\TreeT$ counted with multiplicities. The \Def{variable count}
$\VariableCount(\TreeT)$ of $\TreeT$ is the number of variables of $\TreeT$ counted with
multiplicities. If $\GeneratingSet'$ is a signature and $\phi : \GeneratingSet \to
\GeneratingSet'$ is a graded set morphism, we denote by $\widehat{\phi} :
\SetTerms(\GeneratingSet) \to \SetTerms\Par{\GeneratingSet'}$ the map such that, for any
$\TreeT \in \SetTerms(\GeneratingSet)$, $\widehat{\phi}(\TreeT)$ is the
$\GeneratingSet'$-term obtained by replacing each decoration $\GenG \in \GeneratingSet$ of
an internal node of $\TreeT$ by $\phi(\GenG)$.

For instance, by setting $\GeneratingSet$ as the signature satisfying $\GeneratingSet =
\GeneratingSet(0) \sqcup \GeneratingSet(2) \sqcup \GeneratingSet(3)$ with $\GeneratingSet(0)
= \{\GenA\}$, $\GeneratingSet(2) = \{\GenB, \GenC\}$, and $\GeneratingSet(3) = \{\GenD\}$,
\begin{equation}
    \TreeT :=
    \GenD \Superposition{
        \GenB \Superposition{\GenD \Superposition{\VarX_1, \GenA, \VarX_1}, \VarX_3},
        \GenA,
        \GenD \Superposition{\GenC \Superposition{\VarX_5, \VarX_3}, \VarX_4, \GenA}}
\end{equation}
is a $\GeneratingSet$-term. The treelike representation of $\TreeT$ is
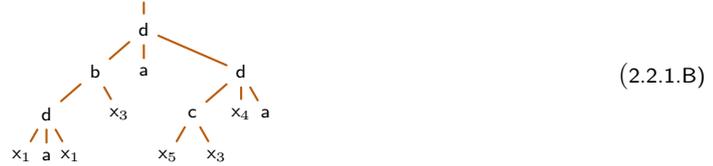
\begin{equation} \label{equ:example_term}
    \begin{tikzpicture}[Centering,xscale=0.32,yscale=0.16]
        \node[LeafST](0)at(0.00,-10.50){$\VarX_1$};
        \node[LeafST](10)at(8.00,-10.50){$\VarX_3$};
        \node[LeafST](12)at(9.00,-7.00){$\VarX_4$};
        \node[NodeST](13)at(10.00,-7.00){$\GenA$};
        \node[NodeST](2)at(1.00,-10.50){$\GenA$};
        \node[LeafST](3)at(2.00,-10.50){$\VarX_1$};
        \node[LeafST](5)at(4.00,-7.00){$\VarX_3$};
        \node[NodeST](7)at(5.00,-3.50){$\GenA$};
        \node[LeafST](8)at(6.00,-10.50){$\VarX_5$};
        \node[NodeST](1)at(1.00,-7.00){$\GenD$};
        \node[NodeST](11)at(9.00,-3.50){$\GenD$};
        \node[NodeST](4)at(3.00,-3.50){$\GenB$};
        \node[NodeST](6)at(5.00,0.00){$\GenD$};
        \node[NodeST](9)at(7.00,-7.00){$\GenC$};
        \draw[Edge](0)--(1);
        \draw[Edge](1)--(4);
        \draw[Edge](10)--(9);
        \draw[Edge](11)--(6);
        \draw[Edge](12)--(11);
        \draw[Edge](13)--(11);
        \draw[Edge](2)--(1);
        \draw[Edge](3)--(1);
        \draw[Edge](4)--(6);
        \draw[Edge](5)--(4);
        \draw[Edge](7)--(6);
        \draw[Edge](8)--(9);
        \draw[Edge](9)--(11);
        \node(r)at(5.00,3){};
        \draw[Edge](r)--(6);
    \end{tikzpicture}
\end{equation}
This term has $8$ as operation count and $6$ as variable count.

There is at this stage a little subtlety to remark: a $\GeneratingSet$-term $\TreeT$ gives
rise to different elements of the graded set $\SetTerms(\GeneratingSet)$ depending on the
arity attributed to it. For instance, the term defined in~\eqref{equ:example_term} can among
others be an element of $\SetTerms(\GeneratingSet)(5)$ or of $\SetTerms(\GeneratingSet)(6)$,
both distinct from each other.

\subsubsection{Free clones}
Given a signature $\GeneratingSet$, $\TreeT \in \SetTerms(\GeneratingSet)(n)$, $n \geq 0$,
and $\TreeT'_1, \dots, \TreeT'_n \in \SetTerms(\GeneratingSet)(m)$, $m \geq 0$, the
\Def{substitution} of $\TreeT'_1, \dots, \TreeT'_n$ in $\TreeT$ is the $\GeneratingSet$-term
$\TreeT \Superposition{\TreeT'_1, \dots, \TreeT'_n}$ obtained by simultaneously replacing
for all $i \in [n]$ all occurrences of the variables $\VarX_i$ in $\TreeT$ by~$\TreeT'_i$.
For instance, by considering the signature $\GeneratingSet$ defined at the end of
Section~\ref{subsubsec:terms}, we have the substitution
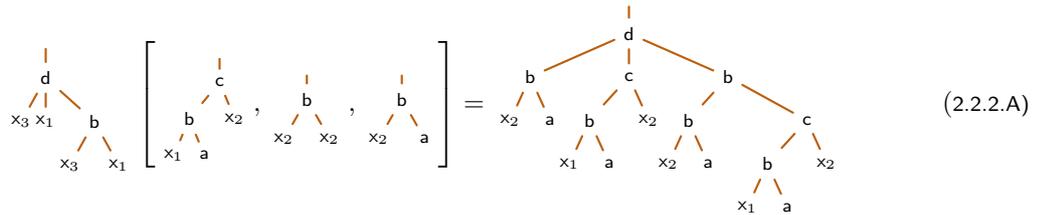
\begin{equation}
    \begin{tikzpicture}[Centering,xscale=0.32,yscale=0.29]
        \node[LeafST](0)at(0.00,-2.00){$\VarX_3$};
        \node[LeafST](2)at(1.00,-2.00){$\VarX_1$};
        \node[LeafST](3)at(2.00,-4.00){$\VarX_3$};
        \node[LeafST](5)at(4.00,-4.00){$\VarX_1$};
        \node[NodeST](1)at(1.00,0.00){$\GenD$};
        \node[NodeST](4)at(3.00,-2.00){$\GenB$};
        \draw[Edge](0)--(1);
        \draw[Edge](2)--(1);
        \draw[Edge](3)--(4);
        \draw[Edge](4)--(1);
        \draw[Edge](5)--(4);
        \node(r)at(1.00,1.75){};
        \draw[Edge](r)--(1);
    \end{tikzpicture}
    \Superposition{
    \LeftCombTwo{\GenB}{\GenC}{\VarX_1}{\GenA}{\VarX_2},
    \BinaryCorolla{\GenB}{\VarX_2}{\VarX_2},
    \BinaryCorolla{\GenB}{\VarX_2}{\GenA}}
    =
    \begin{tikzpicture}[Centering,xscale=0.26,yscale=0.16]
        \node[LeafST](0)at(0.00,-7.20){$\VarX_2$};
        \node[LeafST](11)at(10.00,-10.80){$\GenA$};
        \node[LeafST](13)at(12.00,-14.40){$\VarX_1$};
        \node[LeafST](15)at(14.00,-14.40){$\GenA$};
        \node[LeafST](17)at(16.00,-10.80){$\VarX_2$};
        \node[LeafST](2)at(2.00,-7.20){$\GenA$};
        \node[LeafST](4)at(3.00,-10.80){$\VarX_1$};
        \node[LeafST](6)at(5.00,-10.80){$\GenA$};
        \node[LeafST](8)at(7.00,-7.20){$\VarX_2$};
        \node[LeafST](9)at(8.00,-10.80){$\VarX_2$};
        \node[NodeST](1)at(1.00,-3.60){$\GenB$};
        \node[NodeST](10)at(9.00,-7.20){$\GenB$};
        \node[NodeST](12)at(11.00,-3.60){$\GenB$};
        \node[NodeST](14)at(13.00,-10.80){$\GenB$};
        \node[NodeST](16)at(15.00,-7.20){$\GenC$};
        \node[NodeST](3)at(6.00,0.00){$\GenD$};
        \node[NodeST](5)at(4.00,-7.20){$\GenB$};
        \node[NodeST](7)at(6.00,-3.60){$\GenC$};
        \draw[Edge](0)--(1);
        \draw[Edge](1)--(3);
        \draw[Edge](10)--(12);
        \draw[Edge](11)--(10);
        \draw[Edge](12)--(3);
        \draw[Edge](13)--(14);
        \draw[Edge](14)--(16);
        \draw[Edge](15)--(14);
        \draw[Edge](16)--(12);
        \draw[Edge](17)--(16);
        \draw[Edge](2)--(1);
        \draw[Edge](4)--(5);
        \draw[Edge](5)--(7);
        \draw[Edge](6)--(5);
        \draw[Edge](7)--(3);
        \draw[Edge](8)--(7);
        \draw[Edge](9)--(10);
        \node(r)at(6.00,3){};
        \draw[Edge](r)--(3);
    \end{tikzpicture}
\end{equation}
of $\GeneratingSet$-terms.

The \Def{free clone} on $\GeneratingSet$ is the clone $\SetTerms(\GeneratingSet)$ on the
graded set of the $\GeneratingSet$-terms endowed with the following superposition maps and
projections. Given $\TreeT \in \SetTerms(\GeneratingSet)(n)$, $n \geq 0$, and $\TreeT'_1,
\dots, \TreeT'_n \in \SetTerms(\GeneratingSet)(m)$, $m \geq 0$ the superposition $\TreeT
\Superposition{\TreeT'_1, \dots, \TreeT'_n}$ is the substitution of $\TreeT'_1, \dots,
\TreeT'_n$ in $\TreeT$. Moreover, for any $n \geq 1$ and $i \in [n]$, the projection
$\Unit_{i, n}$ is the $\GeneratingSet$-term~$\VarX_i \in \SetTerms(\GeneratingSet)(n)$.

For any signature $\GeneratingSet'$ and any graded set morphism $\phi : \GeneratingSet \to
\GeneratingSet'$, the graded set morphism $\widehat{\phi} : \SetTerms(\GeneratingSet) \to
\SetTerms\Par{\GeneratingSet'}$ defined in Section~\ref{subsubsec:terms}, by an inductive
argument, becomes a clone morphism w.r.t.\ the clone structure on
$\SetTerms(\GeneratingSet)$ and $\SetTerms\Par{\GeneratingSet'}$ just defined.

\subsection{Clone presentations and varieties}
This preliminary section ends by setting up some notions about varieties of algebras and
clone presentations.

\subsubsection{Evaluation maps}
If $\Clone$ is a clone, $\Clone$ is in particular a graded set and thus, a signature.
Therefore, the free clone on $\Clone$ is a well-defined clone $\SetTerms(\Clone)$. The
\Def{evaluation map} of $\Clone$ is the graded set morphism
\begin{math}
    \Eval_\Clone : \SetTerms(\Clone) \to \Clone
\end{math}
recursively defined, for any $n \geq 1$ and $i \in [n]$ by
\begin{equation}
    \Eval_\Clone\Par{\VarX_i} := \Unit_{i, n},
\end{equation}
and, for any $\GenG \in \Clone(n)$, $n \geq 0$, and $\TreeT_1, \dots, \TreeT_n \in
\SetTerms(\Clone)(m)$, $m \geq 0$, by
\begin{equation} \label{equ:evaluation_map}
    \Eval_\Clone\Par{\GenG \Superposition{\TreeT_1, \dots, \TreeT_n}}
    := \GenG \Superposition{\Eval_\Clone\Par{\TreeT_1}, \dots, \Eval_\Clone\Par{\TreeT_n}},
\end{equation}
where the superposition of the right-hand side of~\eqref{equ:evaluation_map} is the one of
$\Clone$. Note that by induction on the terms, $\Eval_\Clone$ is a clone morphism.

\subsubsection{Varieties and presentations}
A \Def{variety} is a pair $\Variety := (\GeneratingSet, \RelationSet)$ such that
$\GeneratingSet$ is a signature and $\RelationSet$ is a graded set binary relation on
$\SetTerms(\GeneratingSet)$. Any pair $\Par{\TreeT, \TreeT'}$ of $\GeneratingSet$-terms such
that $\TreeT \ \RelationSet \ \TreeT'$ is an \Def{equation} of~$\Variety$. The \Def{clone
congruence generated} by $\RelationSet$ is the smallest clone congruence
$\Equiv_\RelationSet$ of $\SetTerms(\GeneratingSet)$ containing $\RelationSet$.

A \Def{presentation} of a clone $\Clone$ is a variety $\Variety := (\GeneratingSet,
\RelationSet)$ such that $\Clone$ is isomorphic as a clone to $\SetTerms(\GeneratingSet)
/_{\Equiv_\RelationSet}$. A presentation $\Variety := (\GeneratingSet, \RelationSet)$ of
$\Clone$ is \Def{finitely equationally axiomatizable} if $\RelationSet$ is finite. An
\Def{algebra} over the variety $\Variety$ is an algebra over the clone admitting $\Variety$
as presentation.

The following statement is an important tool used in the sequel to establish clone
presentations.

\begin{Statement}{Proposition}{prop:clone_presentation}
    Let $\Clone$ be a clone, $\Variety := (\GeneratingSet, \RelationSet)$ be a variety, and
    $\phi : \GeneratingSet \to \Clone$ be a graded set morphism. If $\phi(\GeneratingSet)$
    is a generating set of $\Clone$ and, for any $\TreeT, \TreeT' \in
    \SetTerms(\GeneratingSet)$, $\TreeT \Equiv_\RelationSet \TreeT'$ if and only if
    \begin{math}
        \Eval_\Clone\Par{\widehat{\phi}\Par{\TreeT}}
        =
        \Eval_\Clone\Par{\widehat{\phi}\Par{\TreeT'}},
    \end{math}
    then $\Variety$ is a presentation of $\Clone$.
\end{Statement}
\begin{Proof}
    Let us denote by $\theta : \SetTerms(\GeneratingSet) \to \Clone$ the map $\Eval_\Clone
    \circ \widehat{\phi}$ where $\widehat{\phi}$ is defined as in
    Section~\ref{subsubsec:terms}. Since $\Eval_\Clone : \SetTerms(\Clone) \to \Clone$ is a
    surjective clone morphism and $\phi(\GeneratingSet)$ is a generating set of $\Clone$,
    $\theta$ is a surjective clone morphism. Moreover, the fact that, by hypothesis, for any
    $\GeneratingSet$-terms $\TreeT$ and $\TreeT'$ such that $\TreeT \Equiv_\RelationSet
    \TreeT'$, $\theta(\TreeT) = \theta\Par{\TreeT'}$ holds, $\theta$ induces a well-defined
    surjective clone morphism $\overline{\theta} : \SetTerms(\GeneratingSet)
    /_{\Equiv_\RelationSet} \to \Clone$. Besides, if $\Han{\TreeT}_{\Equiv_\RelationSet}$
    and $\Han{\TreeT'}_{\Equiv_\RelationSet}$ are two $\Equiv_\RelationSet$-equivalence
    classes of $\GeneratingSet$-terms such that
    \begin{math}
        \overline{\theta}\Par{\Han{\TreeT}_{\Equiv_\RelationSet}}
        =
        \overline{\theta}\Par{\Han{\TreeT'}_{\Equiv_\RelationSet}},
    \end{math}
    then for any $\TreeT \in \Han{\TreeT}_{\Equiv_\RelationSet}$ and $\TreeT' \in
    \Han{\TreeT'}_{\Equiv_\RelationSet}$, we have $\theta(\TreeT) = \theta\Par{\TreeT'}$.
    This implies by using the hypothesis of the statement of the proposition that $\TreeT
    \Equiv_\RelationSet \TreeT'$. Therefore, $\Han{\TreeT}_{\Equiv_\RelationSet} =
    \Han{\TreeT'}_{\Equiv_\RelationSet}$, showing that $\overline{\theta}$ is injective. We
    have shown that $\overline{\theta}$ is a clone isomorphism between
    $\SetTerms(\GeneratingSet) /_{\Equiv_\RelationSet}$ and $\Clone$, implying the statement
    of the proposition.
\end{Proof}

The following statement serves as a tool for describing presentations of clones, which are
defined as quotients of other clones whose presentations are already known.

\begin{Statement}{Proposition}{prop:presentation_quotient}
    Let $\GeneratingSet$ be a signature, and $\RelationSet$ and $\RelationSet'$ be two
    graded set binary relations on $\SetTerms(\GeneratingSet)$. Let $\Equiv$ be the clone
    congruence of $\SetTerms(\GeneratingSet) /_{\Equiv_\RelationSet}$ generated by
    $\Han{\TreeT}_{\Equiv_\RelationSet} \Equiv \Han{\TreeT'}_{\Equiv_\RelationSet}$ whenever
    $\TreeT \ \RelationSet' \ \TreeT'$. The map
    \begin{equation}
        \phi :
        \SetTerms(\GeneratingSet) /_{\Equiv_{\RelationSet \cup \RelationSet'}}
        \to
        \Par{\SetTerms(\GeneratingSet) /_{\Equiv_\RelationSet}} /_{\Equiv}
    \end{equation}
    defined, for any $\TreeT \in \SetTerms(\GeneratingSet)$, by
    \begin{equation}
        \phi\Par{\Han{\TreeT}_{\Equiv_{\RelationSet \cup \RelationSet'}}}
        := \Han{\Han{\TreeT}_{\Equiv_{\RelationSet}}}_{\Equiv}
    \end{equation}
    is a clone isomorphism.
\end{Statement}
\begin{Proof}
    Let $g_1 : \SetTerms(\GeneratingSet) \to \SetTerms(\GeneratingSet)
    /_{\Equiv_{\RelationSet}}$ and $g_2 : \SetTerms(\GeneratingSet)
    /_{\Equiv_{\RelationSet}} \to \Par{\SetTerms(\GeneratingSet) /_{\Equiv_{\RelationSet}}}
    /_{\Equiv}$ be the canonical projection maps, and set $f := g_2 \circ g_1$. For the rest
    of this proof, $\TreeT$ and $\TreeT'$ are any two $\GeneratingSet$-terms. By denoting by
    $\ker(f)$ the kernel of $f$, observe that the property $\TreeT \, \ker(f) \, \TreeT'$ is
    equivalent to $f(\TreeT) = f(\TreeT')$, which is in turn equivalent to
    $g_2\Par{g_1(\TreeT)} = g_2\Par{g_1(\TreeT')}$, which is, by definition of $g_2$,
    finally equivalent to $g_1(\TreeT) \, \Equiv \, g_1(\TreeT')$.

    Let us show that $\ker(f)$ is equal to $\Equiv_{\RelationSet \cup \RelationSet'}$. Since
    $f$ is a clone morphism, $\ker(f)$ is a clone congruence on $\SetTerms(\GeneratingSet)$.
    Moreover, $\RelationSet \subseteq \ker(f)$ because $\TreeT \, \RelationSet \, \TreeT'$
    implies $\TreeT \Equiv_{\RelationSet} \TreeT'$ and thus $g_1(\TreeT) = g_1(\TreeT')$.
    Also, if $\TreeT \, \RelationSet' \, \TreeT'$, then by definition of $\Equiv$,
    $g_1(\TreeT)\Equiv g_1(\TreeT')$. Hence, $\TreeT \, \ker(f) \, \TreeT'$. Therefore,
    $\ker(f)$ contains $\RelationSet \cup \RelationSet'$, and, since $\Equiv_{\RelationSet
    \cup \RelationSet'}$ is the smallest clone congruence of $\SetTerms(\GeneratingSet)$
    containing $\RelationSet \cup \RelationSet'$, we have $\Equiv_{\RelationSet \cup
    \RelationSet'} \, \subseteq \, \ker(f)$. Conversely, as $\Equiv_{\RelationSet} \,
    \subseteq \, \Equiv_{\RelationSet \cup \RelationSet'}$, the Correspondence
    Theorem~\cite{BS81} yields a clone congruence $\overline{\Equiv_{\RelationSet \cup
    \RelationSet'}}$ on $\SetTerms(\GeneratingSet) /_{\Equiv_{\RelationSet}}$ such that
    $g_1(\TreeT) \, \overline{\Equiv_{\RelationSet \cup \RelationSet'}} \, g_1(\TreeT')$ if
    and only if $\TreeT \, \Equiv_{\RelationSet \cup \RelationSet'} \, \TreeT'$. Since
    $\Equiv_{\RelationSet \cup \RelationSet'}$ contains $\RelationSet'$, $\TreeT \,
    \RelationSet' \, \TreeT'$ implies $g_1(\TreeT) \, \overline{\Equiv_{\RelationSet \cup
    \RelationSet'}} \, g_1(\TreeT')$. Hence, $\overline{\Equiv_{\RelationSet \cup
    \RelationSet'}}$ contains the generating pairs of $\Equiv$ and since $\Equiv$ is the
    smallest clone congruence on $\SetTerms(\GeneratingSet) /_{\Equiv_{\RelationSet}}$
    containing these pairs, we have $\Equiv \subseteq \overline{\Equiv_{\RelationSet \cup
    \RelationSet'}}$. Therefore, $\TreeT \, \ker(f) \, \TreeT'$ implies $g_1(\TreeT) \Equiv
    g_1(\TreeT')$, which implies $g_1(\TreeT) \, \overline{\Equiv_{\RelationSet \cup
    \RelationSet'}} \, g_1(\TreeT')$, which implies finally $\TreeT \, \Equiv_{\RelationSet
    \cup \RelationSet'} \, \TreeT'$. This proves that $\ker(f) \, \subseteq \,
    \Equiv_{\RelationSet \cup \RelationSet'}$.

    Finally, since $g_1$ and $g_2$ are surjective, $f$ is also surjective. Hence, the image
    of $f$ is $\Par{\SetTerms(\GeneratingSet)/_{\Equiv_{\RelationSet}}}/_{\Equiv}$. By the
    First Isomorphism Theorem~\cite{BS81}, $f$ induces a clone isomorphism $\phi' :
    \SetTerms(\GeneratingSet) /_{\ker(f)} \to \Par{\SetTerms(\GeneratingSet)
    /_{\Equiv_{\RelationSet}}} /_{\Equiv}$ sending $\Han{\TreeT}_{\ker(f)}$ to $f(\TreeT)$.
    By using the property $\ker(f) \, = \, \Equiv_{\RelationSet \cup \RelationSet'}$ we have
    just shown and $f(\TreeT) = \Han{\Han{\TreeT}_{\Equiv_{\RelationSet}}}_{\Equiv}$, we
    obtain that $\phi'$ sends $\Han{\TreeT}_{\Equiv_{\RelationSet \cup \RelationSet'}}$ to
    $\Han{\Han{\TreeT}_{\Equiv_{\RelationSet}}}_{\Equiv}$. Hence, the images of
    $\Han{\TreeT}_{\Equiv_{\RelationSet \cup \RelationSet'}}$ under $\phi$ and $\phi'$
    coincide, so that $\phi' = \phi$.
\end{Proof}

Let $\Clone$ be a clone and $\RelationSetQ$ be a graded set binary relation on $\Clone$.
Proposition~\ref{prop:presentation_quotient} allows us to describe a presentation of $\Clone
/_{\Equiv_{\RelationSetQ}}$ from a presentation $\Par{\GeneratingSet, \RelationSet}$ of
$\Clone$ as follows. First, fix a clone isomorphism $\psi : \SetTerms(\GeneratingSet)
/_{\Equiv_\RelationSet} \to \Clone$. Let $\pi : \SetTerms(\GeneratingSet) \to
\SetTerms(\GeneratingSet) /_{\Equiv_\RelationSet}$ be the canonical projection map and
$\varepsilon : \SetTerms(\GeneratingSet) \to \Clone$ be the map $\varepsilon := \psi \circ
\pi$. Finally, let $\sigma : \Clone \to \SetTerms(\GeneratingSet)$ be a graded section of
$\varepsilon$, that is, a map preserving the arities such that $\varepsilon \circ \sigma$ is
the identity map on $\Clone$. By setting $\RelationSetQ'$ as the graded set binary relation
$\Bra{\Par{\sigma(c),\sigma(c')} : (c,c') \in \, \RelationSetQ}$ on
$\SetTerms(\GeneratingSet)$, it follows from Proposition~\ref{prop:presentation_quotient}
that the quotient $\Clone /_{\Equiv_{\RelationSetQ}}$ admits the presentation
$\Par{\GeneratingSet, \RelationSet \, \cup \RelationSetQ'}$. Note that $\RelationSetQ'$ has
one pair for each pair of $\RelationSetQ$.

\subsubsection{Clone realizations of varieties}
In the other direction, given a variety $\Variety$, any clone admitting $\Variety$ as
presentation is a \Def{clone realization} of $\Variety$ (see~\cite{Neu70}).

For instance, let the variety $\Variety := (\GeneratingSet, \RelationSet)$ where
$\GeneratingSet$ is the signature satisfying $\GeneratingSet = \GeneratingSet(2) =
\{\Meet\}$ and $\RelationSet$ is the binary relation on $\SetTerms(\GeneratingSet)$
satisfying
\begin{equation}
    \Meet \Superposition{\Meet \Superposition{\VarX_1, \VarX_2}, \VarX_3}
    \ \RelationSet \
    \Meet \Superposition{\VarX_1, \Meet \Superposition{\VarX_2, \VarX_3}},
\end{equation}
\begin{equation}
    \Meet \Superposition{\VarX_1, \VarX_2}
    \ \RelationSet \
    \Meet \Superposition{\VarX_2, \VarX_1},
\end{equation}
\begin{equation}
    \Meet \Superposition{\VarX_1, \VarX_1}
    \ \RelationSet \
    \VarX_1.
\end{equation}
This is the variety of semilattices. The clone realization $\Clone :=
\SetTerms(\GeneratingSet) /_{\Equiv_\RelationSet}$ admits the following concrete
description. For any $n \geq 0$, $\Clone(n)$ is a copy of the set of nonempty subsets of
$[n]$. The superposition maps of $\Clone$ satisfy, for any $n \geq 0$, $\SetU \in
\Clone(n)$, and $\SetU'_1, \dots, \SetU'_n \in \Clone(m)$, $m \geq 0$,
\begin{equation}
    \SetU \Superposition{\SetU'_1, \dots, \SetU'_n} = \bigcup_{i \in \SetU} \SetU'_i,
\end{equation}
and for any $n \geq 1$ and $i \in [n]$, the projection $\Unit_{i, n}$ is~$\{i\}$. Any
algebra over $\Clone$ comprises all term operations of a semilattice.

\section{Pigmented monoids and clones of pigmented words} \label{sec:pigmented_monoids}
We introduce here the variety of pigmented monoids which is roughly speaking a variety
wherein algebras are monoids endowed with monoid endomorphisms indexed by another monoid
$\Monoid$ ---the pigments--- with some extra structure. A clone realization $\P(\Monoid)$ of
this variety involving some particular words as main combinatorial objects is described.

\subsection{Pigmented monoids}
Let us describe the variety of pigmented monoids and browse some examples of such
structures having some combinatorial interest.

\subsubsection{Varieties of pigmented monoids} \label{subsubsec:varieties_pigmented_monoids}
Let $(\Monoid, \MonoidProduct, \MonoidUnit)$ be a monoid. Recall that $\MonoidProduct$ is an
associative binary operation and that $\MonoidUnit$ is the unit w.r.t.\ the operation
$\MonoidProduct$. We denote by $\TrivialMonoid$ the trivial monoid, that is the monoid
having $\MonoidUnit$ as unique element.

The \Def{variety of $\Monoid$-pigmented monoids} (or simply \Def{pigmented monoids} when the
context is clear) is the variety $\Par{\GeneratingSet_\Monoid, \RelationSet_\Monoid}$ such
that
\begin{math}
    \GeneratingSet_\Monoid
    := \GeneratingSet_\Monoid(0) \sqcup \GeneratingSet_\Monoid(1)
    \sqcup \GeneratingSet_\Monoid(2)
\end{math}
where
\begin{equation}
    \GeneratingSet_\Monoid(0) := \{\Identity\},
    \quad
    \GeneratingSet_\Monoid(1) := \Bra{\Pigmentation_\alpha : \alpha \in \Monoid},
    \quad
    \GeneratingSet_\Monoid(2) := \{\Product\},
\end{equation}
and $\RelationSet_\Monoid$ is the binary relation on $\SetTerms\Par{\GeneratingSet_\Monoid}$
satisfying
\begin{equation} \label{equ:pigmented_monoids_1}
    \Product\Superposition{\Product\Superposition{\VarX_1, \VarX_2}, \VarX_3}
    \ \RelationSet_\Monoid \
    \Product\Superposition{\VarX_1, \Product\Superposition{\VarX_2, \VarX_3}},
\end{equation}
\begin{equation} \label{equ:pigmented_monoids_2}
    \Product\Superposition{\Identity, \VarX_1}
    \ \RelationSet_\Monoid \
    \VarX_1
    \ \RelationSet_\Monoid \
    \Product\Superposition{\VarX_1, \Identity},
\end{equation}
\begin{equation} \label{equ:pigmented_monoids_3}
    \Pigmentation_\alpha\Superposition{\Product\Superposition{\VarX_1, \VarX_2}}
    \ \RelationSet_\Monoid \
    \Product\Superposition{
        \Pigmentation_\alpha\Superposition{\VarX_1},
        \Pigmentation_\alpha\Superposition{\VarX_2}
    },
\end{equation}
\begin{equation} \label{equ:pigmented_monoids_4}
    \Pigmentation_\alpha\Superposition{\Identity}
    \ \RelationSet_\Monoid \
    \Identity,
\end{equation}
\begin{equation} \label{equ:pigmented_monoids_5}
    \Pigmentation_{\alpha_1}\Superposition{\Pigmentation_{\alpha_2}\Superposition{\VarX_1}}
    \ \RelationSet_\Monoid \
    \Pigmentation_{\alpha_1 \MonoidProduct \alpha_2} \Superposition{\VarX_1},
\end{equation}
\begin{equation} \label{equ:pigmented_monoids_6}
    \Pigmentation_\MonoidUnit\Superposition{\VarX_1} \ \RelationSet_\Monoid \ \VarX_1,
\end{equation}
for any $\alpha, \alpha_1, \alpha_2 \in \Monoid$.

To simplify the notations, since the fundamental operation $\Product$ is binary, by
Equation~\eqref{equ:pigmented_monoids_1} associative, and by
Equation~\eqref{equ:pigmented_monoids_2} admits $\Identity$ as unit, we shall sometimes
treat $\Product$ as an infix operator which associates from right to left. This means that
for any $\TreeT_1, \dots, \TreeT_k \in \SetTerms\Par{\GeneratingSet_\Monoid}$, $k \geq 1$,
\begin{math}
    \TreeT_1 \Product \TreeT_2 \Product \dots \Product \TreeT_{k - 1}
    \Product \TreeT_k
\end{math}
specifies the $\GeneratingSet_\Monoid$-term
\begin{math}
    \Product \Superposition{
        \TreeT_1,
        \Product \Superposition{
            \TreeT_2,
            \Product \Han{
                \dots
                \Product \Superposition{
                    \TreeT_{k - 1},
                    \TreeT_k
                }
                \dots
            }
        }
    },
\end{math}
and for $k = 0$ it stands for~$\Identity$.

Let $(\Algebra, \Op)$ be an algebra over the variety of $\Monoid$-pigmented monoids. By
denoting by $\Product$ the binary product $\Op(\Product)$, by $\Identity$ the operation
symbol $\Op(\Identity)$, and for any $\alpha \in \Monoid$, by $\Pigmentation_\alpha$ the
unary product $\Op\Par{\Pigmentation_\alpha}$, the following properties hold.
\begin{enumerate}[label=(\roman*)]
    \item By~\eqref{equ:pigmented_monoids_1} and~\eqref{equ:pigmented_monoids_2},
    $\Par{\Algebra, \Product, \Identity}$ is a monoid.
    \item By~\eqref{equ:pigmented_monoids_3} and~\eqref{equ:pigmented_monoids_4}, each
    $\Pigmentation_\alpha$, $\alpha \in \Monoid$, is a monoid endomorphism of
    $\Par{\Algebra, \Product, \Identity}$.
    \item By~\eqref{equ:pigmented_monoids_5} and~\eqref{equ:pigmented_monoids_6}, for any
    $\alpha \in \Monoid$, the map $\Action : \Monoid \times \Algebra \to \Algebra$ defined
    by $\alpha \Action x := \Pigmentation_\alpha(x)$ is a left monoid action of $\Monoid$ on
    $\Algebra$.
\end{enumerate}
Any such structure $\Par{\Algebra, \Product, \Identity, \Par{\Pigmentation_\alpha}_{\alpha
\in \Monoid}}$ is an \Def{$\Monoid$-pigmented monoid} (or simply \Def{pigmented monoid} when
the context is clear).

For instance, any $\Z / 2 \Z$-pigmented monoid is a set $\Algebra$ endowed with an
associative product $\Product$ and two unary operations $\Pigmentation_0$ and
$\Pigmentation_1$ such that $\Product$ admits a unit $\Identity \in \Algebra$,
$\Pigmentation_0$ is the identity map on $\Algebra$, and for any $x, x_1, x_2 \in \Algebra$,
\begin{math}
    \Pigmentation_1\Par{x_1 \Product x_2}
    = \Pigmentation_1\Par{x_1} \Product \Pigmentation_1\Par{x_2},
\end{math}
\begin{math}
    \Pigmentation_1(\Identity) = \Identity,
\end{math}
and
\begin{math}
    \Pigmentation_1\Par{\Pigmentation_1\Par{x}} = x.
\end{math}
In other terms, a $\Z / 2 \Z$-pigmented monoid is a monoid endowed with an involutive monoid
endomorphism. Similarly, a $\Par{\{0, 1\}, \times, 1}$-pigmented monoid is a monoid endowed
with an idempotent monoid endomorphism.

There is a connection between $\Monoid$-pigmented monoids and semimodules. Indeed, we can
view an $\Monoid$-pigmented monoid $\Algebra$ as a semimodule where its underlying monoid
$\Algebra$ is noncommutative and its underlying semiring is the monoid $\Monoid$, obtained
by dropping the additive structure of the underlying semiring and the associated
axiomatizing relations.

A variation of $\Monoid$-pigmented monoids has been considered in~\cite{Gir15} (see
also~\cite[Chap.~4]{Gir18}) as algebras over some operads. In this cited work, the
considered variety admits $\GeneratingSet_\Monoid \setminus \{\Identity\}$ as signature and
$\RelationSet_\Monoid$ deprived of Equations~\eqref{equ:pigmented_monoids_2}
and~\eqref{equ:pigmented_monoids_4} as axiomatizing relation.

\subsubsection{Examples}
Let us consider the following examples of pigmented monoids.
\begin{enumerate}[label=(E\arabic*)]
    \item Let $\Algebra := \Par{\N^*, \Conc, \epsilon, \Par{\Pigmentation_\alpha}_{\alpha
    \in \N}}$ where $\Conc$ is the concatenation product and for any $\alpha \in \N$,
    $\Pigmentation_\alpha$ is the map sending any word to its subword made of the letters
    greater than or equal to $\alpha$. This quadruple is an $\Monoid$-pigmented monoid where
    $\Monoid := (\N, \max, 0)$. For instance,
    \begin{equation}
        \Pigmentation_2\Par{0 0 1 \ColMix{15}{5} \ColMix{15}{2} 1 \ColMix{15}{3}
            \Conc \ColMix{15}{4} 1 \ColMix{15}{2} 0 0}
        = \ColMix{15}{52342}
        = \Pigmentation_2\Par{0 0 1 \ColMix{15}{5} \ColMix{15}{2} 1 \ColMix{15}{3}}
            \Conc \Pigmentation_2\Par{\ColMix{15}{4} 1 \ColMix{15}{2} 0 0}.
    \end{equation}
    \item Let $\Algebra := \Par{\Z^*, \Conc, \epsilon, \Par{\Pigmentation_\alpha}_{\alpha
    \in \Z}}$ where $\Conc$ is the concatenation product and for any $\alpha \in \Z$,
    $\Pigmentation_\alpha$ is the map sending any word to the word obtained by incrementing
    by $\alpha$ its letters. This quadruple is an $\Monoid$-pigmented monoid where $\Monoid
    := (\Z, +, 0)$. For instance, by denoting by $\BBar{n}$ any negative integer having $n$
    as absolute value,
    \begin{equation}
        \Pigmentation_{\BBar{3}}\Par{2 4 \BBar{3} 0 \Conc \BBar{2} 6 4}
        = \BBar{1} 1 \BBar{6} \BBar{3} \BBar{5} 3 1
        = \Pigmentation_{\BBar{3}}\Par{2 4 \BBar{3} 0}
            \Conc \Pigmentation_{\BBar{3}}\Par{\BBar{2} 6 4}.
    \end{equation}
    \item Let $\Algebra := \Par{\K \AAngle{\VarZ}, +, 0, \Par{\Pigmentation_\alpha}_{\alpha
    \in \K}}$ where $\K$ is a field with multiplication denoted by $\MonoidProduct$, $\K
    \AAngle{\VarZ}$ is the space of formal power series on the parameter $\VarZ$, and for
    any $\alpha \in \K$, $\Pigmentation_\alpha$ is the map sending any series to the series
    obtained by multiplying its coefficients by $\alpha$. This quadruple is an
    $\Monoid$-pigmented monoid where $\Monoid := (\K, \MonoidProduct, 1)$.
    \item Generalizing the previous example, let $\Algebra := \Par{V, +, 0,
    \Par{\Pigmentation_\alpha}_{\alpha \in \K}}$ where $V$ is a vector space on a field $\K$
    with multiplication denoted by $\MonoidProduct$, and for any $\alpha \in \K$ and $v \in
    V$, $\Pigmentation_\alpha(v) = \alpha \MonoidProduct v$. This quadruple is an
    $\Monoid$-pigmented monoid where $\Monoid := (\K, \MonoidProduct, 1)$.
\end{enumerate}

\subsection{Clone of pigmented words}
We describe now a construction taking as input a monoid $\Monoid$ and outputting a clone
$\P(\Monoid)$ on the graded set of so-called $\Monoid$-pigmented words. We show some first
properties of this construction $\P$, such as the fact that it is a functor from the
category of monoids to the category of clones, and describe a generating set
of~$\P(\Monoid)$.

\subsubsection{Pigmented words}
Let $S$ be a nonempty set. An \Def{$S$-pigmented letter} (or \Def{pigmented letter} when
the context is clear) is a pair $(i, \alpha)$, denoted by $i^\alpha$, where $\alpha \in S$
and $i$ is a positive integer. We call $i$ (resp.\ $\alpha$) the \Def{value} (resp.\ the
\Def{pigment}) of $i^\alpha$. Let $\SetLetters_S$ be the set of $S$-pigmented letters. An
\Def{$S$-pigmented word} (or \Def{pigmented word} when the context is clear) of \Def{arity}
$n$, $n \geq 0$, is a word $\PWord$ on $\SetLetters_S$ such that all values of the pigmented
letters of $\PWord$ belong to $[n]$. The only $S$-pigmented word of arity $0$ is the empty
word, denoted by $\PWordEmpty$ in this context. For instance, $\PWord := 2^a 1^a 1^b 6^a$ is
an $\{a, b, c, d\}$-pigmented word of arity $17$.

\subsubsection{Construction}
Let $(\Monoid, \MonoidProduct, \MonoidUnit)$ be a monoid. Let $\P(\Monoid)$ be the
graded set of $\Monoid$-pigmented words. Let
\begin{math}
    \MonoidProductExtension : \Monoid \times \P(\Monoid) \to \P(\Monoid)
\end{math}
be the map defined for any $\alpha \in \Monoid$ and any $\Monoid$-pigmented word
$i_1^{\alpha_1} \dots i_\ell^{\alpha_\ell}$ by
\begin{equation}
    \alpha \ \MonoidProductExtension \ i_1^{\alpha_1} \dots i_\ell^{\alpha_\ell} :=
    i_1^{\alpha \MonoidProduct \alpha_1} \dots i_\ell^{\alpha \MonoidProduct \alpha_\ell}.
\end{equation}
Observe that this yields a left $\Monoid$-action on $\P(\Monoid)$, which moreover satisfies
\begin{math}
    \alpha \ \MonoidProductExtension \ \Par{\PWord_1 \Conc \dots \Conc \PWord_n}
    =
    \Par{\alpha \ \MonoidProductExtension \PWord_1}
    \Conc \dots \Conc
    \Par{\alpha \ \MonoidProductExtension \PWord_n}
\end{math}
for any $\PWord_1, \dots, \PWord_n \in \P(\Monoid)(m)$, $n, m \geq 0$. Let us moreover endow
$\P(\Monoid)$ with the superposition maps defined for any $i_1^{\alpha_1} \dots
i_\ell^{\alpha_\ell} \in \P(\Monoid)(n)$, $n \geq 0$, and $\PWord_1, \dots, \PWord_n \in
\P(\Monoid)(m)$, $m \geq 0$, by
\begin{equation}
    i_1^{\alpha_1} \dots i_\ell^{\alpha_\ell} \Superposition{\PWord_1, \dots, \PWord_n}
    := \Par{\alpha_1 \ \MonoidProductExtension \ \PWord_{i_1}}
    \ \Conc \ \dots \ \Conc \
        \Par{\alpha_\ell \ \MonoidProductExtension \ \PWord_{i_\ell}}.
\end{equation}

A pigmented word $\PWord = i_1^{\alpha_1} \dots i_\ell^{\alpha_\ell}$ can be interpreted as
a linear combination $\alpha_1 i_1 + \dots + \alpha_\ell i_\ell$, where the addition $+$ is
noncommutative and is expressed multiplicatively. From this perspective, the operation
$\, \MonoidProductExtension \,$ and the superposition maps are modeled, respectively, after
the scalar multiplication and the substitution of such linear combinations.

For instance, by denoting by $A^*$ the free monoid $\Par{A^*, \Conc, \epsilon}$ generated
by $A := \Bra{a, b, c}$, we have in $\P\Par{A^*}$,
\begin{align}
    2^{\ColB{ba}} 2^{\ColB{aa}} 4^{\ColB{baa}} 3^{\ColB{\epsilon}}
    \Superposition{
        2^{\ColA{b}} 1^{\ColA{aa}},
        1^{\ColA{bbb}} 1^\epsilon 2^{\ColA{b}},
        2^{\ColA{aa}} 2^{\ColA{a}},
        \PWordEmpty}
    & =
    \Par{\ColB{ba} \MonoidProductExtension 1^{\ColA{bbb}} 1^\epsilon 2^{\ColA{b}}}
    \Conc \Par{\ColB{aa} \MonoidProductExtension 1^{\ColA{bbb}} 1^\epsilon 2^{\ColA{b}}}
    \Conc \Par{\ColB{baa} \MonoidProductExtension \PWordEmpty}
    \Conc \Par{\ColB{\epsilon} \MonoidProductExtension 2^{\ColA{aa}} 2^{\ColA{a}}}
    \\
    & = 1^{\ColB{ba}\ColA{bbb}} 1^{\ColB{ba}} 2^{\ColB{ba}\ColA{b}}
    \ \Conc \ 1^{\ColB{aa}\ColA{bbb}} 1^{\ColB{aa}} 2^{\ColB{aa}\ColA{b}}
    \ \Conc \ \PWordEmpty
    \ \Conc \ 2^{\ColA{aa}} 2^{\ColA{a}}
    \notag \\
    & = 1^{\ColB{ba}\ColA{bbb}} 1^{\ColB{ba}} 2^{\ColB{ba}\ColA{b}} 1^{\ColB{aa}\ColA{bbb}}
    1^{\ColB{aa}} 2^{\ColB{aa}\ColA{b}} 2^{\ColA{aa}} 2^{\ColA{a}}.
    \notag
\end{align}
We also set, for any $n \geq 1$ and $i \in [n]$, $\Unit_{i, n}$ as the pigmented word
$i^\MonoidUnit$ of length $1$. For instance, by considering the monoid $\Monoid$ of the
previous example, $\Unit_{2, 4}$ is the pigmented word $2^\epsilon$.

Besides, given a monoid morphism $\phi : \Monoid \to \Monoid'$ between two monoids $\Monoid$
and $\Monoid'$, let $\P(\phi) : \P(\Monoid) \to \P\Par{\Monoid'}$ be the map defined for any
$\Monoid$-pigmented word $i_1^{\alpha_1} \dots i_\ell^{\alpha_\ell}$ by
\begin{equation}
    \P(\phi)\Par{
    i_1^{\alpha_1} \dots i_\ell^{\alpha_\ell}}
    :=
    i_1^{\phi\Par{\alpha_1}} \dots i_\ell^{\phi\Par{\alpha_\ell}}.
\end{equation}
For instance, by denoting by $\N$ the additive monoid $(\N, +, 0)$, the map $\phi : A^* \to
\N$ sending each $w \in A^*$ to its length is a monoid morphism. We have in this context
\begin{equation}
    \P(\phi)\Par{2^{ba} 2^{aa} 3^\epsilon}
    = 2^2 2^2 3^0.
\end{equation}

\begin{Statement}{Theorem}{thm:construction_p}
    The construction $\P$ is a functor from the category of monoids to the category of
    clones. Moreover, this functor preserves injections and surjections.
\end{Statement}
\begin{Proof}
    In this proof, we consider two monoids $(\Monoid, \MonoidProduct, \MonoidUnit)$ and
    $\Par{\Monoid', \MonoidProduct', \MonoidUnit'}$.

    Let us first prove that $\P(\Monoid)$ is a clone. For any $n \geq 1$, $i \in [n]$, and
    $\PWord_1, \dots, \PWord_n \in \P(\Monoid)(m)$, $m \geq 0$, since $\MonoidUnit$ is the
    unit of $\Monoid$, we have
    \begin{math}
        i^\MonoidUnit \Superposition{\PWord_1, \dots, \PWord_n} = \PWord_i
    \end{math}
    so that Relation~\eqref{equ:clone_relation_1} is satisfied. Moreover, for any $n \geq 0$
    and $\PWord \in \P(\Monoid)$, again since $\MonoidUnit$ is the unit of $\Monoid$, we
    have
    \begin{math}
        \PWord \Superposition{1^\MonoidUnit, \dots, n^\MonoidUnit} = \PWord
    \end{math}
    so that Relation~\eqref{equ:clone_relation_2} is satisfied. Finally, for any $n \geq 0$,
    $m \geq 0$, $k \geq 0$, $i_1^{\alpha_1} \dots i_\ell^{\alpha_\ell} \in \P(\Monoid)(n)$,
    \begin{math}
        j_{1, 1}^{\beta_{1, 1}} \dots j_{1, k_1}^{\beta_{1, k_1}},
        \dots,
        j_{n, 1}^{\beta_{n, 1}} \dots j_{n, k_n}^{\beta_{n, k_n}}
        \in \P(\Monoid)(m),
    \end{math}
    and $\PWord_1, \dots, \PWord_m \in \P(\Monoid)(k)$, since $\MonoidProduct$ is
    associative, we have
    \begin{align}
        i_1^{\alpha_1} \dots & i_\ell^{\alpha_\ell}
        \Superposition{
            j_{1, 1}^{\beta_{1, 1}} \dots j_{1, k_1}^{\beta_{1, k_1}}, \dots,
            j_{n, 1}^{\beta_{n, 1}} \dots j_{n, k_n}^{\beta_{n, k_n}}
        }
        \Superposition{\PWord_1, \dots, \PWord_m}
        \\
        & =
        \Par{\alpha_1 \MonoidProductExtension
            j_{i_1, 1}^{\beta_{i_1, 1}} \dots j_{i_1, k_{i_1}}^{\beta_{i_1, k_{i_1}}}}
        \dots
        \Par{\alpha_\ell \MonoidProductExtension
            j_{i_\ell, 1}^{\beta_{i_\ell, 1}} \dots
                j_{i_\ell, k_{i_\ell}}^{\beta_{i_\ell, k_{i_\ell}}}}
        \Superposition{\PWord_1, \dots, \PWord_m}
        \notag \\
        & =
        j_{i_1, 1}^{\alpha_1 \MonoidProduct \beta_{i_1, 1}} \dots
            j_{i_1, k_{i_1}}^{\alpha_1 \MonoidProduct \beta_{i_1, k_{i_1}}}
        \dots
        j_{i_\ell, 1}^{\alpha_\ell \MonoidProduct \beta_{i_\ell, 1}} \dots
            j_{i_\ell, k_{i_\ell}}^{\alpha_\ell \MonoidProduct \beta_{i_\ell, k_{i_\ell}}}
        \Superposition{\PWord_1, \dots, \PWord_m}
        \notag \\
        & =
        \Par{\Par{\alpha_1 \MonoidProduct \beta_{i_1, 1}} \MonoidProductExtension
            \PWord_{j_{i_1, 1}}}
        \dots
        \Par{\Par{\alpha_1 \MonoidProduct \beta_{i_1, k_{i_1}}} \MonoidProductExtension
            \PWord_{j_{i_1, {k_{i_1}}}}}
        \notag \\
        & \qquad \dots
        \Par{\Par{\alpha_\ell \MonoidProduct \beta_{i_\ell, 1}} \MonoidProductExtension
            \PWord_{j_{i_\ell, 1}}}
        \dots
        \Par{\Par{\alpha_\ell \MonoidProduct \beta_{i_\ell, k_{i_\ell}}}
            \MonoidProductExtension
            \PWord_{j_{i_\ell, {k_{i_\ell}}}}}
        \notag \\
        & =
        \alpha_1 \MonoidProductExtension \Par{
            \Par{\beta_{i_1, 1} \MonoidProductExtension \PWord_{j_{i_1, 1}}}
            \dots
            \Par{\beta_{i_1, k_{i_1}} \MonoidProductExtension \PWord_{j_{i_1, {k_{i_1}}}}}}
        \dots
        \alpha_\ell \MonoidProductExtension \Par{
            \Par{\beta_{i_\ell, 1} \MonoidProductExtension \PWord_{j_{i_\ell, 1}}}
            \dots
            \Par{\beta_{i_\ell, k_{i_\ell}} \MonoidProductExtension
                \PWord_{j_{i_\ell, {k_{i_\ell}}}}}}
        \notag \\
        & =
        i_1^{\alpha_1} \dots i_\ell^{\alpha_\ell} \Superposition{
            \Par{\beta_{1, 1} \MonoidProductExtension \PWord_{j_{1, 1}}}
            \dots
            \Par{\beta_{1, k_1} \MonoidProductExtension \PWord_{j_{1, k_1}}},
            \dots,
            \Par{\beta_{n, 1} \MonoidProductExtension \PWord_{j_{n, 1}}}
            \dots
            \Par{\beta_{n, k_n} \MonoidProductExtension \PWord_{j_{n, k_n}}}}
        \notag \\
        & =
        i_1^{\alpha_1} \dots i_\ell^{\alpha_\ell} \Superposition{
            j_{1, 1}^{\beta_{1, 1}} \dots j_{1, k_1}^{\beta_{1, k_1}}
                \Superposition{\PWord_1, \dots, \PWord_m},
            \dots,
            j_{n, 1}^{\beta_{n, 1}} \dots j_{n, k_n}^{\beta_{n, k_n}}
                \Superposition{\PWord_1, \dots, \PWord_m}}
        \notag
    \end{align}
    so that Relation~\eqref{equ:clone_relation_3} is satisfied. Therefore, $\P(\Monoid)$ is
    a clone.

    Let $\phi : \Monoid \to \Monoid'$ be a monoid morphism. Let us show that $\P(\phi)$ is a
    clone morphism. First, $\P(\phi)$ is a graded set morphism. Moreover, for any $n \geq 1$
    and $i \in [n]$, since $\phi$ sends the unit of $\Monoid$ to the unit of $\Monoid'$, we
    have
    \begin{math}
        \P(\phi)\Par{i^\MonoidUnit} = i^{\phi(\MonoidUnit)} = i^{\MonoidUnit'}.
    \end{math}
    Finally, for any $n \geq 0$, $m \geq 0$, $i_1^{\alpha_1} \dots i_\ell^{\alpha_\ell} \in
    \P(\Monoid)(n)$, and
    \begin{math}
        j_{1, 1}^{\beta_{1, 1}} \dots j_{1, k_1}^{\beta_{1, k_1}},
        \dots,
        j_{n, 1}^{\beta_{n, 1}} \dots j_{n, k_n}^{\beta_{n, k_n}}
        \in \P(\Monoid)(m),
    \end{math}
    since $\phi$ is a monoid morphism, we have
    \begin{align}
        \P(\phi) & \Par{
        i_1^{\alpha_1} \dots i_\ell^{\alpha_\ell}
        \Superposition{
            j_{1, 1}^{\beta_{1, 1}} \dots j_{1, k_1}^{\beta_{1, k_1}}, \dots,
            j_{n, 1}^{\beta_{n, 1}} \dots j_{n, k_n}^{\beta_{n, k_n}}
        }}
        \\
        & =
        \P(\phi) \Par{
        \Par{\alpha_1 \MonoidProductExtension
            j_{i_1, 1}^{\beta_{i_1, 1}} \dots j_{i_1, k_{i_1}}^{\beta_{i_1, k_{i_1}}}}
        \dots
        \Par{\alpha_\ell \MonoidProductExtension
            j_{i_\ell, 1}^{\beta_{i_\ell, 1}} \dots
                j_{i_\ell, k_{i_\ell}}^{\beta_{i_\ell, k_{i_\ell}}}}}
        \notag \\
        & =
        \P(\phi) \Par{
        j_{i_1, 1}^{\alpha_1 \MonoidProduct \beta_{i_1, 1}} \dots
            j_{i_1, k_{i_1}}^{\alpha_1 \MonoidProduct \beta_{i_1, k_{i_1}}}
        \dots
        j_{i_\ell, 1}^{\alpha_\ell \MonoidProduct \beta_{i_\ell, 1}} \dots
            j_{i_\ell, k_{i_\ell}}^{\alpha_\ell \MonoidProduct \beta_{i_\ell, k_{i_\ell}}}}
        \notag \\
        & =
        j_{i_1, 1}^{\phi\Par{\alpha_1 \MonoidProduct \beta_{i_1, 1}}} \dots
            j_{i_1, k_{i_1}}^{\phi\Par{\alpha_1 \MonoidProduct \beta_{i_1, k_{i_1}}}}
        \dots
        j_{i_\ell, 1}^{\phi\Par{\alpha_\ell \MonoidProduct \beta_{i_\ell, 1}}} \dots
            j_{i_\ell, k_{i_\ell}}^{\phi\Par{\alpha_\ell \MonoidProduct \beta_{i_\ell,
                k_{i_\ell}}}}
        \notag \\
        & =
        j_{i_1, 1}^{\phi\Par{\alpha_1} \MonoidProduct' \phi\Par{\beta_{i_1, 1}}} \dots
            j_{i_1, k_{i_1}}^{\phi\Par{\alpha_1} \MonoidProduct'
                \phi\Par{\beta_{i_1, k_{i_1}}}}
        \dots
        j_{i_\ell, 1}^{\phi\Par{\alpha_\ell} \MonoidProduct' \phi\Par{\beta_{i_\ell, 1}}}
            \dots
            j_{i_\ell, k_{i_\ell}}^{\phi\Par{\alpha_\ell} \MonoidProduct'
                \phi\Par{\beta_{i_\ell, k_{i_\ell}}}}
        \notag \\
        & =
        \Par{\phi\Par{\alpha_1} \MonoidProductExtension'
            j_{i_1, 1}^{\phi\Par{\beta_{i_1, 1}}} \dots
                j_{i_1, k_{i_1}}^{\phi\Par{\beta_{i_1, k_{i_1}}}}}
        \dots
        \Par{\phi\Par{\alpha_\ell} \MonoidProductExtension'
            j_{i_\ell, 1}^{\phi\Par{\beta_{i_\ell, 1}}} \dots
                j_{i_\ell, k_{i_\ell}}^{\phi\Par{\beta_{i_\ell, k_{i_\ell}}}}}
        \notag \\
        & =
        i_1^{\phi\Par{\alpha_1}} \dots i_\ell^{\phi\Par{\alpha_\ell}}
        \Superposition{
            j_{1, 1}^{\phi\Par{\beta_{1, 1}}} \dots j_{1, k_1}^{\phi\Par{\beta_{1, k_1}}},
            \dots,
            j_{n, 1}^{\phi\Par{\beta_{n, 1}}} \dots j_{n, k_n}^{\phi\Par{\beta_{n, k_n}}}
        }
        \notag \\
        & =
        \P(\phi) \Par{
        i_1^{\alpha_1} \dots i_\ell^{\alpha_\ell}}
        \Superposition{
            \P(\phi)\Par{j_{1, 1}^{\beta_{1, 1}} \dots j_{1, k_1}^{\beta_{1, k_1}}}, \dots,
            \P(\phi)\Par{j_{n, 1}^{\beta_{n, 1}} \dots j_{n, k_n}^{\beta_{n, k_n}}}
        }.
        \notag
    \end{align}
    Therefore, $\P(\phi)$ is a clone morphism. Moreover, it is immediate, for any monoid
    $\Monoid''$ and monoid morphism $\phi' : \Monoid' \to \Monoid''$, that $\P\Par{\phi'
    \circ \phi} = \P\Par{\phi'} \circ \P(\phi)$. It is also immediate that if $\IdentityMap
    : \Monoid \to \Monoid$ is the identity map, then $\P(\IdentityMap)$ is the identity map
    on $\P(\Monoid)$. For these reasons, $\P$ is a functor from the category of monoids to
    the category of clones.

    Let us finally prove that $\P$ preserves injections and surjections.
    Assume that $\phi$ is injective. If $i_1^{\alpha_1} \dots i_\ell^{\alpha_\ell}$
    and $j_1^{\beta_1} \dots j_k^{\beta_k}$ are two elements of $\P(\Monoid)$ such that
    \begin{math}
        \P(\phi)\Par{i_1^{\alpha_1} \dots i_\ell^{\alpha_\ell}}
        = \P(\phi)\Par{j_1^{\beta_1} \dots j_k^{\beta_k}},
    \end{math}
    then
    \begin{math}
        i_1^{\phi\Par{\alpha_1}} \dots i_\ell^{\phi\Par{\alpha_\ell}}
        = j_1^{\phi\Par{\beta_1}} \dots j_k^{\phi\Par{\beta_k}}.
    \end{math}
    Thus, $\ell = k$, $i_1 = j_1$, \dots, $i_\ell = j_\ell$, $\phi\Par{\alpha_1} =
    \phi\Par{\beta_1}$, \dots, $\phi\Par{\alpha_\ell} = \phi\Par{\beta_\ell}$. Since $\phi$
    is injective, we have $\alpha_1 = \beta_1$, \dots, $\alpha_\ell = \beta_\ell$, showing
    that $\P(\phi)$ is injective. Assume that $\phi$ is surjective. Let $j_1^{\beta_1} \dots
    j_k^{\beta_k} \in \P\Par{\Monoid'}$. Since $\phi$ is surjective, there are $\alpha_1,
    \dots, \alpha_k \in \Monoid$ such that $\phi\Par{\alpha_1} = \beta_1$, \dots,
    $\phi\Par{\alpha_k} = \beta_k$. Therefore, we have
    \begin{math}
        \P(\phi)\Par{j_1^{\alpha_1} \dots j_k^{\alpha_k}}
        = j_1^{\beta_1} \dots j_k^{\beta_k},
    \end{math}
    showing that $\P(\phi)$ is surjective.
\end{Proof}

\subsubsection{First properties} \label{subsubsec:first_properties}
We describe now a generating set of $\P(\Monoid)$ and show that the map sending any
$\Monoid$-pigmented word to its mirror image is an involutive clone automorphism
of~$\P(\Monoid)$.

\begin{Statement}{Proposition}{prop:generating_set}
    For any monoid $\Monoid$, the graded set
    \begin{math}
        G_\Monoid :=
        G_\Monoid(0)
        \sqcup G_\Monoid(1)
        \sqcup G_\Monoid(2)
    \end{math}
    defined by $G_\Monoid(0) := \Bra{\PWordEmpty}$, $G_\Monoid(1) := \Bra{1^\alpha : \alpha
    \in \Monoid}$, and $G_\Monoid(2) := \Bra{1^\MonoidUnit 2^\MonoidUnit}$ where
    $\MonoidUnit$ is the unit of $\Monoid$ is a generating set of the clone $\P(\Monoid)$.
\end{Statement}
\begin{Proof}
    Let us prove by induction on the length $\ell$ of $\PWord \in \P(\Monoid)(n)$, $n \geq
    0$, that $\PWord \in \P(\Monoid)^{G_\Monoid}$. First, if $\ell = 0$, then $\PWord =
    \PWordEmpty$ and since $\PWordEmpty \in G_\Monoid$, the property holds. If $\ell \geq
    1$, then $\PWord$ decomposes as $\PWord = \PWord' \Conc i^\alpha$ where $\PWord' \in
    \P(\Monoid)(n)$ and $i^\alpha \in \SetLetters_\Monoid$. By definition of the
    superposition maps of $\P(\Monoid)$, $\PWord$ can be expressed as
    \begin{math}
        \PWord =
        1^\MonoidUnit 2^\MonoidUnit
        \Superposition{\PWord', 1^\alpha \Superposition{\Unit_{i, n}}}.
    \end{math}
    Now, since $\Length\Par{\PWord'} = \ell - 1$, by induction hypothesis, $\PWord' \in
    \P(\Monoid)^{G_\Monoid}$. Moreover, since $1^\MonoidUnit 2^\MonoidUnit \in G_\Monoid$
    and $1^\alpha \in G_\Monoid$, this shows the previously stated property.
\end{Proof}

By considering the graded set $G_\Monoid$ introduced by
Proposition~\ref{prop:generating_set}, let $\Interpretation_\Monoid : \GeneratingSet_\Monoid
\to G_\Monoid$ be the graded set morphism defined by $\Interpretation_\Monoid(\Identity) :=
\PWordEmpty$, $\Interpretation_\Monoid\Par{\Pigmentation_\alpha} := 1^\alpha$, $\alpha \in
\Monoid$, and $\Interpretation_\Monoid(\Product) := 1^\MonoidUnit 2^\MonoidUnit$. This
bijective map will be used together with Proposition~\ref{prop:clone_presentation} in order
to establish a presentation of $\P(\Monoid)$.

The map $\Reverse$ sending any word to its mirror image is in particular a well-defined
graded set morphism from $\P(\Monoid)$ to $\P(\Monoid)$. As stated by the following result,
this map has an additional property.

\begin{Statement}{Proposition}{prop:reverse_map_clone_isomorphism}
    For any monoid $\Monoid$, the map $\Reverse : \P(\Monoid) \to \P(\Monoid)$ is an
    involutive clone automorphism.
\end{Statement}
\begin{Proof}
    Let $\MonoidProduct$ be the operation of $\Monoid$ and $\MonoidUnit$ its unit. It is
    first immediate that the projections $i^\MonoidUnit \in \P(\Monoid)(n)$, $n \geq 1$, $i
    \in [n]$, are fixed-points of $\Reverse$. Moreover, as a consequence of the fact that
    for any words $u$ and $v$ on any alphabet, $\Reverse(u \Conc v) = \Reverse(v) \Conc
        \Reverse(u)$, for any $i_1^{\alpha_1} \dots i_\ell^{\alpha_\ell} \in \P(\Monoid)$,
        $n \geq 0$, and $\PWord_1, \dots, \PWord_n \in \P(\Monoid)(m)$, $m \geq 0$, we have
    \begin{equation}
        \Reverse\Par{i_1^{\alpha_1} \dots i_\ell^{\alpha_\ell}
            \Superposition{\PWord_1, \dots, \PWord_n}}
        = \Reverse\Par{\alpha_\ell \MonoidProductExtension \PWord_{i_\ell}}
        \dots
        \Reverse\Par{\alpha_1 \MonoidProductExtension \PWord_{i_1}}
        = \Reverse\Par{i_1^{\alpha_1} \dots i_\ell^{\alpha_\ell}}
            \Superposition{\Reverse\Par{\PWord_1}, \dots, \Reverse\Par{\PWord_n}}.
    \end{equation}
    Therefore, $\Reverse$ is a clone morphism. Finally, since $\Reverse$ is an
    involution, the statement of the proposition follows.
\end{Proof}

\subsection{Clone realization} \label{subsec:clone_realization}
This last part of the present section is devoted to establish its main result, namely the
fact that $\P(\Monoid)$ is a clone realization of the variety of $\Monoid$-pigmented
monoids. For this, we shall use a method consisting in building a specific system of
representatives for the quotient $\SetTerms\Par{\GeneratingSet_\Monoid}
/_{\Equiv_{\RelationSet_\Monoid}}$ which is in one-to-one correspondence with the graded set
of $\Monoid$-pigmented words. Other approaches are possible as well including those using
term rewrite systems~\cite{BN98,BKV03} and proofs for their termination and confluence.

\subsubsection{Properties of the equation set} \label{subsubsec:properties_equation_set}
We begin with two elementary properties satisfied by the equivalence relation
$\Equiv_{\RelationSet_\Monoid}$.

\begin{Statement}{Lemma}{lem:variable_count_equivalence}
    For any monoid $\Monoid$ and any $\TreeT, \TreeT' \in
    \SetTerms\Par{\GeneratingSet_\Monoid}$, $\TreeT \Equiv_{\RelationSet_\Monoid} \TreeT'$
    implies that $\TreeT$ and $\TreeT'$ have the same variable count.
\end{Statement}
\begin{Proof}
    For any equation $\Par{\TreeT, \TreeT'}$ of the variety $\Par{\GeneratingSet_\Monoid,
    \RelationSet_\Monoid}$ (see
    Relations~\eqref{equ:pigmented_monoids_1}---\eqref{equ:pigmented_monoids_6}), we can
    observe that $\VariableCount(\TreeT) = \VariableCount\Par{\TreeT'}$ and that each
    variable appears at most once in $\TreeT$ and $\TreeT'$. Since by definition,
    $\Equiv_{\RelationSet_\Monoid}$ is the smallest clone congruence containing
    $\RelationSet_\Monoid$, the statement of the lemma follows.
\end{Proof}

The \Def{frontier map} is the map $\Frontier_\Monoid : \SetTerms\Par{\GeneratingSet_\Monoid}
\to \P(\Monoid)$ defined by $\Frontier_\Monoid := \Eval_{\P(\Monoid)} \circ
\widehat{\Interpretation_\Monoid}$, where the graded set morphism $\Interpretation_\Monoid :
\GeneratingSet_\Monoid \to \P(\Monoid)$ is as defined in
Section~\ref{subsubsec:first_properties}. Clearly, $\Frontier_\Monoid$ is a clone morphism
since $\Eval_{\P(\Monoid)}$ and $\widehat{\Interpretation_\Monoid}$ are. For instance, by
considering the free monoid $\Par{A^*, \Conc, \epsilon}$ generated by $A := \Bra{a, b}$, we
have in $\P\Par{A^*}$,
\begin{align}
    \Frontier_{\P\Par{A^*}}
    & \Par{
    \Product \Superposition{
        \Pigmentation_a\Superposition{
            \Product\Superposition{\VarX_3, \Pigmentation_b\Superposition{\VarX_2}}},
        \Product\Superposition{\VarX_1, \Pigmentation_b\Superposition{\VarX_2}}}}
    \\
    & =
    \Eval_{\P\Par{A^*}}\Par{\widehat{\Interpretation_{A^*}}\Par{
    \Product \Superposition{
        \Pigmentation_a\Superposition{
            \Product\Superposition{\VarX_3, \Pigmentation_b\Superposition{\VarX_2}}},
        \Product\Superposition{\VarX_1, \Pigmentation_b\Superposition{\VarX_2}}}}}
    \notag \\
    & =
    \Eval_{\P\Par{A^*}}\Par{
    1^\epsilon 2^\epsilon \Superposition{
        1^a\Superposition{1^\epsilon 2^\epsilon\Superposition{
            \VarX_3, 1^b\Superposition{\VarX_2}}},
        1^\epsilon 2^\epsilon \Superposition{\VarX_1, 1^b\Superposition{\VarX_2}}}}
    \notag \\
    & = 3^a 2^{ab} 1^\epsilon 2^b.
    \notag
\end{align}

\begin{Statement}{Lemma}{lem:equivalence_implies_same_evaluations}
    For any monoid $\Monoid$ and any $\TreeT, \TreeT' \in
    \SetTerms\Par{\GeneratingSet_\Monoid}$, $\TreeT \Equiv_{\RelationSet_\Monoid} \TreeT'$
    implies $\Frontier_\Monoid(\TreeT) = \Frontier_\Monoid\Par{\TreeT'}$.
\end{Statement}
\begin{Proof}
    Let $\MonoidProduct$ be the operation of $\Monoid$ and $\MonoidUnit$ is its unit. For
    any $\alpha, \alpha_1, \alpha_2 \in \Monoid$, we have
    \begin{equation}
        \Frontier_\Monoid\Par{\Product\Superposition{
            \Product\Superposition{\VarX_1, \VarX_2}, \VarX_3}}
        = 1^\MonoidUnit 2^\MonoidUnit 3^\MonoidUnit
        = \Frontier_\Monoid\Par{\Product\Superposition{
            \VarX_1, \Product\Superposition{\VarX_2, \VarX_3}}},
    \end{equation}
    \begin{equation}
        \Frontier_\Monoid\Par{\Product\Superposition{\Identity, \VarX_1}}
        = 1^\MonoidUnit
        = \Frontier_\Monoid\Par{\VarX_1}
        = 1^\MonoidUnit
        = \Frontier_\Monoid\Par{\Product\Superposition{\VarX_1, \Identity}},
    \end{equation}
    \begin{equation}
        \Frontier_\Monoid\Par{\Pigmentation_\alpha\Superposition{
            \Product\Superposition{\VarX_1, \VarX_2}}}
        = 1^\alpha 2^\alpha
        = \Frontier_\Monoid\Par{\Product\Superposition{
            \Pigmentation_\alpha\Superposition{\VarX_1},
            \Pigmentation_\alpha\Superposition{\VarX_2}}},
    \end{equation}
    \begin{equation}
        \Frontier_\Monoid\Par{\Pigmentation_\alpha \Superposition{\Identity}}
        = \PWordEmpty
        = \Frontier_\Monoid\Par{\Identity},
    \end{equation}
    \begin{equation}
        \Frontier_\Monoid\Par{\Pigmentation_{\alpha_1}
            \Superposition{\Pigmentation_{\alpha_2}\Superposition{\VarX_1}}}
        = 1^{\alpha_1 \MonoidProduct \alpha_2}
        = \Frontier_\Monoid\Par{\Par{\Pigmentation_{\alpha_1 \MonoidProduct \alpha_2}}
            \Superposition{\VarX_1}},
    \end{equation}
    \begin{equation}
        \Frontier_\Monoid\Par{\Pigmentation_\MonoidUnit\Superposition{\VarX_1}}
        = 1^\MonoidUnit
        = \Frontier_\Monoid\Par{\VarX_1}.
    \end{equation}
    Since by definition, $\Equiv_{\RelationSet_\Monoid}$ is the smallest clone congruence
    containing $\RelationSet_\Monoid$ and, as we have seen here, for any $\Par{\TreeT,
    \TreeT'} \in \RelationSet_{\Monoid}$, we have $\Frontier_\Monoid(\TreeT) =
    \Frontier_\Monoid\Par{\TreeT'}$, and because $\Frontier_\Monoid$ is a clone morphism
    (that is, its kernel is a clone congruence), the statement of the lemma follows.
\end{Proof}

\subsubsection{Right comb factorization} \label{subsubsec:right_comb_factorization}
We describe now a way to encode any $\Monoid$-pigmented word as a particular
$\GeneratingSet_\Monoid$-term having some important properties.

The \Def{right comb factorization map} is the map $\RightComb_\Monoid : \P(\Monoid) \to
\SetTerms\Par{\GeneratingSet_\Monoid}$ recursively defined, for any $\PWord \in
\P(\Monoid)$, by
\begin{equation}
    \RightComb_\Monoid\Par{\PWord}
    :=
    \begin{cases}
        \Identity & \text{if } \PWord = \PWordEmpty, \\
        \Product \Superposition{
            \Pigmentation_\alpha \Superposition{\VarX_i}, \RightComb_\Monoid\Par{\PWord'}}
        & \text{otherwise, where } \PWord = i^\alpha \Conc \PWord',
    \end{cases}
\end{equation}
where $i^\alpha \in \SetLetters_\Monoid$, $\VarX_i \in \SetVariables$, and $\PWord' \in
\P(\Monoid)$. For instance, for the free monoid $\Par{A^*, \Conc, \epsilon}$ where $A$ is
the alphabet $\Bra{a, b, c}$, we have
\begin{align}
    \RightComb_{A^*}\Par{1^{ab} 3^{aa} 2^\epsilon 2^b}
    & =
    \Product \Superposition{\Pigmentation_{ab} \Superposition{\VarX_1},
    \Product \Superposition{\Pigmentation_{aa} \Superposition{\VarX_3},
    \Product \Superposition{\Pigmentation_\epsilon \Superposition{\VarX_2},
    \Product \Superposition{\Pigmentation_b \Superposition{\VarX_2}, \Identity}}}}
    \\
    & =
    \Pigmentation_{ab} \Superposition{\VarX_1}
    \Product \Pigmentation_{aa} \Superposition{\VarX_3}
    \Product \Pigmentation_{\epsilon} \Superposition{\VarX_2}
    \Product \Pigmentation_{b} \Superposition{\VarX_2}
    \Product \Identity.
    \notag
\end{align}

\begin{Statement}{Lemma}{lem:evaluation_right_comb}
    For any monoid $\Monoid$ and any $\PWord \in \P(\Monoid)$,
    $\Frontier_\Monoid(\RightComb_\Monoid(\PWord)) = \PWord$.
\end{Statement}
\begin{Proof}
    Let $\MonoidProduct$ be the operation of $\Monoid$ and $\MonoidUnit$ is its unit. We
    proceed by induction on the length $\ell$ of $\PWord$. If $\ell = 0$, then $\PWord =
    \PWordEmpty$ and since
    \begin{math}
        \Frontier_\Monoid(\RightComb_\Monoid(\PWordEmpty))
        = \Frontier_\Monoid(\Identity)
        = \PWordEmpty,
    \end{math}
    the property holds. If $\ell \geq 1$, $\PWord$ decomposes as $\PWord = i^\alpha \Conc
    \PWord'$ where $i^\alpha \in \SetLetters_\Monoid$ and $\PWord' \in \P(\Monoid)$. By
    definition of $\RightComb_\Monoid$ and by induction hypothesis,
    \begin{align}
        \Frontier_\Monoid(\RightComb_\Monoid(\PWord))
        = \Frontier_\Monoid\Par{\RightComb_\Monoid\Par{i^\alpha \Conc \PWord'}}
        & = \Frontier_\Monoid\Par{\Product \Superposition{
            \Pigmentation_\alpha \Superposition{\VarX_i},
            \RightComb_\Monoid\Par{\PWord'}}}
        \\
        & = 1^\MonoidUnit 2^\MonoidUnit \Superposition{
            1^\alpha \Superposition{i^\MonoidUnit},
            \Frontier_\Monoid\Par{\RightComb_\Monoid\Par{\PWord'}}}
        = 1^\MonoidUnit 2^\MonoidUnit \Superposition{
            1^\alpha \Superposition{i^\MonoidUnit}, \PWord'}
        = i^\alpha \Conc \PWord'
        = \PWord.
        \notag
    \end{align}
    Therefore, the stated property holds.
\end{Proof}

As a consequence of Lemma~\ref{lem:evaluation_right_comb}, $\Frontier_\Monoid :
\SetTerms\Par{\GeneratingSet_\Monoid} \to \P(\Monoid)$ is a surjective clone morphism and
$\RightComb_\Monoid : \P(\Monoid) \to \SetTerms\Par{\GeneratingSet_\Monoid}$ is an injective
map.

\begin{Statement}{Lemma}{lem:equivalence_right_comb_evaluation}
    For any monoid $\Monoid$ and any $\TreeT \in \SetTerms\Par{\GeneratingSet_\Monoid}$,
    there exists $\TreeT' \in \RightComb_\Monoid(\P(\Monoid))$ such that $\TreeT
    \Equiv_{\RelationSet_\Monoid} \TreeT'$.
\end{Statement}
\begin{Proof}
    Let $\MonoidProduct$ be the operation of $\Monoid$ and $\MonoidUnit$ is its unit. We
    proceed by induction on the pairs $(\ell, d)$ ordered lexicographically, where $\ell$ is
    the variable count of $\TreeT$ and $d$ is the operation count of $\TreeT$.
    \begin{enumerate}[label={\it (\Roman*)}]
        \item If $\ell = 0$, then $\TreeT$ has no variable.
        By~\eqref{equ:pigmented_monoids_2} and~\eqref{equ:pigmented_monoids_4}, $\TreeT
        \Equiv_{\RelationSet_\Monoid} \Identity$. Since $\Identity$ belongs to
        $\RightComb_\Monoid(\P(\Monoid))$, the stated property is satisfied.
        \item If $\ell \geq 1$, we have three sub-cases to explore depending on the general
        form of $\TreeT$.
        \begin{enumerate}[label={\it (\alph*)}]
            \item If $\TreeT = \VarX_i$ where $\VarX_i \in \SetVariables$,
            by~\eqref{equ:pigmented_monoids_2},
            \begin{math}
                \TreeT \Equiv_{\RelationSet_\Monoid}
                \Product \Superposition{\VarX_i, \Identity}.
            \end{math}
            By~\eqref{equ:pigmented_monoids_6},
            \begin{math}
                \TreeT \Equiv_{\RelationSet_\Monoid}
                \Product \Superposition{\Pigmentation_\MonoidUnit \Par{\VarX_i}, \Identity}.
            \end{math}
            Since $\Product \Superposition{\Pigmentation_\MonoidUnit \Par{\VarX_i},
            \Identity}$ belongs to $\RightComb_\Monoid(\P(\Monoid))$, the stated property is
            satisfied.
            \item If $\TreeT = \Pigmentation_\alpha \Superposition{\TreeS}$ where $\alpha
            \in \Monoid$ and $\TreeS \in \SetTerms\Par{\GeneratingSet_\Monoid}$, since
            $\VariableCount(\TreeS) = \VariableCount(\TreeT)$ and $\OperationCount(\TreeS) <
            \OperationCount(\TreeT)$, by induction hypothesis, there exists $\TreeS' \in
            \RightComb_\Monoid(\P(\Monoid))$ such that $\TreeS \Equiv_{\RelationSet_\Monoid}
            \TreeS'$. By definition of $\RightComb_\Monoid$, $\TreeS'$ can have two
            different forms.
            \begin{enumerate}[label={\it (\roman*)}]
                \item If $\TreeS' = \Identity$, we have $\TreeT
                \Equiv_{\RelationSet_\Monoid} \Pigmentation_\alpha
                \Superposition{\Identity}$. By~\eqref{equ:pigmented_monoids_4},
                \begin{math}
                    \TreeT \Equiv_{\RelationSet_\Monoid} \Identity.
                \end{math}
                Since $\Identity$ belongs to $\RightComb_\Monoid(\P(\Monoid))$, the stated
                property is satisfied.
                \item Otherwise, $\TreeS' = \Product \Superposition{\Pigmentation_{\alpha'}
                \Superposition{\VarX_i}, \TreeR}$ where $\alpha' \in \Monoid$, $\VarX_i \in
                \SetVariables$, and $\TreeR \in \SetTerms\Par{\GeneratingSet_\Monoid}$. We
                have
                \begin{math}
                    \TreeT \Equiv_{\RelationSet_\Monoid}
                    \Pigmentation_\alpha \Superposition{\Product
                        \Superposition{\Pigmentation_{\alpha'}\Superposition{\VarX_i},
                        \TreeR}}.
                \end{math}
                By~\eqref{equ:pigmented_monoids_3}, we have
                \begin{math}
                    \TreeT \Equiv_{\RelationSet_\Monoid}
                        \Product \Superposition{
                        \Pigmentation_\alpha \Superposition{\Pigmentation_{\alpha'}
                        \Superposition{\VarX_i}},
                        \Pigmentation_\alpha \Superposition{\TreeR}}
                \end{math}
                and by~\eqref{equ:pigmented_monoids_5}, we have
                \begin{math}
                    \TreeT \Equiv_{\RelationSet_\Monoid}
                    \Product \Superposition{
                        \Pigmentation_{\alpha \MonoidProduct \alpha'}
                        \Superposition{\VarX_i},
                        \Pigmentation_\alpha \Superposition{\TreeR}}.
                \end{math}
                Now, by Lemma~\ref{lem:variable_count_equivalence},
                $\VariableCount\Par{\Pigmentation_\alpha \Superposition{\TreeR}} <
                \VariableCount(\TreeT)$. Thus, by induction hypothesis, there exists
                $\TreeR' \in \RightComb_\Monoid(\P(\Monoid))$ such that
                $\Pigmentation_\alpha \Superposition{\TreeR} \Equiv_{\RelationSet_\Monoid}
                \TreeR'$. Therefore,
                \begin{math}
                    \TreeT \Equiv_{\RelationSet_\Monoid}
                    \Product \Superposition{
                        \Pigmentation_{\alpha \MonoidProduct \alpha'}
                        \Superposition{\VarX_i},
                        \TreeR'}.
                \end{math}
                By definition of $\RightComb_\Monoid$,
                \begin{math}
                    \Product \Superposition{
                        \Pigmentation_{\alpha \MonoidProduct \alpha'}
                        \Superposition{\VarX_i}, \TreeR'}
                \end{math}
                belongs to $\RightComb_\Monoid(\P(\Monoid))$ so that the stated property is
                satisfied.
            \end{enumerate}
            \item Otherwise, $\TreeT = \Product \Superposition{\TreeS_1, \TreeS_2}$ where
            $\TreeS_1 \in \SetTerms\Par{\GeneratingSet_\Monoid}$ and $\TreeS_2 \in
            \SetTerms\Par{\GeneratingSet_\Monoid}$. Since $\VariableCount\Par{\TreeS_1} \leq
            \VariableCount(\TreeT)$, $\OperationCount\Par{\TreeS_1} <
            \OperationCount(\TreeT)$, $\VariableCount\Par{\TreeS_2} \leq
            \VariableCount(\TreeT)$, and $\OperationCount\Par{\TreeS_2} <
            \OperationCount(\TreeT)$, by induction hypothesis, there exist $\TreeS'_1,
            \TreeS'_2 \in \RightComb_\Monoid(\P(\Monoid))$ such that $\TreeS_1
            \Equiv_{\RelationSet_\Monoid} \TreeS'_1$ and $\TreeS_2
            \Equiv_{\RelationSet_\Monoid} \TreeS'_2$. By definition of $\RightComb_\Monoid$,
            $\TreeS'_1$ and $\TreeS'_2$ decompose respectively as
            \begin{math}
                \TreeS'_1 =
                \Pigmentation_{\alpha_{1, 1}} \Superposition{\VarX_{i_{1, 1}}}
                \Product \dots \Product
                \Pigmentation_{\alpha_{1, k_1}} \Superposition{\VarX_{i_{1, k_1}}}
                \Product \Identity
            \end{math}
            and
            \begin{math}
                \TreeS'_2 =
                \Pigmentation_{\alpha_{2, 1}} \Superposition{\VarX_{i_{2, 1}}}
                \Product \dots \Product
                \Pigmentation_{\alpha_{2, k_2}} \Superposition{\VarX_{i_{2, k_2}}}
                \Product \Identity
            \end{math}
            for some $\alpha_{1, 1}, \dots, \alpha_{1, k_1}, \alpha_{2, 1}, \dots,
            \alpha_{2, k_2} \in \Monoid$, $\VarX_{i_{1, 1}}, \dots, \VarX_{i_{1, k_1}},
            \VarX_{i_{2, 1}}, \dots, \VarX_{i_{2, k_2}} \in \SetVariables$, $k_1 \geq 0$,
            and $k_2 \geq 0$. Now, by~\eqref{equ:pigmented_monoids_1}
            and~\eqref{equ:pigmented_monoids_2}, we have
            \begin{align} \label{equ:equivalence_right_comb_evaluation_1}
                \TreeT
                & \Equiv_{\RelationSet_\Monoid}
                \TreeS'_1 \Product \TreeS'_2
                \\
                & =
                \Par{
                \Pigmentation_{\alpha_{1, 1}} \Superposition{\VarX_{i_{1, 1}}}
                \Product \dots \Product
                \Pigmentation_{\alpha_{1, k_1}} \Superposition{\VarX_{i_{1, k_1}}}
                \Product \Identity
                }
                \Product
                \Par{
                \Pigmentation_{\alpha_{2, 1}} \Superposition{\VarX_{i_{2, 1}}}
                \Product \dots \Product
                \Pigmentation_{\alpha_{2, k_2}} \Superposition{\VarX_{i_{2, k_2}}}
                \Product \Identity
                }
                \notag \\
                & \Equiv_{\RelationSet_\Monoid}
                \Par{
                \Pigmentation_{\alpha_{1, 1}} \Superposition{\VarX_{i_{1, 1}}}
                \Product \dots \Product
                \Pigmentation_{\alpha_{1, k_1}} \Superposition{\VarX_{i_{1, k_1}}}
                }
                \Product
                \Par{
                \Pigmentation_{\alpha_{2, 1}} \Superposition{\VarX_{i_{2, 1}}}
                \Product \dots \Product
                \Pigmentation_{\alpha_{2, k_2}} \Superposition{\VarX_{i_{2, k_2}}}
                \Product \Identity
                }
                \notag \\
                & \Equiv_{\RelationSet_\Monoid}
                \Pigmentation_{\alpha_{1, 1}} \Superposition{\VarX_{i_{1, 1}}}
                \Product \dots \Product
                \Pigmentation_{\alpha_{1, k_1}} \Superposition{\VarX_{i_{1, k_1}}}
                \Product
                \Pigmentation_{\alpha_{2, 1}} \Superposition{\VarX_{i_{2, 1}}}
                \Product \dots \Product
                \Pigmentation_{\alpha_{2, k_2}} \Superposition{\VarX_{i_{2, k_2}}}
                \Product \Identity
                \notag.
            \end{align}
            By definition of $\RightComb_\Monoid$, the last term
            of~\eqref{equ:equivalence_right_comb_evaluation_1} belongs to
            $\RightComb_\Monoid(\P(\Monoid))$ so that the stated property is satisfied.
        \end{enumerate}
    \end{enumerate}
\end{Proof}

\subsubsection{Clone presentation}
We use now the tools developed in the previous sections to prove that $\P(\Monoid)$ is a
clone realization of the variety of $\Monoid$-pigmented monoids.

\begin{Statement}{Lemma}{lem:same_evaluations_implies_equivalence}
    For any monoid $\Monoid$ and any $\TreeT, \TreeT' \in
    \SetTerms\Par{\GeneratingSet_\Monoid}$, $\Frontier_\Monoid(\TreeT) =
    \Frontier_\Monoid\Par{\TreeT'}$ implies $\TreeT \Equiv_{\RelationSet_\Monoid} \TreeT'$.
\end{Statement}
\begin{Proof}
    Assume that $\Frontier_\Monoid(\TreeT) = \Frontier_\Monoid\Par{\TreeT'}$. By
    Lemma~\ref{lem:equivalence_right_comb_evaluation}, there exist $\PWord, \PWord' \in
    \P(\Monoid)$ such that $\TreeT \Equiv_{\RelationSet_\Monoid} \RightComb_\Monoid(\PWord)$
    and $\TreeT' \Equiv_{\RelationSet_\Monoid} \RightComb_\Monoid\Par{\PWord'}$. By
    Lemma~\ref{lem:equivalence_implies_same_evaluations}, $\Frontier_\Monoid(\TreeT) =
    \Frontier_\Monoid(\RightComb_\Monoid(\PWord))$ and $\Frontier_\Monoid\Par{\TreeT'} =
    \Frontier_\Monoid(\RightComb_\Monoid\Par{\PWord'})$. By
    Lemma~\ref{lem:evaluation_right_comb}, $\Frontier_\Monoid(\TreeT) = \PWord$ and
    $\Frontier_\Monoid\Par{\TreeT'} = \PWord'$. Since $\Frontier_\Monoid(\TreeT) =
    \Frontier_\Monoid\Par{\TreeT'}$, we have $\PWord = \PWord'$. This shows that
    \begin{math}
        \TreeT
        \Equiv_{\RelationSet_\Monoid} \RightComb_\Monoid(\PWord)
        = \RightComb_\Monoid\Par{\PWord'}
        \Equiv_{\RelationSet_\Monoid} \TreeT',
    \end{math}
    so that $\TreeT \Equiv_{\RelationSet_\Monoid} \TreeT'$.
\end{Proof}

Here is the main result of the section.

\begin{Statement}{Theorem}{thm:clone_presentation_pigmented_monoids}
    For any monoid $\Monoid$, the clone $\P(\Monoid)$ is a clone realization of the variety
    of $\Monoid$-pigmented monoids.
\end{Statement}
\begin{Proof}
    By Lemmas~\ref{lem:same_evaluations_implies_equivalence}
    and~\ref{lem:equivalence_implies_same_evaluations}, for any $\TreeT, \TreeT' \in
    \SetTerms\Par{\GeneratingSet_\Monoid}$, $\TreeT \Equiv_{\RelationSet_\Monoid} \TreeT'$
    if and only if $\Frontier_\Monoid(\TreeT) = \Frontier_\Monoid\Par{\TreeT'}$. Moreover,
    by Proposition~\ref{prop:generating_set}, $G_\Monoid =
    \Interpretation_\Monoid\Par{\GeneratingSet_\Monoid}$ is a generating set of
    $\P(\Monoid)$. Therefore, by Proposition~\ref{prop:clone_presentation}, these two
    properties imply that the variety $\Par{\GeneratingSet_\Monoid, \RelationSet_\Monoid}$
    of $\Monoid$-pigmented monoids is a presentation of $\P(\Monoid)$.
\end{Proof}

By Theorem~\ref{thm:clone_presentation_pigmented_monoids}, for any monoid $\Monoid$, all
algebras over $\P(\Monoid)$ are $\Monoid$-pigmented monoids. Recall that $\T$ is a functor
from the category of monoids to the category of operads introduced in~\cite{Gir15}. Since
all algebras over the operad $\T(\Monoid)$ can be seen as specialized versions of
$\Monoid$-pigmented monoids, we can see the construction $\P$ as a generalization of the
construction $\T$ at the level of clones.

\subsubsection{Clone presentations of quotients} \label{subsubsec:remark_presentations}
Let us provide a remark, useful in the sequel when we study several quotients of
$\P(\Monoid)$. Let $\Equiv$ be a clone congruence of $\P(\Monoid)$ generated by a graded set
binary relation $\mathcal{R}$ on the set of $\Monoid$-pigmented words. By
Theorem~\ref{thm:clone_presentation_pigmented_monoids},
Proposition~\ref{prop:presentation_quotient}, and
Lemmas~\ref{lem:equivalence_implies_same_evaluations}
and~\ref{lem:same_evaluations_implies_equivalence}, the quotient $\P(\Monoid) /_{\Equiv}$
admits the presentation $\Par{\GeneratingSet_\Monoid, \RelationSet_\Monoid \cup
\RelationSet}$, where $\RelationSet$ is the graded set binary relation on
$\SetTerms\Par{\GeneratingSet_\Monoid}$ satisfying $\RightComb_\Monoid(\PWord) \
\RelationSet \ \RightComb_\Monoid\Par{\PWord'}$ whenever $\PWord \ \mathcal{R} \ \PWord'$.
By Lemmas \ref{lem:equivalence_implies_same_evaluations}, \ref{lem:evaluation_right_comb},
and~\ref{lem:equivalence_right_comb_evaluation}, $\RelationSet$ is also the graded set
binary relation where for any $\GeneratingSet_\Monoid$-terms $\TermT$ and $\TermT'$, $\TermT
\ \RelationSet \ \TermT'$ whenever $\Frontier_\Monoid(\TermT) \ \mathcal{R} \
\Frontier_\Monoid\Par{\TermT'}$, which is the preimage of $\mathcal{R}$ under the
isomorphism of Theorem~\ref{thm:clone_presentation_pigmented_monoids} presenting
$\P(\Monoid)$ as a quotient of $\SetTerms\Par{\GeneratingSet_\Monoid}$.

\section{Construction of quotients} \label{sec:construction_quotients}
The clones $\P(\Monoid)$ are very large and contain a lot of subclones and quotients worth
investigating. We present here some tools to construct quotients of $\P(\Monoid)$ through
so-called $\PSymbol$-symbols which are here particular maps from $\P(\Monoid)$ to itself.
Results about the description of the elements of such quotients are provided. As a direct
application, we construct in this section the quotients clones $\WInc(\Monoid)$,
$\Arra_k(\Monoid)$, and $\Inc_k$ of $\P(\Monoid)$.

\subsection{\texorpdfstring{$\PSymbol$}{P}-symbols and quotient clones}
\label{subsec:p_symbols}
A $\PSymbol$-symbol of $\P(\Monoid)$ is a map $\PSymbol : \P(\Monoid) \to \P(\Monoid)$
satisfying some properties, presented hereafter. Such maps will be used in this work to
build quotients of $\P(\Monoid)$ and describe explicitly their projections and superposition
maps.

As it is usually the case in the description of $\PSymbol$-symbols, it is always possible to
provide an iterative description of such maps through algorithms by setting
$\PSymbol(\PWordEmpty) := \PWordEmpty$ and by computing $\PSymbol\Par{\PWord \Conc
i^\alpha}$ as the insertion of the $\Monoid$-pigmented letter $i^\alpha$ into the
$\Monoid$-pigmented word $\PSymbol(\PWord)$. As a side remark, most $\PSymbol$-symbols
appearing in the literature map words to other combinatorial objects (like Young
tableaux~\cite[Chap.~5]{Lot02}, binary trees~\cite{HNT05}, or pairs of twin binary
trees~\cite{Gir12}). They lead to the constructions of various monoids as explicit quotients
of free monoids. Here, our notion of $\PSymbol$-symbol is very specific to our purposes.

\subsubsection{$\PSymbol$-symbols} \label{subsubsec:p_symbols}
Let $X$ be a set and $\Equiv$ be an equivalence relation on~$X$. A \Def{$\PSymbol$-symbol}
for $\Equiv$ is a map $\phi : X \to X$ such that
\begin{enumerate}[label=(\roman*)]
    \item \label{item:p_symbol_1}
    for any $x \in X$, $x \Equiv \phi(x)$;
    \item \label{item:p_symbol_2}
    for any $x, x' \in X$, $x \Equiv x'$ implies $\phi(x) = \phi\Par{x'}$.
\end{enumerate}
By extension, given $x \in X$, $\phi(x)$ is the \Def{$\PSymbol$-symbol} of $x$. Besides, for
any $x \in X$, by~\ref{item:p_symbol_1}, $x \Equiv \phi(x)$, and by~\ref{item:p_symbol_2},
this implies that $\phi(x) = \phi(\phi(x))$. For this reason, $\phi$ is idempotent.
Moreover, observe that for any $x, x' \in X$, if $\phi(x) = \phi\Par{x'}$, then
by~\ref{item:p_symbol_1},
\begin{math}
    x \Equiv \phi(x) = \phi\Par{x'} \Equiv x',
\end{math}
which implies $x \Equiv x'$. Therefore, the converse of~\ref{item:p_symbol_2} holds.

In the other direction, given a map $\phi : X \to X$, the \Def{kernel} of $\phi$ is the
equivalence relation $\Equiv$ on $X$ such that for any $x, x' \in X$, $x \Equiv x'$ whenever
$\phi(x) = \phi\Par{x'}$.

\begin{Statement}{Proposition}{prop:p_symbol_fiber_equivalence_relations}
    Let $X$ be a set and $\phi : X \to X$ be a map. If $\phi$ is idempotent, then the map
    $\phi$ is a $\PSymbol$-symbol for the kernel of $\phi$.
\end{Statement}
\begin{Proof}
    Let $\Equiv$ be the kernel of $\phi$. The map $\phi$ satisfies
    Condition~\ref{item:p_symbol_2} immediately by construction of $\Equiv$. Besides, since
    $\phi$ is idempotent, for any $x \in X$, we have $\phi(x) = \phi(\phi(x))$ so that $x
    \Equiv \phi(x)$. Therefore, Condition~\ref{item:p_symbol_1} holds.
\end{Proof}

\subsubsection{Quotient clones}
Let us now consider $\PSymbol$-symbols in the context of clones, with the aim of
constructing and studying clone congruences.

\begin{Statement}{Proposition}{prop:p_symbol_congruence}
    Let $\Clone$ be a clone, $\Equiv$ be a graded set equivalence relation on $\Clone$, and
    $\phi$ be a $\PSymbol$-symbol for $\Equiv$. The equivalence relation $\Equiv$ is a clone
    congruence of $\Clone$ if and only if for any $x \in \Clone(n)$, $n \geq 0$, and $x'_1,
    \dots, x'_n \in \Clone(m)$, $m
    \geq 0$,
    \begin{equation} \label{equ:p_symbol_congruence}
        x \Superposition{x'_1, \dots, x'_n}
        \Equiv
        \phi(x) \Superposition{\phi\Par{x'_1}, \dots, \phi\Par{x'_n}}.
    \end{equation}
\end{Statement}
\begin{Proof}
    If $\Equiv$ is a clone congruence of $\Clone$, \eqref{equ:p_symbol_congruence} holds by
    the fact that since $\phi$ is a $\PSymbol$-symbol for $\Equiv$, $\phi$ satisfies
    Condition~\ref{item:p_symbol_1} of $\PSymbol$-symbols.

    Conversely, let us assume that~\eqref{equ:p_symbol_congruence} holds. Let $x \in
    \Clone(n)$, $y \in \Clone(n)$, $n \geq 0$, and $x'_1, \dots, x'_n \in \Clone(m)$, $y'_1,
    \dots, y'_n \in \Clone(m)$, $m \geq 0$, such that $x \Equiv y$ and $x'_i \Equiv y'_i$
    for all $i \in [n]$. Therefore, by Condition~\ref{item:p_symbol_2} of
    $\PSymbol$-symbols, $\phi(x) = \phi(y)$ and $\phi\Par{x'_i} = \phi\Par{y'_i}$ for all $i
    \in [n]$, so that
    \begin{math}
        \phi(x) \Superposition{\phi\Par{x'_1}, \dots, \phi\Par{x'_n}}
        = \phi(y) \Superposition{\phi\Par{y'_1}, \dots, \phi\Par{y'_n}}.
    \end{math}
    By~\eqref{equ:p_symbol_congruence}, this implies that
    \begin{math}
        x \Superposition{x'_1, \dots, x'_n} \Equiv y \Superposition{y'_1, \dots, y'_n}
    \end{math}
    and shows as expected that $\Equiv$ is a clone congruence of~$\Clone$.
\end{Proof}

The following result provides a concrete description of the quotient $\Clone /_{\Equiv}$ of
a clone $\Clone$ by a clone congruence $\Equiv$, assuming the existence of a
$\PSymbol$-symbol for~$\phi$.

\begin{Statement}{Proposition}{prop:clone_realization_p_symbol}
    Let $\Clone$ be a clone, $\Equiv$ be a clone congruence of $\Clone$, and $\phi$ be a
    $\PSymbol$-symbol for $\Equiv$. The clone $\Clone /_{\Equiv}$ is isomorphic to the clone
    on $\phi(\Clone)$ with superposition maps defined, for any $x \in \phi(\Clone)(n)$, $n
    \geq 0$, and $x'_1, \dots, x'_n \in \phi(\Clone)(m)$, $m \geq 0$, by
    \begin{equation} \label{equ:superposition_p_symbol}
        x \Superposition{x'_1, \dots, x'_n}
        := \phi\Par{x \Superposition{x'_1, \dots, x'_n}},
    \end{equation}
    where the superposition map of the right-hand side of~\eqref{equ:superposition_p_symbol}
    is the one of $\Clone$, and the projections of this clone are the images by $\phi$ of
    the projections of $\Clone$.
\end{Statement}
\begin{Proof}
    Let $\bar{\phi} : \Clone /_{\Equiv} \to \phi(\Clone)$ be the map defined for any $x \in
    \Clone$ by $\bar{\phi}\Par{\Han{x}_{\Equiv}} := \phi(x)$. Since $\Equiv$ and $\phi$
    satisfy~\ref{item:p_symbol_2}, $\bar{\phi}$ is a well-defined map. Moreover,
    by~\ref{item:p_symbol_1}, $\bar{\phi}$ is surjective, and it follows from the converse
    of~\ref{item:p_symbol_2} (which holds, as noticed in Section~\ref{subsubsec:p_symbols})
    that $\bar{\phi}$ is injective. Let us prove that $\bar{\phi}$ is a clone morphism.
    For any $x \in \phi(\Clone)(n)$, $n
    \geq 0$, and $x'_1, \dots, x'_n \in \phi(\Clone)(m)$, $m \geq 0$, we have
    \begin{align} \label{equ:clone_realization_p_symbol}
        \bar{\phi}\Par{
            \Han{x}_{\Equiv}
            \Superposition{
                \Han{x'_1}_{\Equiv}, \dots, \Han{x'_n}_{\Equiv}
            }
        }
        & =
        \bar{\phi}\Par{\Han{x \Superposition{x'_1, \dots, x'_n}}_{\Equiv}}
        \\
        & = \phi\Par{x \Superposition{x'_1, \dots, x'_n}}
        \notag
        \\
        & = x \Superposition{x'_1, \dots, x'_n}
        \notag
        \\
        & = \phi(x) \Superposition{\phi\Par{x'_1}, \dots, \phi\Par{x'_n}}
        \notag
        \\
        & = \bar{\phi}\Par{\Han{x}_{\Equiv}}
        \Superposition{
            \bar{\phi}\Par{\Han{x'_1}_{\Equiv}},
            \dots,
            \bar{\phi}\Par{\Han{x'_n}_{\Equiv}}
        }.
        \notag
    \end{align}
    The first equality of~\eqref{equ:clone_realization_p_symbol} comes from the fact
    $\Equiv$ is a clone congruence, the second and fifth are by definition of $\bar{\phi}$,
    the third is by definition of the superposition maps of $\phi(\Clone)$ provided by the
    statement of the proposition, and the fourth comes from the fact that since $\phi$ is
    idempotent, each element of $\phi(\Clone)$ is a fixed-point of~$\phi$. Therefore,
    $\bar{\phi}$ is a clone isomorphism and the statement of the proposition follows.
\end{Proof}

\subsubsection{Composition of $\PSymbol$-symbols}
Let us focus now on the compositions of $\PSymbol$-symbols and on the properties of the
resulting maps.

\begin{Statement}{Proposition}{prop:p_symbol_composition}
    Let $\Clone$ be a clone, $\Equiv_1$ and $\Equiv_2$ be two clone congruences of $\Clone$,
    and $\phi_1$ and $\phi_2$ be two $\PSymbol$-symbols, respectively for $\Equiv_1$ and
    $\Equiv_2$. If $\phi_1$ and $\phi_2$ commute w.r.t.\ the composition of maps, then by
    setting $\phi_{12}$ as the map $\phi_1 \circ \phi_2 = \phi_2 \circ \phi_1$ and $\Equiv$
    as the kernel of $\phi_{12}$,
    \begin{enumerate}[label=(\roman*)]
        \item \label{item:p_symbol_composition_1}
        the map $\phi_{12}$ is a $\PSymbol$-symbol for $\Equiv$;
        \item \label{item:p_symbol_composition_2}
        the equivalence relation $\Equiv$ is a clone congruence of $\Clone$;
        \item \label{item:p_symbol_composition_3}
        the clone $\Clone /_{\Equiv}$ is a quotient of both $\Clone /_{\Equiv_1}$ and
        $\Clone /_{\Equiv_2}$.
    \end{enumerate}
\end{Statement}
\begin{Proof}
    In this proof, in order to lighten the notation, for any word $w \in [2]^*$, we denote
    by~$\phi_w$ the map
    \begin{math}
        \phi_{w(1)} \circ \dots \circ \phi_{w(\Length(w))}.
    \end{math}

    Let us first show~\ref{item:p_symbol_composition_1}. Since $\phi_1$ and $\phi_2$ are
    $\PSymbol$-symbols, they are idempotent. Moreover, by hypothesis, they commute w.r.t\
    the composition of maps. Thus, we have
    \begin{math}
        \phi_{12} \circ \phi_{12} = \phi_{1212} = \phi_{1122} = \phi_{12}.
    \end{math}
    Therefore, $\phi_{12}$ is idempotent, implying by
    Proposition~\ref{prop:p_symbol_fiber_equivalence_relations} that $\phi_{12}$ is a
    $\PSymbol$-symbol for~$\Equiv$.

    Let us prove~\ref{item:p_symbol_composition_2}. Since $\phi_1$ and $\phi_2$ are
    respectively $\PSymbol$-symbols for the congruences $\Equiv_1$ and $\Equiv_2$ of
    $\Clone$, and $\phi_1$ and $\phi_2$ commute w.r.t.\ the composition of maps, by
    Proposition~\ref{prop:p_symbol_congruence}, for any $x \in \Clone(n)$, $n \geq
    0$, and $x'_1, \dots, x'_n \in \Clone(m)$, $m \geq 0$, we have
    \begin{align}
        \phi_{12} \Par{x \Superposition{x'_1, \dots x'_n}}
        & = \phi_{12} \Par{\phi_2(x) \Superposition{\phi_2\Par{x'_1}, \dots,
        \phi_2\Par{x'_n}}}
        \\
        & = \phi_{21} \Par{\phi_2(x) \Superposition{\phi_2\Par{x'_1}, \dots,
        \phi_2\Par{x'_n}}}
        \notag \\
        & = \phi_{21} \Par{\phi_{12}(x) \Superposition{\phi_{12}\Par{x'_1},
            \dots, \phi_{12}\Par{x'_n}}}
        \notag \\
        & = \phi_{12} \Par{\phi_{12}(x) \Superposition{\phi_{12}\Par{x'_1},
            \dots, \phi_{12}\Par{x'_n}}}.
        \notag
    \end{align}
    Therefore,
    \begin{math}
        x \Superposition{x'_1, \dots x'_n}
        \Equiv \phi_{12}(x) \Superposition{\phi_{12}\Par{x'_1},
        \dots, \phi_{12}\Par{x'_n}},
    \end{math}
    and by~\ref{item:p_symbol_1} and Proposition~\ref{prop:p_symbol_congruence}, $\Equiv$ is
    a clone congruence of $\Clone$.

    To show~\ref{item:p_symbol_composition_3}, let $x, x' \in \Clone(n)$, $n \geq 0$, such
    that $x \Equiv_1 x'$. Since $\phi_1$ is a $\PSymbol$-symbol for $\Equiv_1$, we have
    $\phi_1(x) = \phi_1\Par{x'}$, so that $\phi_{21}(x) = \phi_{21}\Par{x'}$. Since $\phi_1$
    and $\phi_2$ commute, this shows that $\phi_{12}(x) = \phi_{12}\Par{x'}$. Hence, we have
    $x \Equiv x'$. The same argument shows that $x \Equiv_2 x'$ implies $x \Equiv x'$.
    Therefore, as equivalence relations, $\Equiv$ is coarser than both $\Equiv_1$ and
    $\Equiv_2$. By~\ref{item:p_symbol_composition_2}, $\Equiv$ is a clone congruence of
    $\Clone$ so that $\Clone /_{\Equiv}$ is a well-defined quotient of $\Clone$. The
    statement follows from the first isomorphism theorem.
\end{Proof}

\subsection{Congruences of the clone of pigmented words} \label{subsec:clone_congruences}
Two maps $\Sort_{\Leq}$ and $\First_k$ from $\P(\Monoid)$ to $\P(\Monoid)$ are introduced.
These maps and some of their compositions lead through their kernels to clone congruences of
$\P(\Monoid)$.

In this section, $\Monoid$ is any monoid but in order to give concrete examples here, we
shall consider $\Monoid$ as the free monoid $\Par{A^*, \Conc, \epsilon}$ where $A$ is the
alphabet $\{a, b, c\}$.

\subsubsection{Reversions of congruences}
We start by introducing an involutive transformation on clone congruences of $\P(\Monoid)$.
For any graded set equivalence relation $\Equiv$ of $\P(\Monoid)$, the \Def{reversion} of
$\Equiv$ is the equivalence relation $\Equiv^\Reverse$ on $\P(\Monoid)$ satisfying, for any
$\PWord, \PWord' \in \P(\Monoid)$, $\PWord \Equiv^\Reverse \PWord'$ if $\Reverse(\PWord)
\Equiv \Reverse\Par{\PWord'}$.

\begin{Statement}{Proposition}{prop:reversed_congruence}
    Let $\Monoid$ be a monoid. If $\Equiv$ is a clone congruence of $\P(\Monoid)$, then
    \begin{enumerate}[label=(\roman*)]
        \item \label{item:reversed_congruence_1}
        the equivalence relation $\Equiv^\Reverse$ is a clone congruence of $\P(\Monoid)$;
        \item \label{item:reversed_congruence_2}
        the map $\bar{\Reverse} : \P(\Monoid) /_{\Equiv} \to \P(\Monoid)
        /_{\Equiv^\Reverse}$ defined for any $\PWord \in \P(\Monoid)$ by
        $\bar{\Reverse}\Par{\Han{\PWord}_{\Equiv}} := \Bra{\Reverse\Par{\PWord'} : \PWord
        \Equiv \PWord'}$, is a clone isomorphism.
    \end{enumerate}
\end{Statement}
\begin{Proof}
    Let $\PWord, \QWord \in \P(\Monoid)(n)$, $n \geq 0$, and $\PWord'_1, \dots, \PWord'_n,
    \QWord'_1, \dots, \QWord'_n \in \P(\Monoid)(m)$, $m \geq 0$, such that $\PWord
    \Equiv^\Reverse \QWord$ and $\PWord'_i \Equiv^\Reverse \QWord'_i$ for all $i \in [n]$.
    By definition of $\Equiv^\Reverse$ and since $\Reverse$ is an involution, we have
    $\Reverse(\PWord) \Equiv \Reverse(\QWord)$ and $\Reverse\Par{\PWord'_i} \Equiv
    \Reverse\Par{\QWord'_i}$ for all $i \in [n]$. Now, since $\Equiv$ is a clone of
    congruence of $\P(\Monoid)$,
    \begin{equation}
        \Reverse(\PWord) \Superposition{
            \Reverse\Par{\PWord'_1}, \dots, \Reverse\Par{\PWord'_n}}
        \Equiv
        \Reverse(\QWord) \Superposition{
            \Reverse\Par{\QWord'_1}, \dots, \Reverse\Par{\QWord'_n}}.
    \end{equation}
    This implies, since by Proposition~\ref{prop:reverse_map_clone_isomorphism}, $\Reverse$
    is a clone isomorphism of $\P(\Monoid)$, that
    \begin{equation}
        \Reverse\Par{\PWord \Superposition{\PWord'_1, \dots, \PWord'_n}}
        \Equiv
        \Reverse\Par{\QWord \Superposition{\QWord'_1, \dots, \QWord'_n}}.
    \end{equation}
    Therefore, by definition of $\Equiv^\Reverse$, this shows that $\PWord
    \Superposition{\PWord'_1, \dots, \PWord'_n}$ is $\Equiv^\Reverse$-equivalent to $\QWord
    \Superposition{\QWord'_1, \dots, \QWord'_n}$,
    establishing~\ref{item:reversed_congruence_1}.

    To prove~\ref{item:reversed_congruence_2}, observe first that since $\Reverse$ is an
    involution of $\P(\Monoid)$, by definition of $\Equiv^\Reverse$, for any $\PWord \in
    \P(\Monoid)$,
    \begin{equation} \label{equ:reversed_congruence_1}
        \bar{\Reverse}\Par{\Han{\PWord}_{\Equiv}}
        = \Bra{\Reverse\Par{\PWord'} : \PWord \Equiv \PWord'}
        = \Bra{\Reverse\Par{\PWord'} :
            \Reverse(\PWord) \Equiv^\Reverse \Reverse\Par{\PWord'}}
        = \Bra{\PWord' : \Reverse(\PWord) \Equiv^\Reverse \PWord'}
        = \Han{\Reverse(\PWord)}_{\Equiv^\Reverse}.
    \end{equation}
    Therefore, the map $\bar{\Reverse}$ from $\P(\Monoid) /_{\Equiv}$ to $\P(\Monoid)
    /_{\Equiv^\Reverse}$ is well-defined and is bijective. Now, by using consecutively the
    fact that $\Equiv$ is a clone congruence of $\P(\Monoid)$,
    Relation~\eqref{equ:reversed_congruence_1}, the fact that by
    Proposition~\ref{prop:reverse_map_clone_isomorphism}, $\Reverse$ is an endomorphism of
    $\P(\Monoid)$, and the fact that by~\ref{item:reversed_congruence_1}, $\Equiv^\Reverse$
    is a clone congruence of $\P(\Monoid)$, for any $\PWord \in \P(\Monoid)(n)$, $n \geq 0$,
    and $\PWord_1', \dots, \PWord_n' \in \P(\Monoid)(m)$, $m \geq 0$, we have
    \begin{align}
        \bar{\Reverse}\Par{\Han{\PWord}_{\Equiv}
            \Superposition{\Han{\PWord_1'}_{\Equiv}, \dots, \Han{\PWord_n'}_{\Equiv}}}
        & = \bar{\Reverse}\Par{
            \Han{\PWord \Superposition{\PWord_1', \dots, \PWord_n'}}_{\Equiv}}
        \\
        & = \Han{\Reverse\Par{\PWord \Superposition{
            \PWord_1', \dots, \PWord_n'}}}_{\Equiv^\Reverse}
        \notag \\
        & = \Han{\Reverse\Par{\PWord} \Superposition{\Reverse\Par{\PWord_1'}, \dots,
            \Reverse\Par{\PWord_n'}}}_{\Equiv^\Reverse}
        \notag \\
        & = \Han{\Reverse(\PWord)}_{\Equiv^\Reverse}
            \Superposition{
                \Han{\Reverse\Par{\PWord_1'}}_{\Equiv^\Reverse},
                \dots,
                \Han{\Reverse\Par{\PWord_n'}}_{\Equiv^\Reverse}}.
        \notag
    \end{align}
    Observe also that, by denoting by $\MonoidUnit$ the unit of $\Monoid$, for any
    $i^\MonoidUnit \in \P(\Monoid)(n)$, $n \geq 1$, $i \in [n]$,
    \begin{math}
        \bar{\Reverse}\Par{\Han{i^\MonoidUnit}_{\Equiv}}
        = \Han{i^\MonoidUnit}_{\Equiv^\Reverse}.
    \end{math}
    Therefore, $\bar{\Reverse}$ is a clone isomorphism from $\P(\Monoid) /_{\Equiv}$ to
    $\P(\Monoid) /_{\Equiv^\Reverse}$.
\end{Proof}

For any clone $\Clone := \P(\Monoid) /_{\Equiv}$ where $\Equiv$ is a clone congruence of
$\P(\Monoid)$, we denote by $\Clone^\Reverse$ the clone $\P(\Monoid) /_{\Equiv^\Reverse}$.
This clone is, by Proposition~\ref{prop:reversed_congruence}, well-defined and isomorphic to
$\Clone$.

\subsubsection{Sorting congruence}
For any total order relation $\Leq$ on $\Monoid$, let $\Sort_{\Leq} : \P(\Monoid) \to
\P(\Monoid)$ be the map sending any $\PWord \in \P(\Monoid)$ to the $\Monoid$-pigmented word
obtained by rearranging the values of $\PWord$ in weakly increasing way w.r.t.\ the total
order relation $\PLeq$ on the set of the $\Monoid$-pigmented letters satisfying
$i_1^{\alpha_1} \PLeq i_2^{\alpha_2}$ if $i_1 < i_2$, or $i_1 = i_2$ and $\alpha_1 \Leq
\alpha_2$. For instance, in $\P\Par{A^*}$, where $\Leq$ is the lexicographic order on $A^*$
satisfying $a \Leq b \Leq c$, we have
\begin{equation} \label{equ:example_sort}
    \Sort_{\Leq}\Par{3^\epsilon 1^b 3^a 1^a 4^{ab} 2^b 3^\epsilon 1^\epsilon}
    = 1^\epsilon 1^a 1^b 2^b 3^\epsilon 3^\epsilon 3^a 4^{ab}.
\end{equation}
Let $\Equiv_{\Sort_{\Leq}}$ be the kernel of $\Sort_{\Leq}$. By
Proposition~\ref{prop:p_symbol_fiber_equivalence_relations}, since $\Sort_{\Leq}$ is
idempotent, $\Sort_{\Leq}$ is a $\PSymbol$-symbol for $\Equiv_{\Sort_{\Leq}}$. Observe
moreover that for any $\PWord, \PWord' \in \P(\Monoid)$, we have $\PWord
\Equiv_{\Sort_{\Leq}} \PWord'$ if and only if the multisets of pigmented letters of $\PWord$
and $\PWord'$ coincide. For this reason, the equivalence relation $\Equiv_{\Sort_{\Leq}}$
does not depend on the total order relation $\Leq$. Therefore, we denote this equivalence
relation simply by~$\Equiv_\Sort$.

\begin{Statement}{Proposition}{prop:congruence_sort}
    For any monoid $\Monoid$, the equivalence relation $\Equiv_\Sort$ is a clone congruence
    of~$\P(\Monoid)$.
\end{Statement}
\begin{Proof}
    Let $\Leq$ be any total order relation on $\Monoid$ and $\PWord \in \P(\Monoid)$. For
    any $i^\alpha \in \SetLetters_\Monoid$, $\PWord$ and $\Sort_{\Leq}(\PWord)$ admit the
    same number of occurrences of $i^\alpha$. For this reason and by the definition of the
    superposition maps of $\P(\Monoid)$, the $\PSymbol$-symbol $\Sort_{\Leq}$ for
    $\Equiv_\Sort$ satisfies the prerequisites of
    Proposition~\ref{prop:p_symbol_congruence}. This implies the statement of the
    proposition.
\end{Proof}

\subsubsection{First occurrences congruence}
For any $k \geq 0$ and any $\PWord \in \P(\Monoid)$, a position $j \in [\Length(\PWord)]$ is
a \Def{left $k$-witness} of $\PWord$ if in $\PWord(1, j - 1)$, there are at most $k - 1$
$\Monoid$-pigmented letters having as value the one of $\PWord(j)$. Similarly, a position $j
\in [\Length(\PWord)]$ is a \Def{right $k$-witness} of $\PWord$ if in $\PWord(j + 1,
\Length(\PWord))$, there are at most $k - 1$ $\Monoid$-pigmented letters having as value the
one of $\PWord(j)$.

We shall highlight these properties by putting a segment with a circle on the left (resp.\
right) under each $\Monoid$-pigmented letter such that its position is a left (resp.\ right)
$k$-witness. In the opposite case, we shall put a cross on the left (resp.\ right) edge of
the segment to highlight the fact that this position is not a left (resp.\ right)
$k$-witness when it is the case. For instance, by setting
\begin{math}
    \PWord := 2^{aa} 2^b 3^a 1^a 3^{ba} 2^b 3^\epsilon,
\end{math}
the left and right $1$-witnesses of $\PWord$ are highlighted as
\begin{equation}
    \Witness{y}{n}{2^{aa}} \Witness{n}{n}{2^b} \Witness{y}{n}{3^a} \Witness{y}{y}{1^a}
    \Witness{n}{n}{3^{ba}} \Witness{n}{y}{2^b} \Witness{n}{y}{3^\epsilon}
\end{equation}
and the left and right $2$-witnesses of $\PWord$ are highlighted as
\begin{equation}
    \Witness{y}{n}{2^{aa}} \Witness{y}{y}{2^b} \Witness{y}{n}{3^a} \Witness{y}{y}{1^a}
    \Witness{y}{y}{3^{ba}} \Witness{n}{y}{2^b} \Witness{n}{y}{3^\epsilon}.
\end{equation}
Moreover, a left (resp.\ right) edge of a segment having neither a circle nor a cross
specifies the fact that the status of this position is unknown. For instance, for a fixed $k
\geq 0$, the notation
\begin{equation}
    \PWord_1 \Conc \Witness{y}{u}{1^{ba}} \Conc \PWord_2 \Conc \Witness{u}{n}{1^{ab}}
    \Witness{u}{y}{1^{b}} \Conc \PWord_3
\end{equation}
where $\PWord_1$, $\PWord_2$, and $\PWord_3$ are some $A^*$-pigmented words specifies an
$A^*$-pigmented word such that the position of the shown $\Monoid$-pigmented letter $1^{ba}$
is a left $k$-witness and may or may not be a right $k$-witness, that the position of the
shown $\Monoid$-pigmented letter $1^{ab}$ may or may not be a left $k$-witness and is not a
right $k$-witness, and that the position of the shown $\Monoid$-pigmented letter $1^{b}$ may
or may not be a left $k$-witness and is a right $k$-witness.

Now, let $\First_k : \P(\Monoid) \to \P(\Monoid)$ be the map sending any $\PWord \in
\P(\Monoid)$ to the $\Monoid$-pigmented word defined as the subword of $\PWord$ consisting
of the letters whose positions are left $k$-witnesses. For instance,
\begin{equation}
    \First_1\Par{
        \Witness{y}{u}{1^\epsilon} \Witness{y}{u}{3^{ab}} \Witness{n}{u}{1^b}
        \Witness{n}{u}{3^b} \Witness{n}{u}{1^{aa}} \Witness{n}{u}{3^\epsilon}
        \Witness{y}{u}{2^{aa}} \Witness{n}{u}{3^{bba}}}
    = \Witness{y}{u}{1^\epsilon} \Witness{y}{u}{3^{ab}} \Witness{y}{u}{2^{aa}},
\end{equation}
\begin{equation}
    \First_2\Par{
        \Witness{y}{u}{1^\epsilon} \Witness{y}{u}{3^{ab}} \Witness{y}{u}{1^b}
        \Witness{y}{u}{3^b} \Witness{n}{u}{1^{aa}} \Witness{n}{u}{3^\epsilon}
        \Witness{y}{u}{2^{aa}} \Witness{n}{u}{3^{bba}}}
    = \Witness{y}{u}{1^\epsilon} \Witness{y}{u}{3^{ab}} \Witness{y}{u}{1^b}
    \Witness{y}{u}{3^b} \Witness{y}{u}{2^{aa}}.
\end{equation}
Let $\Equiv_{\First_k}$ be the kernel of $\First_k$. By
Proposition~\ref{prop:p_symbol_fiber_equivalence_relations}, since $\First_k$ is idempotent,
$\First_k$ is a $\PSymbol$-symbol for $\Equiv_{\First_k}$.

Observe that for any $0 \leq k \leq k'$ and any $\PWord, \PWord' \in \P(\Monoid)$, $\PWord
\Equiv_{\First_{k'}} \PWord'$ implies $\PWord \Equiv_{\First_k} \PWord'$. Hence, the
equivalence relation $\Equiv_{\First_{k'}}$ is a refinement of $\Equiv_{\First_k}$.

\begin{Statement}{Proposition}{prop:congruence_first}
    For any monoid $\Monoid$ and any $k \geq 0$, the equivalence relation
    $\Equiv_{\First_k}$ is a clone congruence of $\P(\Monoid)$.
\end{Statement}
\begin{Proof}
    From the definitions of $\First_k$ and of the superposition maps of $\P(\Monoid)$, for
    any $\PWord \in \P(\Monoid)(n)$, $n \geq 0$, and $\PWord_1, \dots, \PWord_n \in
    \P(\Monoid)(m)$, $m \geq 0$, we have
    \begin{equation}
        \First_k\Par{\PWord \Superposition{\PWord_1, \dots, \PWord_n}}
        = \First_k\Par{\First_k\Par{\PWord} \Superposition{\PWord_1, \dots, \PWord_n}}
    \end{equation}
    and, for any $j \in [n]$,
    \begin{equation}
        \First_k\Par{\PWord \Superposition{\PWord_1, \dots, \PWord_n}}
        = \First_k\Par{\PWord \Superposition{\PWord_1, \dots, \PWord_{j - 1},
            \First_k\Par{\PWord_j}, \PWord_{j + 1}, \dots, \PWord_n}}.
    \end{equation}
    These two properties imply that the $\PSymbol$-symbol $\First_k$ for $\Equiv_{\First_k}$
    satisfies the prerequisites of Proposition~\ref{prop:p_symbol_congruence}. This
    establishes the statement of the proposition.
\end{Proof}

For any $k \geq 0$, let us denote by $\First_k^\Reverse : \P(\Monoid) \to \P(\Monoid)$ the
map defined for any $\PWord \in \P(\Monoid)$ by $\First_k^\Reverse(\PWord) :=
\Reverse\Par{\First_k\Par{\Reverse(\PWord)}}$. In this way, for any $\PWord \in
\P(\Monoid)$, $\First_k^\Reverse(\PWord)$ is the subword of $\PWord$ consisting of the
letters whose positions are right $k$-witnesses. It is straightforward to prove that
$\First_k^\Reverse$ is idempotent and that the kernel of $\First_k^\Reverse$ is the
equivalence relation $\Equiv_{\First_k}^\Reverse$. By
Propositions~\ref{prop:congruence_first} and~\ref{prop:reversed_congruence},
$\Equiv_{\First_k}^\Reverse$ is a clone congruence of~$\P(\Monoid)$.

\subsubsection{Compositions}
We consider here some compositions of the maps $\Sort_{\Leq}$, $\First_k$, and
$\First_{k'}^\Reverse$, $k, k' \geq 0$. Directly from the definition of the map $\First_k$,
for any $k, k' \geq 0$, $\First_k \circ \First_{k'} = \First_{\min \Bra{k, k'}}$. Moreover,
for any $k, k' \geq 0$ such that $k \leq k'$, $\First_{k'} \circ \First_k^\Reverse =
\First_k^\Reverse$ and $\First_{k'}^\Reverse \circ \First_k = \First_k$. Observe also
that the maps $\First_k$ and $\First_{k'}$, $k, k' \geq 0$ do not commute. Indeed, in
$\P(\TrivialMonoid)$, where $\TrivialMonoid = \Bra{\MonoidUnit}$ is the trivial monoid, we
have
\begin{equation}
    \First_1\Par{\First_2^\Reverse\Par{
        2^\MonoidUnit 1^\MonoidUnit 2^\MonoidUnit 1^\MonoidUnit 2^\MonoidUnit}}
    =
    1^\MonoidUnit 2^\MonoidUnit
    \ne
    2^\MonoidUnit 1^\MonoidUnit
    =
    \First_2^\Reverse\Par{\First_1\Par{
        2^\MonoidUnit 1^\MonoidUnit 2^\MonoidUnit 1^\MonoidUnit 2^\MonoidUnit}}.
\end{equation}

\begin{Statement}{Proposition}{prop:commutation_sort_first}
    For any monoid $\Monoid$, any $k \geq 1$, and any total order relation $\Leq$ on
    $\Monoid$, the maps $\Sort_{\Leq}$ and $\First_k$ on $\P(\Monoid)$ commute if and only
    if $\Monoid$ is the trivial monoid $\TrivialMonoid$.
\end{Statement}
\begin{Proof}
    Let $\PWord \in \P(\TrivialMonoid)(n)$, $n \geq 0$. By definition of $\Sort$ and of
    $\First_k$, $\Sort_{\Leq}\Par{\First_k(\PWord)}$ is the $\TrivialMonoid$-pigmented word
    $\QWord$ such that for any $j \in [\Length(\QWord) - 1]$, $\QWord(j) \PLeq \QWord(j +
    1)$, and for any $i^\MonoidUnit \in \SetLetters_\TrivialMonoid$, $\QWord$ has exactly
    $\min \Bra{\Brr{\PWord}_{i^\MonoidUnit}, k}$ occurrences of $i^\MonoidUnit$, where
    $\Brr{\PWord}_{i^\MonoidUnit}$ is the number of occurrences of $i^\MonoidUnit$ in
    $\PWord$. Since $\First_k\Par{\Sort_{\Leq}(\PWord)}$ satisfies the same property, we
    have $\Sort_{\Leq}\Par{\First_k(\PWord)} = \First_k\Par{\Sort_{\Leq}(\PWord)}$.

    Conversely, assume that $\Monoid$ is not trivial. Thus, $\Monoid$ contains two distinct
    elements $\alpha_1$ and $\alpha_2$. By considering without loss of generality that
    $\alpha_1 \Leq \alpha_2$, we have in particular, since $k \geq 1$,
    \begin{align}
        \Sort_{\Leq}\Par{\First_k\Par{
            \underbrace{1^{\alpha_2} \dots 1^{\alpha_2}}_{k \text{ times}} 1^{\alpha_1}
        }}
        & = \underbrace{1^{\alpha_2} \dots 1^{\alpha_2}}_{k \text{ times}}
        \\
        & \ne 1^{\alpha_1} \underbrace{1^{\alpha_2} \dots 1^{\alpha_2}}_{k - 1 \text{ times}}
        = \First_k\Par{\Sort_{\Leq}\Par{
            \underbrace{1^{\alpha_2} \dots 1^{\alpha_2}}_{k \text{ times}} 1^{\alpha_1}
        }}.
        \notag
    \end{align}
    This shows that $\Sort_{\Leq}$ and $\First_k$ do not commute.
\end{Proof}

\subsection{Three non-complicated quotients}
We use the clone congruences introduced in the previous section to build three quotients
$\WInc(\Monoid)$, $\Arra_k(\Monoid)$, and $\Inc_k$ of $\P(\Monoid)$. Each of these clones
admits finitely equationally axiomatizable presentations: the first clone is a clone
realization of a generalization of the variety of commutative monoids, the second one is a
clone realization of a generalization of the variety of left-regular bands, and the last one
is a clone realization of a generalization of the variety of bounded semilattices.

\subsubsection{On pigmented weakly increasing words}
Let
\begin{equation}
    \WInc(\Monoid) := \P(\Monoid) /_{\Equiv_\Sort}.
\end{equation}
By Proposition~\ref{prop:congruence_sort}, $\WInc(\Monoid)$ is a well-defined quotient clone
of $\P(\Monoid)$.

Since $\Sort_{\Leq}$ is a $\PSymbol$-symbol for $\Equiv_\Sort$ where $\Leq$ is any total
order relation on $\Monoid$, the clone $\WInc(\Monoid)$ is described by
Proposition~\ref{prop:clone_realization_p_symbol}. Hence, by definition of $\Sort_{\Leq}$,
$\WInc(\Monoid)$ is a clone on the graded set of \Def{weakly $\Leq$-increasing}
$\Monoid$-pigmented words, which are the $\Monoid$-pigmented words $\PWord$ such that, for
any $j \in [\Length(\PWord) - 1]$, $\PWord(j) \PLeq \PWord(j + 1)$. Equivalently, the
elements of $\WInc(\Monoid)$ can be seen as multisets of $\Monoid$-pigmented letters. For
instance, in $\WInc\Par{A^*}$, where $\Leq$ is the order on $A^*$ used
in~\eqref{equ:example_sort}, we have, up to isomorphism,
\begin{align}
    2^{ab} 3^\epsilon 3^a 4^b 4^b \Superposition{
        1^{ab} 2^{ba},
        1^b 2^{ba} 3^\epsilon 3^b,
        1^\epsilon 2^b,
        3^b}
    & =
    \Sort_{\Leq}\Par{
        1^{abb} 2^{abba} 3^{ab} 3^{abb} 1^\epsilon 2^b 1^a 2^{ab} 3^{bb} 3^{bb}}
    \\
    & =
    1^\epsilon 1^a 1^{abb} 2^{ab} 2^{abba} 2^b 3^{ab} 3^{abb} 3^{bb} 3^{bb}.
    \notag
\end{align}
Besides, the clone $\WInc(\Monoid)$ is not combinatorial because
\begin{math}
    \Bra{\PWordEmpty, 1^\MonoidUnit, 1^\MonoidUnit 1^\MonoidUnit, \dots}
    \subseteq \WInc(\Monoid)(1)
\end{math}
where $\MonoidUnit$ is the unit of $\Monoid$.

\begin{Statement}{Proposition}{prop:presentation_winc}
    For any monoid $\Monoid$, the clone $\WInc(\Monoid)$ admits the presentation
    $\Par{\GeneratingSet_\Monoid, \RelationSet_\Monoid'}$ where $\RelationSet_\Monoid'$ is
    the set $\RelationSet_\Monoid$ from Section~\ref{subsubsec:varieties_pigmented_monoids}
    augmented with the $\GeneratingSet_\Monoid$-equation
    \begin{equation} \label{equ:presentation_winc}
        \RightComb_\Monoid\Par{1^\MonoidUnit 2^\MonoidUnit}
        \ \RelationSet_\Monoid' \
        \RightComb_\Monoid\Par{2^\MonoidUnit 1^\MonoidUnit}
    \end{equation}
    that is,
    \begin{math}
        \Pigmentation_\MonoidUnit\Superposition{\VarX_1} \Product
        \Pigmentation_\MonoidUnit\Superposition{\VarX_2} \Product
        \Identity
        \ \RelationSet_\Monoid' \
        \Pigmentation_\MonoidUnit\Superposition{\VarX_2} \Product
        \Pigmentation_\MonoidUnit\Superposition{\VarX_1} \Product
        \Identity,
    \end{math}
    where $\MonoidUnit$ is the unit of $\Monoid$.
\end{Statement}
\begin{Proof}
    Let $\Equiv'$ be the clone congruence of $\P(\Monoid)$ generated by the pair
    \begin{equation} \label{equ:presentation_winc_1}
        1^\MonoidUnit 2^\MonoidUnit \Equiv' 2^\MonoidUnit 1^\MonoidUnit.
    \end{equation}
    Let us show that the clone congruences $\Equiv'$ and $\Equiv_{\Sort}$ of $\P(\Monoid)$
    are equal. This will imply, by the remark stated in
    Section~\ref{subsubsec:remark_presentations}, that $\P(\Monoid) /_{\Equiv'} =
    \P(\Monoid) /_{\Equiv_{\Sort}} = \WInc(\Monoid)$ admits the stated presentation.

    For this, let us introduce some intermediate binary relations on $\P(\Monoid)$. Let
    $\Leq$ be any total order on $\Monoid$ and $\Covering$ be the binary relation on
    $\P(\Monoid)$ satisfying
    \begin{equation}
        \PWord \Conc i_1^{\alpha_1} i_2^{\alpha_2} \Conc \PWord'
        \enspace \Covering \enspace
        \PWord \Conc i_2^{\alpha_2} i_1^{\alpha_1} \Conc \PWord'
        \qquad
        \text{if } i_1^{\alpha_1} \ne i_2^{\alpha_2} \text{ and }
        i_2^{\alpha_2} \PLeq i_1^{\alpha_1},
    \end{equation}
    where $\PWord, \PWord' \in \P(\Monoid)$ and $i_1^{\alpha_1}, i_2^{\alpha_2} \in
    \SetLetters_\Monoid$. Let $\EquivCovering$ be the reflexive, symmetric, and transitive
    closure of $\Covering$ and let us show that $\EquivCovering$ is equal to
    $\Equiv_{\Sort}$. First, observe that directly from the definition of $\Covering$, for
    any $\RWord, \RWord' \in \P(\Monoid)$, $\RWord \Covering \RWord'$ implies
    $\Sort_{\Leq}(\RWord) = \Sort_{\Leq}\Par{\RWord'}$. Hence, we have $\RWord
    \Equiv_{\Sort} \RWord'$, and since $\EquivCovering$ is the smallest equivalence relation
    containing $\Covering$, $\RWord \EquivCovering \RWord'$ implies $\RWord \Equiv_{\Sort}
    \RWord'$. Conversely, assume that $\RWord \Equiv_{\Sort} \RWord'$ for $\RWord, \RWord'
    \in \P(\Monoid)$. By definition of $\Sort_{\Leq}$, for any $\QWord \in \P(\Monoid)$, the
    process consisting in swapping iteratively and as long as possible two adjacent
    $\Monoid$-pigmented letters $i_1^{\alpha_1}$ and $i_2^{\alpha_2}$ of $\QWord$ such that
    $i_1^{\alpha_1} \ne i_2^{\alpha_2}$ and $i_2^{\alpha_2} \PLeq i_1^{\alpha_1}$ finally
    produces the $\Monoid$-pigmented word $\Sort_{\Leq}(\QWord)$. Moreover, observe that by
    definition of $\Covering$, for any $\QWord', \QWord'' \in \P(\Monoid)$, the property
    $\QWord' \Covering \QWord''$ is equivalent to the fact that $\QWord''$ is obtained from
    $\QWord'$ by swapping two adjacent $\Monoid$-pigmented letters $i_1^{\alpha_1}$ and
    $i_2^{\alpha_2}$ such that $i_1^{\alpha_1} \ne i_2^{\alpha_2}$ and $i_2^{\alpha_2} \PLeq
    i_1^{\alpha_1}$. Due to the fact that $\EquivCovering$ is the smallest equivalence
    relation containing $\Covering$, $\RWord \EquivCovering \RWord'$ holds.

    Now, let us show that $\Equiv'$ is equal to $\EquivCovering$. First, since the left-hand
    and the right-hand sides of~\eqref{equ:presentation_winc_1} are
    $\EquivCovering$-equivalent and since $\EquivCovering = \Equiv_{\Sort}$ is a clone
    congruence by Proposition~\ref{prop:congruence_sort}, $\Equiv'$ is contained in
    $\EquivCovering$. Conversely, for any $\PWord, \PWord' \in \P(\Monoid)$ and
    $i_1^{\alpha_1}, i_2^{\alpha_2} \in \SetLetters_\Monoid$, since $\Equiv'$ is a clone
    congruence of $\P(\Monoid)$ containing the pair in~\eqref{equ:presentation_winc_1}, we
    have
    \begin{equation}
        \PWord \Conc i_1^{\alpha_1} i_2^{\alpha_2} \Conc \PWord'
        =
        1^\MonoidUnit 2^\MonoidUnit 3^\MonoidUnit \Superposition{
            \PWord,
            1^\MonoidUnit 2^\MonoidUnit \Superposition{i_1^{\alpha_1}, i_2^{\alpha_2}},
            \PWord'}
        \Equiv'
        1^\MonoidUnit 2^\MonoidUnit 3^\MonoidUnit \Superposition{
            \PWord,
            2^\MonoidUnit 1^\MonoidUnit \Superposition{i_1^{\alpha_1}, i_2^{\alpha_2}},
            \PWord'}
        =
        \PWord \Conc i_2^{\alpha_2} i_1^{\alpha_1} \Conc \PWord'.
    \end{equation}
    This shows that for any $\RWord, \RWord' \in \P(\Monoid)$, $\RWord \Covering \RWord'$
    implies $\RWord \Equiv' \RWord'$. Since $\EquivCovering$ is the smallest equivalence
    relation containing $\Covering$, $\EquivCovering$ is contained in $\Equiv'$. This
    establishes the statement of the proposition.
\end{Proof}

By Proposition~\ref{prop:presentation_winc}, any $\WInc(\Monoid)$-algebra is, up to term
equivalence, an $\Monoid$-pigmented monoid $\Par{\Algebra, \Product, \Identity,
\Par{\Pigmentation_\alpha}_{\alpha \in \Monoid}}$ where~$\Product$ is commutative. In
particular, $\WInc(\TrivialMonoid)$ is a clone realization of the variety of commutative
monoids equipped with an additional unary fundamental operation forced to operate as the
identity map on the monoid.

\subsubsection{On pigmented arrangements}
For any $k \geq 0$, let
\begin{equation}
    \Arra_k(\Monoid) := \P(\Monoid) /_{\Equiv_{\First_k}}.
\end{equation}
By Proposition~\ref{prop:congruence_first}, $\Arra_k(\Monoid)$ is a well-defined quotient
clone of $\P(\Monoid)$. Since for any $0 \leq k \leq k'$, $\Equiv_{\First_{k'}}$ is a
refinement of $\Equiv_{\First_k}$, $\Arra_k(\Monoid)$ is isomorphic to a quotient of
$\Arra_{k'}(\Monoid)$. Moreover, since $\Equiv_{\First_0}$ is the coarsest clone congruence
of $\P(\Monoid)$, $\Arra_0(\Monoid)$ is the trivial clone $\TrivialClone$. Besides, the
clone
\begin{math}
    \Arra_k^\Reverse(\Monoid)
    := \Arra_k(\Monoid)^\Reverse
    = \P(\Monoid) /_{\Equiv_{\First_k}^\Reverse}
\end{math}
is by Proposition~\ref{prop:reversed_congruence} isomorphic to $\Arra_k(\Monoid)$.

Since $\First_k$ is a $\PSymbol$-symbol for $\Equiv_{\First_k}$, the clone
$\Arra_k(\Monoid)$ is described by Proposition~\ref{prop:clone_realization_p_symbol}. Hence,
by definition of $\First_k$, $\Arra_k(\Monoid)$ is isomorphic to a clone on the graded set
of \Def{$\Monoid$-pigmented $k$-arrangements}, which are the $\Monoid$-pigmented words
$\PWord$ such that for any value $i$, there are at most $k$ $\Monoid$-pigmented letters of
$\PWord$ having $i$ as value. For instance, in $\Arra_1\Par{A^*}$, up to isomorphism,
\begin{align}
    2^\epsilon 3^{aa} 1^b 4^{ca}
    \Superposition{3^\epsilon 1^a, 2^{bb}, 2^b 1^a 3^{a}, 1^c 2^c}
    & =
    \First_1\Par{
        2^{bb}
        2^{aab} 1^{aaa} 3^{aaa}
        3^b 1^{ba}
        1^{cac} 2^{cac}
    }
    \\
    & = 2^{bb} 1^{aaa} 3^{aaa},
    \notag
\end{align}
and in $\Arra_2\Par{A^*}$,
\begin{align}
    2^\epsilon 3^{aa} 1^b 4^{ca}
    \Superposition{3^\epsilon 1^a, 2^{bb}, 2^b 1^a 3^{a}, 1^c 2^c}
    & =
   \First_2\Par{
        2^{bb}
        2^{aab} 1^{aaa} 3^{aaa}
        3^b 1^{ba}
        1^{cac} 2^{cac}
    }
    \\
    & = 2^{bb} 2^{aab} 1^{aaa} 3^{aaa} 3^b 1^{ba}.
    \notag
\end{align}
Besides, when $\Monoid$ is finite, $\Arra_k(\Monoid)$ is combinatorial and for any $n \geq
0$,
\begin{equation} \label{equ:enumeration_arra_k}
    \# \Arra_k(\Monoid)(n)
    = \sum_{u \in \HanL{k}^n}
    \frac{\Par{u(1) + \dots + u(n)}!}{u(1)! \dots u(n)!}
    (\# \Monoid)^{u(1) + \dots + u(n)},
\end{equation}
where, as a reminder, $\HanL{k}$ is the set $\Bra{0, 1, \dots, k}$. Let us
explain~\eqref{equ:enumeration_arra_k}. An $\Monoid$-pigmented $k$-arrangement $\PWord$ of
arity $n \geq 0$ and length $\ell \geq 0$ is specified by
\begin{enumerate}[label=(S\arabic*)]
    \item \label{item:enumeration_arra_k_1}
    a word $u \in \HanL{k}^n$ such that $\ell = u(1) + \dots + u(n)$ and for any $i \in
    [n]$, $u(i)$ is the number of occurrences of pigmented letters having $i$ as value in
    $\PWord$;
    \item \label{item:enumeration_arra_k_2}
    a partition $P := \Bra{P_1, \dots, P_n}$ of the set $[\ell]$ of positions of $\PWord$
    where some parts may be empty and such that for any $i \in [n]$, $\# P_i = u(i)$ and
    $P_i$ is the set of positions of pigmented letters having $i$ as value in $\PWord$;
    \item \label{item:enumeration_arra_k_3}
    a word $v$ of length $\ell$ on $\Monoid$ such that for any $j \in [\ell]$, $v(j)$ is the
    pigment of the letter $\PWord(j)$.
\end{enumerate}
The cardinality of $\Arra_k(\Monoid)(n)$ expressed by~\eqref{equ:enumeration_arra_k} follows
from this specification. Indeed, \ref{item:enumeration_arra_k_1} gives rise to the sum over
all possible such words $u$, \ref{item:enumeration_arra_k_2} gives rise to the fraction
which is the multinomial coefficient enumerating all possible such partitions $P$,
and~\ref{item:enumeration_arra_k_3} gives rise to the last term enumerating all possible
such words~$v$.

Moreover, we have in particular
\begin{equation} \label{equ:enumeration_arra_1}
    \# \Arra_1(\Monoid)(n) = \sum_{j \in \HanL{n}} \binom{n}{j} j! (\# \Monoid)^j.
\end{equation}
The sequences of sizes of $\Arra_k(\TrivialMonoid)$ for $k \in \HanL{2}$ start by
\begin{equation}
    1, 1, 1, 1, 1, 1, 1, 1, 1,
    \qquad k = 0,
\end{equation}
\begin{equation}
    1, 2, 5, 16, 65, 326, 1957, 13700, 109601,
    \qquad k = 1,
\end{equation}
\begin{equation}
    1, 3, 19, 271, 7365, 326011, 21295783, 1924223799, 229714292041,
    \qquad k = 2.
\end{equation}
The second and third ones are Sequences~\OEIS{A000522} and~\OEIS{A003011}
of~\cite{Slo}, respectively.

\begin{Statement}{Proposition}{prop:presentation_arra}
    For any monoid $\Monoid$ and any $k \geq 0$, the clone $\Arra_k(\Monoid)$ admits the
    presentation $\Par{\GeneratingSet_\Monoid, \RelationSet_\Monoid'}$ where
    $\RelationSet_\Monoid'$ is the set $\RelationSet_\Monoid$ augmented with the
    $\GeneratingSet_\Monoid$-equations
    \begin{equation} \label{equ:presentation_arra}
        \RightComb_\Monoid\Par{
            1^{\alpha_1} 2^\MonoidUnit 1^{\alpha_2} 3^\MonoidUnit \dots 1^{\alpha_k}
            {(k + 1)}^\MonoidUnit 1^{\alpha_{k + 1}}}
        \ \RelationSet_\Monoid' \
        \RightComb_\Monoid\Par{
            1^{\alpha_1} 2^\MonoidUnit 1^{\alpha_2} 3^\MonoidUnit \dots 1^{\alpha_k}
            {(k + 1)}^\MonoidUnit}
    \end{equation}
    for any $\alpha_1, \alpha_2, \dots, \alpha_k, \alpha_{k + 1} \in \Monoid$ where
    $\MonoidUnit$ is the unit of $\Monoid$.
\end{Statement}
\begin{Proof}
    Let $\Equiv'$ be the clone congruence of $\P(\Monoid)$ generated by the pairs
    \begin{equation} \label{equ:presentation_arra_1}
        1^{\alpha_1} 2^\MonoidUnit 1^{\alpha_2} 3^\MonoidUnit \dots 1^{\alpha_k}
            {(k + 1)}^\MonoidUnit 1^{\alpha_{k + 1}}
        \Equiv'
        1^{\alpha_1} 2^\MonoidUnit 1^{\alpha_2} 3^\MonoidUnit \dots 1^{\alpha_k}
            {(k + 1)}^\MonoidUnit
    \end{equation}
    where $\alpha_1, \alpha_2, \dots, \alpha_k, \alpha_{k + 1} \in \Monoid$. Let us show
    that the clone congruences $\Equiv'$ and $\Equiv_{\First_k}$ of $\P(\Monoid)$ are equal.
    This will imply, by the remark stated in Section~\ref{subsubsec:remark_presentations},
    that $\P(\Monoid) /_{\Equiv'} = \P(\Monoid) /_{\Equiv_{\First_k}} = \Arra_k(\Monoid)$
    admits the stated presentation.

    For this, let us introduce some intermediate binary relations on $\P(\Monoid)$. Let
    $\Covering$ be the binary relation on $\P(\Monoid)$ satisfying
    \begin{equation}
        \PWord \Conc i^{\alpha_1} \Conc \QWord_1 \Conc i^{\alpha_2} \Conc \QWord_2 \Conc
            \ \dots \ \Conc i^{\alpha_k} \Conc \QWord_k \Conc i^{\alpha_{k + 1}} \Conc
            \PWord'
        \enspace \Covering \enspace
        \PWord \Conc i^{\alpha_1} \Conc \QWord_1 \Conc i^{\alpha_2} \Conc \QWord_2 \Conc
            \ \dots \ \Conc i^{\alpha_k} \Conc \QWord_k \Conc \PWord'
    \end{equation}
    where $\PWord, \QWord_1, \QWord_2, \dots, \QWord_k, \PWord' \in \P(\Monoid)$ and
    $i^{\alpha_1}, i^{\alpha_2}, \dots, i^{\alpha_k}, i^{\alpha_{k + 1}} \in
    \SetLetters_\Monoid$. Let $\EquivCovering$ be the reflexive, symmetric, and transitive
    closure of $\Covering$ and let us show that $\EquivCovering$ is equal to
    $\Equiv_{\First_k}$. First, observe that directly from the definition of $\Covering$,
    for any $\RWord, \RWord' \in \P(\Monoid)$, $\RWord \Covering \RWord'$ implies
    $\First_k(\RWord) = \First_k\Par{\RWord'}$. Hence, we have $\RWord \Equiv_{\First_k}
    \RWord'$, and since $\EquivCovering$ is the smallest equivalence relation containing
    $\Covering$, $\RWord \EquivCovering \RWord'$ implies $\RWord \Equiv_{\First_k} \RWord'$.
    Conversely, assume that $\RWord \Equiv_{\First_k} \RWord'$ for $\RWord, \RWord' \in
    \P(\Monoid)$. By definition of $\First_k$, for any $\QWord \in \P(\Monoid)$, the process
    consisting in deleting iteratively and as long as possible each letter of $\QWord$ which
    is not a left $k$-witness finally produces the $\Monoid$-pigmented word
    $\First_k(\QWord)$. Moreover, observe that by definition of $\Covering$, for any
    $\QWord', \QWord'' \in \P(\Monoid)$, the property $\QWord' \Covering \QWord''$ is
    equivalent to the fact that $\QWord''$ is obtained from $\QWord'$ by deleting a letter
    which is not a left $k$-witness. Due to the fact that $\EquivCovering$ is the smallest
    equivalence relation containing $\Covering$, $\RWord \EquivCovering \RWord'$ holds.

    Now, let us show that $\Equiv'$ is equal to $\EquivCovering$. First, since the left-hand
    and right-hand sides of~\eqref{equ:presentation_arra_1} are $\EquivCovering$-equivalent
    and since $\EquivCovering = \Equiv_{\First_k}$ is a clone congruence by
    Proposition~\ref{prop:congruence_first}, $\Equiv'$ is contained in $\EquivCovering$.
    Conversely, for any $\PWord, \QWord_1, \QWord_2, \dots, \QWord_k, \PWord' \in
    \P(\Monoid)$ and $i^{\alpha_1}, i^{\alpha_2}, \dots, i^{\alpha_k}, i^{\alpha_{k + 1}}
    \in \SetLetters_\Monoid$, we have
    \begin{align}
        \PWord \Conc i^{\alpha_1} \Conc \QWord_1 \Conc i^{\alpha_2} \Conc \QWord_2
        \Conc \ \dots \ \Conc i^{\alpha_k} & \Conc \QWord_k \Conc i^{\alpha_{k + 1}} \Conc
        \PWord'
        \\
        & =
        1^\MonoidUnit 2^\MonoidUnit 3^\MonoidUnit \Superposition{
            \PWord,
                1^{\alpha_1} 2^\MonoidUnit 1^{\alpha_2} 3^\MonoidUnit \dots 1^{\alpha_k}
                {(k + 1)}^\MonoidUnit 1^{\alpha_{k + 1}} \Superposition{
                    i^\MonoidUnit, \QWord_1, \QWord_2, \dots, \QWord_k},
            \PWord'}
        \notag \\
        & \Equiv'
        1^\MonoidUnit 2^\MonoidUnit 3^\MonoidUnit \Superposition{
            \PWord,
                1^{\alpha_1} 2^\MonoidUnit 1^{\alpha_2} 3^\MonoidUnit \dots 1^{\alpha_k}
                {(k + 1)}^\MonoidUnit \Superposition{
                    i^\MonoidUnit, \QWord_1, \QWord_2, \dots, \QWord_k},
            \PWord'}
        \notag \\
        & =
        \PWord \Conc i^{\alpha_1} \Conc \QWord_1 \Conc i^{\alpha_2} \Conc \QWord_2
        \Conc \ \dots \ \Conc i^{\alpha_k} \Conc \QWord_k \Conc \PWord'.
        \notag
    \end{align}
    This shows that for any $\RWord, \RWord' \in \P(\Monoid)$, $\RWord \Covering \RWord'$
    implies $\RWord \Equiv' \RWord'$. Since $\EquivCovering$ is the smallest equivalence
    relation containing $\Covering$, $\EquivCovering$ is contained in $\Equiv'$. This
    establishes the statement of the proposition.
\end{Proof}

By Proposition~\ref{prop:presentation_arra}, any $\Arra_k(\Monoid)$-algebra is, up to term
equivalence, an $\Monoid$-pigmented monoid $\Par{\Algebra, \Product, \Identity,
\Par{\Pigmentation_\alpha}_{\alpha \in \Monoid}}$ where~$\Product$ and
$\Pigmentation_\alpha$ satisfy, by spelling out~\eqref{equ:presentation_arra} and
simplifying it modulo the background equational theory $\Equiv_{\RelationSet_\Monoid}$,
\begin{align}
    \Pigmentation_{\alpha_1} \Par{x_1} \Product x_2 \Product
    \Pigmentation_{\alpha_2} \Par{x_1} \Product x_3 \Product
    \cdots \Product \Pigmentation_{\alpha_k} \Par{x_1} \Product x_{k + 1} \Product
    \Pigmentation_{\alpha_{k + 1}} \Par{x_1}
    \\
    =
    \Pigmentation_{\alpha_1} \Par{x_1} \Product x_2 \Product
    \Pigmentation_{\alpha_2} \Par{x_1} \Product x_3 \Product
    \cdots \Product \Pigmentation_{\alpha_k} \Par{x_1} \Product x_{k + 1}
    \notag
\end{align}
for any $x_1, \dots, x_{k + 1} \in \Algebra$ and $\alpha_1, \dots, \alpha_k, \alpha_{k + 1}
\in \Monoid$. In particular, $\Arra_1(\TrivialMonoid)$ is a clone realization of the variety
of left-regular bands equipped with an additional unary operation acting identically. A
left-regular band is a monoid $(\Algebra, \Product, \Identity)$ such that $\Product$
satisfies $x_1 \Product x_2 \Product x_1 = x_1 \Product x_2$ for any~$x_1, x_2 \in
\Algebra$.

\subsubsection{On increasing monochrome words}
Let us denote by $\Leq$ the unique order relation on the trivial monoid $\TrivialMonoid$. By
Proposition~\ref{prop:congruence_sort} (resp.\ \ref{prop:congruence_first}), $\Equiv_\Sort$
(resp.\ $\Equiv_{\First_k}$) is a clone congruence of $\P(\TrivialMonoid)$, and by
Proposition~\ref{prop:p_symbol_fiber_equivalence_relations}, $\Sort_{\Leq}$ (resp.\
$\First_k$) is a $\PSymbol$-symbol for $\Equiv_\Sort$ (resp.\ $\Equiv_{\First_k}$).
Therefore, by Propositions~\ref{prop:commutation_sort_first}
and~\ref{prop:p_symbol_composition}, the map $\phi_k := \Sort_{\Leq} \circ \First_k =
\First_k \circ \Sort_{\Leq}$ is a $\PSymbol$-symbol for the kernel $\Equiv_{\phi_k}$ of
$\phi_k$, and $\Equiv_{\phi_k}$ is a clone congruence of $\P(\TrivialMonoid)$.

For any $k \geq 0$, let
\begin{equation}
    \Inc_k := \P(\TrivialMonoid) /_{\Equiv_{\phi_k}}.
\end{equation}
For the previous reasons, $\Inc_k$ is a well-defined quotient of $\P(\TrivialMonoid)$.
Moreover, since for any $0 \leq k \leq k'$ and any $\PWord, \PWord' \in \P(\TrivialMonoid)$,
$\PWord \Equiv_{\phi_{k'}} \PWord'$ implies $\PWord \Equiv_{\phi_k} \PWord'$, the
equivalence relation $\Equiv_{\phi_{k'}}$ is a refinement of $\Equiv_{\phi_k}$. Therefore,
$\Inc_k$ is isomorphic to a quotient of $\Inc_{k'}$. Besides, since $\Equiv_{\phi_0}$ is the
coarsest clone congruence of $\P(\TrivialMonoid)$, $\Inc_0$ is the trivial clone
$\TrivialClone$.

Since $\phi_k$ is a $\PSymbol$-symbol for $\Equiv_{\phi_k}$, the clone $\Inc_k$ is described
by Proposition~\ref{prop:clone_realization_p_symbol}. Hence, by definition of $\phi_k$,
$\Inc_k$ is isomorphic to a clone on the set of \Def{monochrome $k$-increasing words}, which
are the $\TrivialMonoid$-pigmented words $\PWord$ such that $\PWord$ is weakly
$\PLeq$-increasing and for any value $i$, $\PWord$ has at most $k$ occurrences of
$i^\MonoidUnit$. Equivalently, the elements of $\Inc_k$ can be seen as multisets of positive
integers where each element has multiplicity at most~$k$. For instance, in $\Inc_1$, up to
isomorphism,
\begin{equation}
    1^\MonoidUnit 3^\MonoidUnit \Superposition{
        2^\MonoidUnit 4^\MonoidUnit,
        1^\MonoidUnit 3^\MonoidUnit 4^\MonoidUnit,
        2^\MonoidUnit
     }
     = 2^\MonoidUnit 4^\MonoidUnit,
\end{equation}
and in $\Inc_2$,
\begin{equation}
    1^\MonoidUnit 3^\MonoidUnit \Superposition{
        2^\MonoidUnit 4^\MonoidUnit,
        1^\MonoidUnit 3^\MonoidUnit 4^\MonoidUnit,
        2^\MonoidUnit
     }
     = 2^\MonoidUnit 2^\MonoidUnit 4^\MonoidUnit.
\end{equation}
Besides, $\Inc_k$ is combinatorial and for any $n \geq 0$,
\begin{math}
    \# \Inc_k(n) = (k + 1)^n.
\end{math}

The clone $\Inc_k$ is not parameterized by a monoid $\Monoid$ since, as shown by
Proposition~\ref{prop:commutation_sort_first}, we had to choose $\Monoid = \TrivialMonoid$
to ensure that $\Sort_{\Leq}$ and $\First_k$ commute in order to guarantee that
$\Equiv_{\phi_k}$ is a clone congruence of~$\P(\Monoid)$.

\begin{Statement}{Proposition}{prop:presentation_inc}
    For any $k \geq 0$, the clone $\Inc_k$ admits the presentation
    $\Par{\GeneratingSet_\TrivialMonoid, \RelationSet_\TrivialMonoid'}$ where
    $\RelationSet_\TrivialMonoid'$ is the set $\RelationSet_\TrivialMonoid$ augmented with
    the $\GeneratingSet_\TrivialMonoid$-equations
    \begin{equation} \label{equ:presentation_inc_1}
        \RightComb_\TrivialMonoid\Par{1^\MonoidUnit 2^\MonoidUnit}
        \ \RelationSet_\TrivialMonoid' \
        \RightComb_\TrivialMonoid\Par{2^\MonoidUnit 1^\MonoidUnit},
    \end{equation}
    \begin{equation} \label{equ:presentation_inc_2}
        \RightComb_\TrivialMonoid\Par{
            \underbrace{1^\MonoidUnit \dots 1^\MonoidUnit}_{k + 1 \text{ times}}
        }
        \ \RelationSet_\TrivialMonoid' \
        \RightComb_\TrivialMonoid\Par{
            \underbrace{1^\MonoidUnit \dots 1^\MonoidUnit}_{k \text{ times}}
        }
    \end{equation}
    where $\MonoidUnit$ is the unique element of $\TrivialMonoid$.
\end{Statement}
\begin{Proof}
    Let $\Equiv'$ be the clone congruence of $\P(\TrivialMonoid)$ generated by
    \begin{equation} \label{equ:presentation_inc_1_1}
        1^\MonoidUnit 2^\MonoidUnit \Equiv' 2^\MonoidUnit 1^\MonoidUnit,
    \end{equation}
    \begin{equation} \label{equ:presentation_inc_1_2}
        \underbrace{1^\MonoidUnit \dots 1^\MonoidUnit}_{k + 1 \text{ times}}
        \Equiv'
        \underbrace{1^\MonoidUnit \dots 1^\MonoidUnit}_{k \text{ times}}.
    \end{equation}
    Let us show that the clone congruences $\Equiv'$ and $\Equiv_{\phi_k}$ of
    $\P(\TrivialMonoid)$ are equal. This will imply, by the remark stated in
    Section~\ref{subsubsec:remark_presentations}, that $\P(\TrivialMonoid) /_{\Equiv'} =
    \P(\TrivialMonoid) /_{\phi_k} = \Inc_k$ admits the stated presentation.

    For this, let us introduce some intermediate binary relations on $\P(\TrivialMonoid)$.
    Let $\Covering$ be the binary relation on $\P(\TrivialMonoid)$ satisfying
    \begin{equation}
        \PWord \Conc i_1^\MonoidUnit i_2^\MonoidUnit \Conc \PWord'
        \enspace \Covering \enspace
        \PWord \Conc i_2^\MonoidUnit i_1^\MonoidUnit \Conc \PWord'
        \qquad
        \text{if } i_2 < i_1,
    \end{equation}
    \begin{equation}
        \PWord \Conc
        \underbrace{1^\MonoidUnit \dots 1^\MonoidUnit}_{k + 1 \text{ times}}
        \Conc \PWord'
        \enspace \Covering \enspace
        \PWord \Conc
        \underbrace{1^\MonoidUnit \dots 1^\MonoidUnit}_{k \text{ times}}
        \Conc \PWord',
    \end{equation}
    where $\PWord, \PWord' \in \P(\TrivialMonoid)$ and $i^\MonoidUnit, i_1^\MonoidUnit,
    i_2^\MonoidUnit \in \SetLetters_\TrivialMonoid$. Let $\EquivCovering$ be the reflexive,
    symmetric, and transitive closure of $\Covering$ and let us show that $\EquivCovering$
    is equal to $\Equiv_{\phi_k}$. First, observe that directly from the definition of
    $\Covering$, for any $\RWord, \RWord' \in \P(\TrivialMonoid)$, $\RWord \Covering
    \RWord'$ implies $\phi_k(\RWord) = \phi_k\Par{\RWord'}$. Hence, we have $\RWord
    \Equiv_{\phi_k} \RWord'$, and since $\EquivCovering$ is the smallest equivalence
    relation containing $\Covering$, $\RWord \EquivCovering \RWord'$ implies $\RWord
    \Equiv_{\phi_k} \RWord'$. Conversely, assume that $\RWord \Equiv_{\phi_k} \RWord'$ for
    $\RWord, \RWord' \in \P(\TrivialMonoid)$. By definition of $\phi_k$, for any $\QWord \in
    \P(\Monoid)$, the process consisting in swapping iteratively and as long as possible two
    adjacent $\TrivialMonoid$-pigmented letters $i_1^\MonoidUnit$ and $i_2^\MonoidUnit$ of
    $\QWord$ such that $i_2 < i_1$ and then by deleting iteratively and as long as possible
    each $\TrivialMonoid$-pigmented letter $i^\MonoidUnit$ having on its left $k$
    occurrences of $i^\MonoidUnit$ finally produces the $\TrivialMonoid$-pigmented word
    $\phi_k(\QWord)$. Moreover, observe that by definition of $\Covering$, for any $\QWord',
    \QWord'' \in \P(\TrivialMonoid)$, the property $\QWord' \Covering \QWord''$ is
    equivalent to the fact that $\QWord''$ is obtained from $\QWord'$ swapping two adjacent
    $\TrivialMonoid$-pigmented letters $i_1^\MonoidUnit$ and $i_2^\MonoidUnit$ such that
    $i_2 < i_1$ or by deleting iteratively each $\TrivialMonoid$-pigmented letter
    $i^\MonoidUnit$ having on its left $k$ occurrences of $i^\MonoidUnit$. Due to the fact
    that $\EquivCovering$ is the smallest equivalence relation containing $\Covering$,
    $\RWord \EquivCovering \RWord'$ holds.

    Now, let us show that $\Equiv'$ is equal to $\EquivCovering$. First, since the left-hand
    and right-hand sides of~\eqref{equ:presentation_inc_1_1} (resp.\
    \eqref{equ:presentation_inc_1_2}) are $\EquivCovering$-equivalent and since
    $\EquivCovering = \Equiv_{\phi_k}$ is a clone congruence by
    Proposition~\ref{prop:p_symbol_composition}, $\Equiv'$ is contained in $\EquivCovering$.
    Conversely, for any  $\PWord, \PWord' \in \P(\TrivialMonoid)$ and $i^\MonoidUnit,
    i_1^\MonoidUnit, i_2^\MonoidUnit \in \SetLetters_\TrivialMonoid$, we have
    \begin{equation}
        \PWord \Conc i_1^\MonoidUnit i_2^\MonoidUnit \Conc \PWord'
        =
        1^\MonoidUnit 2^\MonoidUnit 3^\MonoidUnit \Superposition{
            \PWord,
            1^\MonoidUnit 2^\MonoidUnit \Superposition{i_1^\MonoidUnit, i_2^\MonoidUnit},
            \PWord'}
        \Equiv'
        1^\MonoidUnit 2^\MonoidUnit 3^\MonoidUnit \Superposition{
            \PWord,
            2^\MonoidUnit 1^\MonoidUnit \Superposition{i_1^\MonoidUnit, i_2^\MonoidUnit},
            \PWord'}
        =
        \PWord \Conc i_2^\MonoidUnit i_1^\MonoidUnit \Conc \PWord'
    \end{equation}
    and
    \begin{equation}
        \PWord \Conc \Par{i^\MonoidUnit}^k i^\MonoidUnit \Conc \PWord'
        =
        1^\MonoidUnit 2^\MonoidUnit 3^\MonoidUnit \Superposition{
            \PWord,
            \Par{1^\MonoidUnit}^k 1^\MonoidUnit \Superposition{i^\MonoidUnit},
            \PWord'}
        \Equiv'
        1^\MonoidUnit 2^\MonoidUnit 3^\MonoidUnit \Superposition{
            \PWord,
            \Par{1^\MonoidUnit}^k \Superposition{i^\MonoidUnit},
            \PWord'}
        =
        \PWord \Conc \Par{i^\MonoidUnit}^k \Conc \PWord'.
    \end{equation}
    This shows that for any $\RWord, \RWord' \in \P(\TrivialMonoid)$, $\RWord \Covering
    \RWord'$ implies $\RWord \Equiv' \RWord'$. Since $\EquivCovering$ is the smallest
    equivalence relation containing $\Covering$, $\EquivCovering$ is contained in
    $\Equiv'$. This establishes the statement of the proposition.
\end{Proof}

By Proposition~\ref{prop:presentation_inc}, any $\Inc_k$-algebra is, up to term equivalence,
a monoid $\Par{\Algebra, \Product, \Identity}$ equipped with an additional unary operation
acting identically, where $\Product$ is commutative and satisfies, by spelling
out~\eqref{equ:presentation_inc_2} and simplifying it modulo the background equational
theory $\Equiv_{\RelationSet_{\TrivialMonoid}}$,
\begin{equation}
    \underbrace{x_1 \Product \cdots \Product x_1}_{k + 1 \text{ times}}
    =
    \underbrace{x_1 \Product \cdots \Product x_1}_{k \text{ times}}.
\end{equation}
In particular, $\Inc_1$ is a clone realization of the variety of meet-semilattices admitting
a greatest element (also known as bounded semilattices).

\section{A hierarchy of clones} \label{sec:hierarchy}
We use the construction $\P$ and intersections of the clone congruences $\Equiv_\Sort$,
$\Equiv_{\First_k}$, and $\Equiv_{\First_{k'}^\Reverse}$ introduced in the previous section
to build a hierarchy of clone quotients of $\P(\Monoid)$. Figure~\ref{fig:clone_diagram}
contains the full diagram of the constructed clones. The clones located on the bottom three
lines of the diagram have been constructed and studied in
Section~\ref{sec:construction_quotients}.
\begin{figure}[ht]
    \centering
    \begin{equation*}
        \scalebox{0.9}{
        \begin{tikzpicture}[Centering,xscale=1.15,yscale=0.85,font=\footnotesize]
            \node(P)at(0,2){$\P(\Monoid)$};
            \node(Pill)at(0,0){$\Pill_{k, k'}(\Monoid)$};
            \node(Magn)at(0,-2){$\Magn_{k, k'}(\Monoid)$};
            \node(Stal)at(-2,-2){$\Stal_k(\Monoid)$};
            \node(StalRev)at(2,-2){$\Stal_{k'}^\Reverse(\Monoid)$};
            \node(WInc)at(0,-4){$\WInc(\Monoid)$};
            \node(Arra)at(-2,-4){$\Arra_k(\Monoid)$};
            \node(ArraRev)at(2,-4){$\Arra_{k'}^\Reverse(\Monoid)$};
            \node(Inc)at(0,-6){$\Inc_{\min \Bra{k, k'}}$};
            \node(Inc0)at(0,-8){$\Inc_0$};
            \draw[Surjection](P)--(Pill);
            \draw[Surjection](Pill)--(Magn);
            \draw[Surjection](Pill)--(Stal);
            \draw[Surjection](Pill)--(StalRev);
            \draw[Surjection](Stal)--(Arra);
            \draw[Surjection](StalRev)--(ArraRev);
            \draw[Surjection](Magn)--(Arra);
            \draw[Surjection](Magn)--(ArraRev);
            \draw[Surjection](Stal)--(WInc);
            \draw[Surjection](StalRev)--(WInc);
            \draw[Surjection](Arra)--(Inc);
            \draw[Surjection](ArraRev)--(Inc);
            \draw[Surjection](WInc)--(Inc);
            \draw[Surjection](Inc)--(Inc0);
        \end{tikzpicture}}
    \end{equation*}
    \caption{
        The full diagram of the considered quotients of the clone $\P(\Monoid)$ where
        $\Monoid$ is a monoid and $k, k' \geq 0$. The arrows denote surjective clone
        morphisms.
    }
    \label{fig:clone_diagram}
\end{figure}
The clones constructed in the following sections are clone realizations of varieties
generalizing some special classes of monoids, including regular bands. These structures
allow us to solve the word problem in the corresponding varieties. The algorithms are
described in terms of $\PSymbol$-symbols and are similar to the ones solving the word
problem in idempotent semigroups by using conditional string rewrite
systems~\cite{SS82,NS00}.

In this section, $\Monoid$ is a (finite or infinite) monoid endowed with a total order
relation $\Leq$. To give concrete examples, we shall consider $\Monoid$ as the free monoid
$\Par{A^*, \Conc, \epsilon}$ where $A$ is the alphabet $\{a, b, c\}$ and $\Leq$ is the
lexicographic order on $A^*$ satisfying $a \Leq b \Leq c$.

\subsection{On pigmented magnets} \label{subsec:magn}
By considering the intersection of the clone congruences $\Equiv_{\First_k}$, $k \geq 0$,
and their reversions $\Equiv_{\First_{k'}}^\Reverse$, $k' \geq 0$, we construct a quotient
clone $\Magn_{k, k'}(\Monoid)$ of $\P(\Monoid)$. This clone is studied in detail for the
case $k = 1 = k'$. A description through new combinatorial objects named $\Monoid$-pigmented
magnets is introduced and a finitely equationally generated presentation is described. The
algebras over this clone are generalizations of regular bands. These results are based on
the introduction of a $\PSymbol$-symbol for the underlying equivalence relation.

\subsubsection{Clone construction}
For any parameters $k, k' \geq 0$, let $\Equiv_{k, k'}$ be the clone congruence
$\Equiv_{\First_k} \cap \Equiv_{\First_{k'}}^\Reverse$, and let
\begin{equation}
    \Magn_{k, k'}(\Monoid) := \P(\Monoid) /_{\Equiv_{k, k'}}.
\end{equation}
By Propositions~\ref{prop:congruence_first} and~\ref{prop:reversed_congruence}, $\Magn_{k,
k'}(\Monoid)$ is a well-defined clone, and $\Arra_k(\Monoid)$ and
${\Arra_{k'}^\Reverse(\Monoid)}$ are both isomorphic to quotients of $\Magn_{k, k'}$. Since
for any $0 \leq k \leq k''$ and $0 \leq k' \leq k'''$, $\Equiv_{k'', k'''}$ is a refinement
of $\Equiv_{k, k'}$, $\Magn_{k, k'}(\Monoid)$ is isomorphic to a quotient of $\Magn_{k'',
k'''}(\Monoid)$. Moreover, since $\Equiv_{0, 0}$ is the coarsest clone congruence of
$\P(\Monoid)$, $\Magn_{0, 0}$ is the trivial clone $\TrivialClone$. Besides, the clone
\begin{math}
    \Magn_{k, k'}^\Reverse(\Monoid)
    := \Magn_{k, k'}(\Monoid)^\Reverse
    = \P(\Monoid) /_{\Equiv_{k, k'}^\Reverse}
\end{math}
is by Proposition~\ref{prop:reversed_congruence} isomorphic to $\Magn_{k, k'}(\Monoid)$.
Since the reversion operation on congruences is involutive, it follows that ${\Equiv_{k,
k'}}^\Reverse = \Equiv_{k', k}$ and thus, the clones $\Magn_{k, k'}^\Reverse(\Monoid)$ and
$\Magn_{k', k}(\Monoid)$ are identical and $\Magn_{k, k'}(\Monoid)$ and $\Magn_{k',
k}(\Monoid)$ are isomorphic.

\subsubsection{Equivalence relation}
To lighten the notation, we denote by $\Equiv$ the equivalence relation $\Equiv_{1, 1}$ on
$\P(\Monoid)$. By definition, for any $\PWord, \PWord' \in \P(\Monoid)$, $\PWord \Equiv
\PWord'$ holds if and only if
\begin{math}
    \Par{\First_1\Par{\PWord}, \First^\Reverse_1\Par{\PWord}}
    = \Par{\First_1\Par{\PWord'}, \First^\Reverse_1\Par{\PWord'}}.
\end{math}

In order to obtain properties about the clone $\Magn_{1, 1}(\Monoid)$, we introduce an
alternative equivalence relation $\EquivCovering$ for which we will show that it is equal to
$\Equiv$. Let $\Covering^{(1)}$, $\Covering^{(2)}$, and $\Covering^{(3)}$ be the three
binary relations on $\P(\Monoid)$ satisfying
\begin{equation} \label{equ:covering_relation_magn_1}
    \PWord \Conc \Witness{n}{n}{i^\alpha} \Conc \PWord'
    \enspace \Covering^{(1)} \enspace
    \PWord \Conc \PWord',
\end{equation}
\begin{equation} \label{equ:covering_relation_magn_2}
    \PWord \Conc \Witness{u}{n}{i_1^{\alpha_1}} \Witness{n}{y}{i_2^{\alpha_2}}
        \Conc \PWord'
    \enspace \Covering^{(2)} \enspace
    \PWord \Conc \Witness{n}{y}{i_2^{\alpha_2}} \Witness{u}{n}{i_1^{\alpha_1}}
        \Conc \PWord'
    \qquad
    \text{ where }
    i_1 \ne i_2,
\end{equation}
\begin{equation} \label{equ:covering_relation_magn_3}
    \PWord \Conc \Witness{u}{n}{i^\alpha} \Witness{n}{u}{i^\alpha} \Conc \PWord'
    \enspace \Covering^{(3)} \enspace
    \PWord \Conc \Witness{u}{u}{i^\alpha} \Conc \PWord',
\end{equation}
where $\PWord, \PWord' \in \P(\Monoid)$ and $i^\alpha, i_1^{\alpha_1}, i_2^{\alpha_2} \in
\SetLetters_\Monoid$. Let $\CoveringRT^{(j)}$, $j \in [3]$, be the reflexive and transitive
closure of~$\Covering^{(j)}$, $\Covering$ be the union $\Covering^{(1)} \cup \Covering^{(2)}
\cup \Covering^{(3)}$, and $\EquivCovering$ be the reflexive, symmetric, and transitive
closure of~$\Covering$.

As a side remark, let up emphasize the fact that, despite appearances, $\Covering^{(1)}$,
$\Covering^{(2)}$, and $\Covering^{(3)}$ cannot be studied as rewrite rules of string
rewrite systems~\cite{BO93,BN98,BKV03}. Indeed, since we could have for instance $\PWord
\Covering^{(2)} \PWord'$ but not $\PWord \Conc \QWord \Covering^{(2)} \PWord' \Conc \QWord$
for some $\PWord, \PWord', \QWord \in \P(\Monoid)$, the compatibility with the context
required by string rewrite systems is not satisfied.

\begin{Statement}{Lemma}{lem:covering_relation_magn_monoid_congruence}
    For any monoid $\Monoid$, the equivalence relation $\EquivCovering$ is a monoid
    congruence of the monoid $\Par{\P(\Monoid), \Conc, \PWordEmpty}$.
\end{Statement}
\begin{Proof}
    To prove this statement, since $\EquivCovering$ is the smallest equivalence relation
    containing $\Covering^{(1)}$, $\Covering^{(2)}$, and $\Covering^{(3)}$, it is enough to
    prove that for any $j \in [3]$ and $\QWord, \QWord', \RWord \in \P(\Monoid)$, if $\QWord
    \Covering^{(j)} \QWord'$ then $\QWord \Conc \RWord \EquivCovering \QWord' \Conc \RWord$
    and $\RWord \Conc \QWord \EquivCovering \RWord \Conc \QWord'$.

    Directly from the definitions of $\Covering^{(1)}$, $\Covering^{(2)}$, and
    $\Covering^{(3)}$, for any $j \in [3]$, $\QWord \Covering^{(j)} \QWord'$ implies $\RWord
    \Conc \QWord \Covering^{(j)} \RWord \Conc \QWord'$. This is due to the fact that
    in~\eqref{equ:covering_relation_magn_1}, \eqref{equ:covering_relation_magn_2},
    and~\eqref{equ:covering_relation_magn_3}, adding more letters on the left of $\QWord$
    and $\QWord'$ preserves the required conditions on the left $1$-witnesses of the
    involved $\Monoid$-pigmented words. Moreover, directly from the definitions of
    $\Covering^{(1)}$ and $\Covering^{(3)}$, for any $j \in \{1, 3\}$, $\QWord
    \Covering^{(j)} \QWord'$ implies $\QWord \Conc \RWord \Covering^{(j)} \QWord' \Conc
    \RWord$. This is due to the fact that in~\eqref{equ:covering_relation_magn_1}
    and~\eqref{equ:covering_relation_magn_3}, adding more letters on the right of $\QWord$
    and $\QWord'$ preserves the required conditions on the right $1$-witnesses of the
    involved $\Monoid$-pigmented words. The remaining case to explore happens when $\QWord
    \Covering^{(2)} \QWord'$. In this case, $\QWord$ and $\QWord'$ decompose as $\QWord =
    \PWord \Conc \Witness{u}{n}{i_1^{\alpha_1}} \Witness{n}{y}{i_2^{\alpha_2}} \Conc
    \PWord'$ and $\QWord' = \PWord \Conc \Witness{n}{y}{i_2^{\alpha_2}}
    \Witness{u}{n}{i_1^{\alpha_1}} \Conc \PWord'$ where $\PWord, \PWord' \in \P(\Monoid)$,
    $i_1^{\alpha_1}, i_2^{\alpha_2} \in \SetLetters_\Monoid$, and $i_2 \ne i_1$. As the
    position $\Length(\PWord) + 2$ of $\QWord \Conc \RWord$ is a right $1$-witness if and
    only if there is no $\Monoid$-pigmented letter of value $i_2$ in $\RWord$, we have two
    cases to explore. If this position is a right $1$-witness, then
    \begin{equation}
        \QWord \Conc \RWord
        = \PWord \Conc \Witness{u}{n}{i_1^{\alpha_1}} \Witness{n}{y}{i_2^{\alpha_2}}
            \Conc \PWord'\Conc \RWord
        \Covering^{(2)}
        \PWord \Conc \Witness{n}{y}{i_2^{\alpha_2}}
        \Witness{u}{n}{i_1^{\alpha_1}} \Conc \PWord'\Conc \RWord
        = \QWord' \Conc \RWord.
    \end{equation}
    Otherwise, we have
    \begin{equation}
        \QWord \Conc \RWord
        = \PWord \Conc \Witness{u}{n}{i_1^{\alpha_1}} \Witness{n}{n}{i_2^{\alpha_2}}
            \Conc \PWord'\Conc \RWord
        \Covering^{(1)}
        \PWord \Conc \Witness{u}{n}{i_1^{\alpha_1}} \Conc \PWord'\Conc \RWord
    \end{equation}
    and
    \begin{equation}
        \QWord' \Conc \RWord
        = \PWord \Conc \Witness{n}{n}{i_2^{\alpha_2}} \Witness{u}{n}{i_1^{\alpha_1}}
            \Conc \PWord'\Conc \RWord
        \Covering^{(1)}
        \PWord \Conc \Witness{u}{n}{i_1^{\alpha_1}} \Conc \PWord'\Conc \RWord.
    \end{equation}
    This shows that $\QWord \Conc \RWord \EquivCovering \QWord' \Conc \RWord$.
\end{Proof}

\begin{Statement}{Lemma}{lem:covering_relation_magn_first}
    For any monoid $\Monoid$ and any $\PWord \in \P(\Monoid)$,
    \begin{equation}
        \PWord \EquivCovering \First_1(\PWord) \Conc \First_1^\Reverse(\PWord).
    \end{equation}
\end{Statement}
\begin{Proof}
    Let us first show that $\PWord \EquivCovering \PWord \Conc \PWord$ by induction on $\ell
    := \Length(\PWord)$. If $\ell = 0$, then $\PWord = \PWordEmpty$ and since $\PWord \Conc
    \PWord = \PWordEmpty$, the stated property holds. Assume now that $\ell \geq 1$. In this
    case, $\PWord$ decomposes as $\PWord = \Witness{y}{u}{i^\alpha} \Conc \PWord'$ where
    $i^\alpha \in \SetLetters_\Monoid$ and $\PWord' \in \P(\Monoid)$. We have now
    \begin{math}
        \PWord \Conc \PWord
        = \Witness{y}{n}{i^\alpha} \Conc \PWord' \Conc
            \Witness{n}{u}{i^\alpha} \Conc \PWord'
    \end{math}
    and two cases to explore depending on whether the position $\ell + 1$ in $\PWord \Conc
    \PWord$ is a right $1$-witness.
    \begin{enumerate}[label={\it (\Roman*)}]
        \item If it is the case, then
        \begin{math}
            \PWord \Conc \PWord
            = \Witness{y}{n}{i^\alpha} \Conc \PWord' \Conc
            \Witness{n}{y}{i^\alpha} \Conc \PWord'.
        \end{math}
        Since there is no occurrence of any $\Monoid$-pigmented letter having $i$ as value
        in $\PWord'$, and additionally, there is no position $j \in [\ell]$ in $\PWord \Conc
        \PWord$ which is a right $1$-witness, we have
        \begin{math}
            \Witness{y}{n}{i^\alpha} \Conc \PWord' \Conc
            \Witness{n}{y}{i^\alpha} \Conc \PWord'
            \Covering^{(2)} \dots \Covering^{(2)}
            \Witness{y}{n}{i^\alpha} \Witness{n}{y}{i^\alpha} \Conc \PWord' \Conc \PWord'.
        \end{math}
        \item Otherwise,
        \begin{math}
            \PWord \Conc \PWord
            = \Witness{y}{n}{i^\alpha} \Conc \PWord' \Conc
            \Witness{n}{n}{i^\alpha} \Conc \PWord'.
        \end{math}
        Since there are occurrences of letters having $i$ as value in $\PWord'$, we have
        \begin{math}
            \Witness{y}{n}{i^\alpha} \Conc \PWord' \Conc
                \Witness{n}{n}{i^\alpha} \Conc \PWord'
            \Covering^{(1)}
           \Witness{y}{n}{i^\alpha} \Conc \PWord' \Conc \PWord'
        \end{math}
        and
        \begin{math}
            \Witness{y}{n}{i^\alpha} \Witness{n}{n}{i^\alpha} \Conc \PWord'
                 \Conc \PWord'
            \Covering^{(1)}
            \Witness{y}{n}{i^\alpha} \Conc \PWord' \Conc \PWord'.
        \end{math}
    \end{enumerate}
    In both cases, by induction hypothesis and by using the fact that by
    Lemma~\ref{lem:covering_relation_magn_monoid_congruence}, $\EquivCovering$ is a monoid
    congruence, we obtain
    \begin{math}
        \PWord \Conc \PWord
        \EquivCovering \Witness{y}{n}{i^\alpha} \Witness{n}{u}{i^\alpha} \Conc \PWord'
        \Conc \PWord'
        \EquivCovering \Witness{y}{n}{i^\alpha} \Witness{n}{u}{i^\alpha} \Conc \PWord'.
    \end{math}
    Finally, since
    \begin{math}
        \Witness{y}{n}{i^\alpha} \Witness{n}{u}{i^\alpha} \Conc \PWord'
        \Covering^{(3)}
        \Witness{y}{u}{i^\alpha} \Conc \PWord'
        = \PWord,
    \end{math}
    the stated property is established.

    Let us now show that
    \begin{math}
        \PWord \Conc \PWord
        \EquivCovering
        \First_1(\PWord) \Conc \First_1^\Reverse(\PWord).
    \end{math}
    By assuming that $\PWord$ can be written as $\PWord = i_1^{\alpha_1} \dots
    i_\ell^{\alpha_\ell}$, there exists a unique pair $\Par{r_1 \dots r_k, s_1 \dots s_k}$
    of subwords of $1 \dots \ell$ such that
    \begin{math}
        \First_1(\PWord) = i_{r_1}^{\alpha_{r_1}} \dots i_{r_k}^{\alpha_{r_k}}
    \end{math}
    and
    \begin{math}
        \First_1^\Reverse(\PWord) = i_{s_1}^{\alpha_{s_1}} \dots i_{s_k}^{\alpha_{s_k}}.
    \end{math}
    Therefore, we have
    \begin{math}
        \Witness{y}{u}{i_{r_1}^{\alpha_{r_1}}} \Conc \PWord_1 \Conc \
        \dots \ \Conc \Witness{y}{u}{i_{r_k}^{\alpha_{r_k}}} \Conc \PWord_k
        = \PWord =
        \PWord'_k \Conc \Witness{u}{y}{i_{s_1}^{\alpha_{s_1}}} \Conc \
        \dots \ \Conc \PWord'_1 \Conc \Witness{u}{y}{i_{s_k}^{\alpha_{s_k}}}
    \end{math}
    where $\PWord_1, \dots, \PWord_k, \PWord'_k, \dots, \PWord'_1 \in \P(\Monoid)$. Hence,
    \begin{equation}
        \PWord \Conc \PWord =
        \Witness{y}{u}{i_{r_1}^{\alpha_{r_1}}} \Conc \PWord_1 \Conc \
        \dots \ \Conc \Witness{y}{u}{i_{r_k}^{\alpha_{r_k}}} \Conc \PWord_k
        \Conc
        \PWord'_k \Conc \Witness{u}{y}{i_{s_1}^{\alpha_{s_1}}} \Conc \
        \dots \ \Conc \PWord'_1 \Conc \Witness{u}{y}{i_{s_k}^{\alpha_{s_k}}},
    \end{equation}
    and since the positions in $\PWord \Conc \PWord$ of the letters of its factors
    $\PWord_1$, \dots, $\PWord_k$, $\PWord'_k$, \dots, $\PWord'_1$ are neither left
    $1$-witnesses nor right $1$-witnesses, we have
    \begin{equation}
        \PWord \Conc \PWord
        \Covering^{(1)} \dots \Covering^{(1)}
        \Witness{y}{u}{i_{r_1}^{\alpha_{r_1}}} \dots
        \Witness{y}{u}{i_{r_k}^{\alpha_{r_k}}}
        \Witness{u}{y}{i_{s_1}^{\alpha_{s_1}}} \dots
        \Witness{u}{y}{i_{s_k}^{\alpha_{s_k}}}
        = \First_1(\PWord) \Conc \First_1^\Reverse(\PWord).
    \end{equation}

    By putting these $\EquivCovering$-equivalences together, we obtain
    \begin{math}
        \PWord \EquivCovering \PWord \Conc \PWord
        \EquivCovering
        \First_1(\PWord) \Conc \First_1^\Reverse(\PWord)
    \end{math}
    establishing the stated $\EquivCovering$-equivalence.
\end{Proof}

\begin{Statement}{Proposition}{prop:covering_relation_magn_equivalence}
    For any monoid $\Monoid$, the binary relations $\Equiv$ and $\EquivCovering$ on
    $\P(\Monoid)$ are equal.
\end{Statement}
\begin{Proof}
    First, observe that for any $j \in [3]$ and any $\PWord, \PWord' \in \P(\Monoid)$, if
    $\PWord \Covering^{(j)} \PWord'$, then $\First_1\Par{\PWord} = \First_1\Par{\PWord'}$
    and $\First_1^\Reverse\Par{\PWord} = \First_1^\Reverse\Par{\PWord'}$. Hence, and since
    $\EquivCovering$ is the smallest equivalence relation containing $\Covering^{(1)}$,
    $\Covering^{(2)}$, and $\Covering^{(3)}$, we have that $\EquivCovering$ is contained in
    $\Equiv$. Conversely, for any $\PWord, \PWord' \in \P(\Monoid)$ such that $\PWord \Equiv
    \PWord'$, we have $\First_1\Par{\PWord} = \First_1\Par{\PWord'}$ and
    $\First_1^\Reverse\Par{\PWord} = \First_1^\Reverse\Par{\PWord'}$. By
    Lemma~\ref{lem:covering_relation_magn_first},
    \begin{math}
        \PWord \EquivCovering \First_1(\PWord) \Conc \First_1^\Reverse(\PWord)
        = \First_1\Par{\PWord'} \Conc \First_1^\Reverse\Par{\PWord'}
        \EquivCovering \PWord'.
    \end{math}
    For this reason, we have $\PWord \EquivCovering \PWord'$, showing that $\Equiv$ is
    contained in~$\EquivCovering$.
\end{Proof}

\subsubsection{$\PSymbol$-symbol algorithm}
With the aim of describing $\Magn_{1, 1}(\Monoid)$, we propose now a $\PSymbol$-symbol for
$\Equiv$.

\begin{Statement}{Lemma}{lem:covering_relation_magn_orders}
    For any monoid $\Monoid$, the binary relation $\CoveringRT^{(j)}$, $j \in [3]$, is a
    partial order relation on $\P(\Monoid)$. Moreover, for any $\PWord \in \P(\Monoid)$,
    there is exactly one maximal element $\QWord$ of the poset $\Par{\P(\Monoid),
    \CoveringRT^{(j)}}$ such that $\PWord \CoveringRT^{(j)} \QWord$.
\end{Statement}
\begin{Proof}
    Let us consider each binary relation $\CoveringRT^{(j)}$, $j \in [3]$, one by one.
    \begin{enumerate}[label={\it (\Roman*)}]
        \item For any $\PWord, \PWord' \in \P(\Monoid)$, we have $\PWord \CoveringRT^{(1)}
        \PWord'$ if and only if $\PWord'$ can be obtained from $\PWord$ by deleting some
        (possibly none) $\Monoid$-pigmented letters whose positions are neither left
        $1$-witnesses nor right $1$-witnesses. This implies immediately the properties of
        the statement of lemma for~$\CoveringRT^{(1)}$.
        \item For any $\PWord, \PWord' \in \P(\Monoid)$, if $\PWord \Covering^{(2)}
        \PWord'$, then by denoting by $\tau(\PWord)$ (resp.\ $\tau\Par{\PWord'}$) the sum of
        the positions of $\PWord$ (resp.\ $\PWord'$) of the $\Monoid$-pigmented letters
        which are right $1$-witnesses, we have $\tau(\PWord) = \tau\Par{\PWord'} + 1$. Since
        $\CoveringRT^{(2)}$ is the reflexive and transitive closure of $\Covering^{(2)}$,
        this shows that $\CoveringRT^{(2)}$ is antisymmetric. The second property is a
        consequence of the fact that for any $\PWord, \PWord', \PWord'' \in \P(\Monoid)$, if
        $\PWord' \ne \PWord''$, $\PWord \Covering^{(2)} \PWord'$, and $\PWord
        \Covering^{(2)} \PWord''$, then there exists $\PWord''' \in \P(\Monoid)$ such that
        $\PWord' \CoveringRT^{(2)} \PWord'''$ and $\PWord'' \CoveringRT^{(2)} \PWord'''$.
        This property is due to the fact that for any $\RWord, \RWord' \in \P(\Monoid)$ and
        $i_1^{\alpha_1}, i_2^{\alpha_2}, i_3^{\alpha_3} \in \SetLetters_\Monoid$, it is not
        possible to have both
        \begin{math}
            \RWord \Conc \Witness{u}{n}{i_1^{\alpha_1}} \Witness{n}{y}{i_2^{\alpha_2}}
                \Witness{u}{u}{i_3^{\alpha_3}} \Conc \RWord'
            \Covering^{(2)}
            \RWord \Conc \Witness{n}{y}{i_2^{\alpha_2}} \Witness{u}{n}{i_1^{\alpha_1}}
                \Witness{u}{u}{i_3^{\alpha_3}} \Conc \RWord'
        \end{math}
        and
        \begin{math}
            \RWord \Conc \Witness{u}{u}{i_1^{\alpha_1}} \Witness{u}{n}{i_2^{\alpha_2}}
                \Witness{n}{y}{i_3^{\alpha_3}} \Conc \RWord'
            \Covering^{(2)}
            \RWord \Conc \Witness{u}{u}{i_1^{\alpha_1}} \Witness{n}{y}{i_3^{\alpha_3}}
                \Witness{u}{n}{i_2^{\alpha_2}} \Conc \RWord'.
        \end{math}
        Indeed, these two properties would lead to the fact that the position
        $\Length(\RWord) + 2$ of $\RWord \Conc i_1^{\alpha_1} i_2^{\alpha_2} i_3^{\alpha_3}
        \Conc \RWord'$ is a right $1$-witness and, at the same time, is not a right
        $1$-witness. Therefore, the swapping of two consecutive positions that led from
        $\PWord$ to $\PWord'$ must appear in two consecutive positions compared to those
        swapped by the transition from $\PWord$ to $\PWord''$. Consequently, $\PWord'''$ can
        be constructed by first swapping the first two places and then the two other ones of
        $\PWord$.
        \item For any $\PWord, \PWord' \in \P(\Monoid)$, we have $\PWord \CoveringRT^{(3)}
        \PWord'$ if and only if $\PWord'$ can be obtained from $\PWord$ by deleting some
        (possibly none) $\Monoid$-pigmented letters which have a same $\Monoid$-pigmented
        letter as neighbor. In the same way as the first case, this implies immediately the
        properties of the statement of lemma for~$\CoveringRT^{(3)}$.
    \end{enumerate}
\end{Proof}

Let, for any $j \in [3]$, $\Max^{(j)} : \P(\Monoid) \to \P(\Monoid)$ be the map such that
for any $\PWord \in \P(\Monoid)$, $\PWord \Max^{(j)}$ is the maximal element of the poset
$\Par{\P(\Monoid), \CoveringRT^{(j)}}$ comparable with $\PWord$. By
Lemma~\ref{lem:covering_relation_magn_orders}, this map is well-defined.

Let $\PSymbol_{\Equiv} : \P(\Monoid) \to \P(\Monoid)$ be the map defined for any $\PWord \in
\P(\Monoid)$ by
\begin{equation}
    \PSymbol_{\Equiv}(\PWord) := \PWord \Max^{(1)} \Max^{(2)} \Max^{(3)}.
\end{equation}
For instance, in $\P\Par{A^*}$, where $A^*$ is the free monoid over $\Bra{a, b, c}$, we have
\begin{align}
    \PSymbol_{\Equiv}\Par{
        \Witness{y}{n}{2^\epsilon} \Witness{y}{n}{1^b} \Witness{n}{y}{2^\epsilon}
        \Witness{y}{n}{3^a} \Witness{n}{n}{1^{ba}} \Witness{n}{y}{1^b}
        \Witness{n}{y}{3^\epsilon}}
    & =
    \Witness{y}{n}{2^\epsilon} \Witness{y}{n}{1^b} \Witness{n}{y}{2^\epsilon}
    \Witness{y}{n}{3^a} \Witness{n}{n}{1^{ba}} \Witness{n}{y}{1^b}
    \Witness{n}{y}{3^\epsilon} \Max^{(1)} \Max^{(2)} \Max^{(3)}
    \\
    & =
    \Witness{y}{n}{2^\epsilon} \Witness{y}{n}{1^b} \Witness{n}{y}{2^\epsilon}
    \Witness{y}{n}{3^a} \Witness{n}{y}{1^b} \Witness{n}{y}{3^\epsilon} \Max^{(2)} \Max^{(3)}
    \notag \\
    & =
    \Witness{y}{n}{2^\epsilon} \Witness{n}{y}{2^\epsilon} \Witness{y}{n}{1^b}
    \Witness{n}{y}{1^b} \Witness{y}{n}{3^a} \Witness{n}{y}{3^\epsilon} \Max^{(3)}
    \notag \\
    & =
    \Witness{y}{y}{2^\epsilon} \Witness{y}{y}{1^b} \Witness{y}{n}{3^a}
    \Witness{y}{y}{3^\epsilon}
    \notag
\end{align}
and
\begin{align}
    \PSymbol_{\Equiv} & \Par{
        \Witness{y}{n}{4^a} \Witness{y}{n}{2^b} \Witness{y}{n}{1^c} \Witness{n}{y}{1^c}
        \Witness{n}{n}{4^b} \Witness{y}{n}{3^b} \Witness{n}{y}{3^a} \Witness{n}{n}{2^a}
        \Witness{n}{n}{2^a} \Witness{n}{n}{4^a} \Witness{n}{n}{2^c} \Witness{n}{y}{4^a}
        \Witness{n}{y}{2^c}}
    \\
    & \qquad \qquad =
    \Witness{y}{n}{4^a} \Witness{y}{n}{2^b} \Witness{y}{n}{1^c} \Witness{n}{y}{1^c}
    \Witness{n}{n}{4^b} \Witness{y}{n}{3^b} \Witness{n}{y}{3^a} \Witness{n}{n}{2^a}
    \Witness{n}{n}{2^a} \Witness{n}{n}{4^a} \Witness{n}{n}{2^c} \Witness{n}{y}{4^a}
    \Witness{n}{y}{2^c}
    \Max^{(1)} \Max^{(2)} \Max^{(3)}
    \notag \\
    & \qquad \qquad =
    \Witness{y}{n}{4^a} \Witness{y}{n}{2^b} \Witness{y}{n}{1^c} \Witness{n}{y}{1^c}
    \Witness{y}{n}{3^b} \Witness{n}{y}{3^a} \Witness{n}{y}{4^a} \Witness{n}{y}{2^c}
    \Max^{(2)} \Max^{(3)}
    \notag \\
    & \qquad \qquad =
    \Witness{y}{n}{4^a} \Witness{y}{n}{2^b} \Witness{y}{n}{1^c} \Witness{n}{y}{1^c}
    \Witness{y}{n}{3^b} \Witness{n}{y}{3^a} \Witness{n}{y}{4^a} \Witness{n}{y}{2^c}
    \Max^{(3)}
    \notag \\
    & \qquad \qquad =
    \Witness{y}{n}{4^a} \Witness{y}{n}{2^b} \Witness{y}{y}{1^c} \Witness{y}{n}{3^b}
    \Witness{n}{y}{3^a} \Witness{n}{y}{4^a} \Witness{n}{y}{2^c}.
    \notag
\end{align}
Let us emphasize the fact that the maps $\Max^{(1)}$, $\Max^{(2)}$, and $\Max^{(3)}$ do not
commute. Indeed, we have for instance
\begin{equation}
    \PSymbol_{\Equiv} \Par{
        \Witness{y}{n}{1^\epsilon} \Witness{n}{n}{1^a} \Witness{n}{y}{1^\epsilon}}
        = \Witness{y}{y}{1^\epsilon}
        \enspace \ne \enspace \Witness{y}{n}{1^\epsilon} \Witness{n}{y}{1^\epsilon}
        = \Witness{y}{n}{1^\epsilon} \Witness{n}{n}{1^a} \Witness{n}{y}{1^\epsilon}
        \Max^{(2)} \Max^{(3)} \Max^{(1)}.
\end{equation}

\begin{Statement}{Lemma}{lem:p_symbol_magn_equiv_p_symbol}
    For any monoid $\Monoid$ and any $\PWord \in \P(\Monoid)$, $\PWord \EquivCovering
    \PSymbol_{\Equiv}(\PWord)$.
\end{Statement}
\begin{Proof}
    First, since for any $j \in [3]$, $\EquivCovering$ contains $\CoveringRT^{(j)}$, we have
    $\PWord \EquivCovering \PWord \Max^{(j)}$. Moreover, as $\PSymbol_{\Equiv}$ is by
    definition the map composition $\Max^{(3)} \circ \Max^{(2)} \circ \Max^{(1)}$,
    $\PSymbol_{\Equiv}(\PWord)$ is $\EquivCovering$-equivalent to $\PWord$.
\end{Proof}

\begin{Statement}{Lemma}{lem:covering_relation_magn_p_symbol}
    For any monoid $\Monoid$ and any $\PWord, \PWord' \in \P(\Monoid)$, $\PWord
    \EquivCovering \PWord'$ implies $\PSymbol_{\Equiv}(\PWord) =
    \PSymbol_{\Equiv}\Par{\PWord'}$.
\end{Statement}
\begin{Proof}
    Let us show that $\PWord \Covering \PWord'$ for $\PWord, \PWord' \in \P(\Monoid)$
    entails $\PSymbol_{\Equiv}(\PWord) = \PSymbol_{\Equiv}\Par{\PWord'}$. Since the
    equivalence relation $\EquivCovering$ is generated by $\Covering$, this will entail the
    statement of the lemma. We are going to consider the following three cases depending
    whether $\PWord \Covering^{(1)} \PWord'$, $\PWord \Covering^{(2)} \PWord'$, or $\PWord
    \Covering^{(3)} \PWord'$.
    \begin{enumerate}[label={\it (\Roman*)}]
        \item Assume that $\PWord \Covering^{(1)} \PWord'$. By
        Lemma~\ref{lem:covering_relation_magn_orders}, $\PWord \Max^{(1)} = \PWord'
        \Max^{(1)}$. Therefore, by definition of $\PSymbol_{\Equiv}$,
        $\PSymbol_{\Equiv}(\PWord) = \PSymbol_{\Equiv}\Par{\PWord'}$.
        \item Assume that $\PWord \Covering^{(2)} \PWord'$. Hence, $\PWord$ and $\PWord'$
        decompose as $\PWord = \QWord \Conc \Witness{u}{n}{i_1^{\alpha_1}}
        \Witness{n}{y}{i_2^{\alpha_2}} \Conc \RWord$ and $\PWord' = \QWord \Conc
        \Witness{n}{y}{i_2^{\alpha_2}} \Witness{u}{n}{i_1^{\alpha_1}} \Conc \RWord$ where
        $\QWord, \RWord \in \P(\Monoid)$, $i_1^{\alpha_1}, i_2^{\alpha_2} \in
        \SetLetters_\Monoid$, and $i_1 \ne i_2$. If the letter at position $\Length(\QWord)
        + 1$ of $\PWord$ is not a left $1$-witness, then the letter at position
        $\Length(\QWord) + 2$ of $\PWord'$ is not a left $1$-witness and
        \begin{equation}
            \PWord \Max^{(1)}
            = \Par{\QWord \Max^{(1)}} \Conc \Witness{n}{y}{i_2^{\alpha_2}}
                \Conc \Par{\RWord \Max^{(1)}}
            = \PWord' \Max^{(1)}.
        \end{equation}
        Therefore, $\PWord \Max^{(1)} = \PWord' \Max^{(1)}$ and $\PSymbol_{\Equiv}(\PWord) =
        \PSymbol_{\Equiv}\Par{\PWord'}$. Otherwise, when the letter at position
        $\Length(\QWord) + 1$ of $\PWord$ is a left $1$-witness, the letter at position
        $\Length(\QWord) + 2$ of $\PWord'$ is also a left $1$-witness and we have
        \begin{equation}
            \PWord \Max^{(1)}
            = \Par{\QWord \Max^{(1)}} \Conc \Witness{y}{n}{i_1^{\alpha_1}}
                \Witness{n}{y}{i_2^{\alpha_2}} \Conc \Par{\RWord \Max^{(1)}}
            \Covering^{(2)}
            \Par{\QWord \Max^{(1)}} \Conc \Witness{n}{y}{i_2^{\alpha_2}}
                \Witness{y}{n}{i_1^{\alpha_1}} \Conc \Par{\RWord \Max^{(1)}}
            = \PWord' \Max^{(1)}.
        \end{equation}
        By Lemma~\ref{lem:covering_relation_magn_orders}, $\PWord \Max^{(1)} \Max^{(2)} =
        \PWord' \Max^{(1)} \Max^{(2)}$. Therefore, by definition of $\PSymbol_{\Equiv}$,
        $\PSymbol_{\Equiv}(\PWord) = \PSymbol_{\Equiv}\Par{\PWord'}$.
        \item Assume that $\PWord \Covering^{(3)} \PWord'$. Hence, $\PWord$ and $\PWord'$
        decompose as $\PWord = \QWord \Conc \Witness{u}{u}{i^\alpha}
        \Witness{u}{u}{i^\alpha} \Conc \RWord$ and $\PWord' = \QWord \Conc
        \Witness{u}{u}{i^\alpha} \Conc \RWord$ where $\QWord, \RWord \in \P(\Monoid)$ and
        $i^\alpha \in \SetLetters_\Monoid$. If the letters at positions $\Length(\QWord) +
        1$ and $\Length(\QWord) + 2$ of $\PWord$ are neither left $1$-witnesses nor right
        $1$-witnesses, then the letter at position $\Length(\QWord) + 1$ of $\PWord'$ is
        neither a left $1$-witnesses nor a right $1$-witness and
        \begin{equation}
             \PWord \Max^{(1)}
            = \Par{\QWord \Max^{(1)}} \Conc \Par{\RWord \Max^{(1)}}
            = \PWord' \Max^{(1)}.
        \end{equation}
        Therefore, $\PWord \Max^{(1)} = \PWord' \Max^{(1)}$ and $\PSymbol_{\Equiv}(\PWord) =
        \PSymbol_{\Equiv}\Par{\PWord'}$. Otherwise, if there is exactly one position among
        $\Length(\QWord) + 1$ and $\Length(\QWord) + 2$ of $\PWord$ which is neither a left
        $1$-witness nor a right $1$-witness, then, since the letter at position
        $\Length(\QWord) + 1$ of $\PWord$ cannot be a right $1$-witness and the other one
        cannot be a left $1$-witness, we have
        \begin{equation}
            \PWord \Max^{(1)}
            = \Par{\QWord \Max^{(1)}} \Conc \Witness{u}{u}{i^\alpha} \Conc \Par{\RWord
            \Max^{(1)}} = \PWord' \Max^{(1)}.
        \end{equation}
        Therefore, $\PWord \Max^{(1)} = \PWord' \Max^{(1)}$ and $\PSymbol_{\Equiv}(\PWord) =
        \PSymbol_{\Equiv}\Par{\PWord'}$. The last possibility happens when the letter at
        position $\Length(\QWord) + 1$ of $\PWord$ is a left $1$-witness and the letter at
        position $\Length(\QWord) + 2$ of $\PWord$ is a right $1$-witness. In this case,
        \begin{align}
            \PWord \Max^{(1)} \Max^{(2)}
            & =
            \Par{\QWord \Max^{(1)} \Max^{(2)}} \Conc \Witness{y}{n}{i^\alpha}
                \Witness{n}{y}{i^\alpha} \Conc \Par{\RWord \Max^{(1)} \Max^{(2)}}
            \\
            & \Covering^{(3)}
            \Par{\QWord \Max^{(1)} \Max^{(2)}} \Conc \Witness{y}{y}{i^\alpha}
                \Conc \Par{\RWord \Max^{(1)} \Max^{(2)}}
            = \PWord' \Max^{(1)} \Max^{(2)}.
            \notag
        \end{align}
        By Lemma~\ref{lem:covering_relation_magn_orders},
        \begin{math}
            \PWord \Max^{(1)} \Max^{(2)} \Max^{(3)} = \PWord' \Max^{(1)} \Max^{(2)}
            \Max^{(3)}. \end{math}
        Therefore, by definition of $\PSymbol_{\Equiv}$, $\PSymbol_{\Equiv}(\PWord) =
        \PSymbol_{\Equiv}\Par{\PWord'}$.
    \end{enumerate}
\end{Proof}

By Proposition~\ref{prop:covering_relation_magn_equivalence} and
Lemmas~\ref{lem:p_symbol_magn_equiv_p_symbol} and~\ref{lem:covering_relation_magn_p_symbol},
$\PSymbol_{\Equiv}$ is a $\PSymbol$-symbol for $\Equiv$.

\subsubsection{Description}
An \Def{$\Monoid$-pigmented magnet} (or simply \Def{pigmented magnet} when the context is
clear) of arity $n \geq 0$ is an $\Monoid$-pigmented word $\PWord$ of arity $n$ which is a
maximal element at the same time in the posets $\Par{\P(\Monoid), \CoveringRT^{(1)}}$,
$\Par{\P(\Monoid), \CoveringRT^{(2)}}$, and $\Par{\P(\Monoid), \CoveringRT^{(3)}}$. For
instance, in $\P\Par{A^*}$, where $A^*$ is the free monoid over $\Bra{a, b}$,
\begin{equation}
    \Witness{y}{n}{1^{aa}} \Witness{n}{n}{1^{b}} \Witness{y}{y}{2^{ab}}
    \Witness{n}{y}{1^{b}}
    \quad \text{and} \quad
    \Witness{y}{n}{2^{ba}} \Witness{y}{n}{3^{\epsilon}} \Witness{n}{y}{2^{ab}}
    \Witness{y}{y}{1^{a}} \Witness{n}{y}{3^{ba}}
\end{equation}
are not $A^*$-pigmented magnets. In contrast,
\begin{equation}
    \Witness{y}{y}{3^{b}} \Witness{y}{n}{2^{ba}} \Witness{y}{y}{4^{ba}}
    \Witness{y}{y}{1^{a}} \Witness{n}{y}{2^{ab}}
    \quad \text{and} \quad
    \Witness{y}{n}{2^{bb}} \Witness{y}{n}{1^{a}} \Witness{n}{y}{1^{aa}}
    \Witness{n}{y}{2^{bb}}
\end{equation}
are $A^*$-pigmented magnets.

\begin{Statement}{Lemma}{lem:p_symbol_magn_magnets}
    For any monoid $\Monoid$ and any $\PWord \in \P(\Monoid)$, $\PSymbol_{\Equiv}(\PWord)$
    is an $\Monoid$-pigmented magnet.
\end{Statement}
\begin{Proof}
    Let $\PWord_1 := \PWord \Max^{(1)}$, $\PWord_2 := \PWord_1 \Max^{(2)}$, and $\PWord_3 :=
    \PWord_2 \Max^{(3)}$. By definition of $\PSymbol_{\Equiv}$, $\PWord_3 =
    \PSymbol_{\Equiv}\Par{\PWord}$. Let us show that $\PWord_3$ is a maximal element w.r.t.\
    the partial order relations $\CoveringRT^{(1)}$, $\CoveringRT^{(2)}$, and
    $\CoveringRT^{(3)}$ at the same time.

    Observe first that $\PWord_1$ does not contain any letters $i^\alpha \in
    \SetLetters_\Monoid$ of the form $\Witness{n}{n}{i^\alpha}$. The construction of
    $\PWord_2$ from $\PWord_1$ iteratively swaps adjacent letters
    $\Witness{y}{n}{i_1^{\alpha_1}} \Witness{n}{y}{i_2^{\alpha_2}}$ with $i_1, i_2 \in
    \SetLetters_\Monoid$ and $i_1 \ne i_2$ into $\Witness{n}{y}{i_2^{\alpha_2}}
    \Witness{y}{n}{i_1^{\alpha_1}}$ thereby not modifying the left or right $1$-witnesses
    states of any other position of the word. Therefore, $\PWord_2$ neither has any letters
    $i^\alpha \in \SetLetters_\Monoid$ of the form $\Witness{n}{n}{i^\alpha}$ nor any
    adjacent letters $i_1, i_2 \in \SetLetters_\Monoid$ of the form
    $\Witness{y}{n}{i_1^{\alpha_1}} \Witness{n}{y}{i_2^{\alpha_2}}$ with $i_1 \ne i_2$.
    Thus, $\PWord_2$ is at the same time maximal w.r.t.\ $\CoveringRT^{(1)}$ and
    $\CoveringRT^{(2)}$. Since $\PWord_2$ does not contain letters $i^\alpha \in
    \SetLetters_\Monoid$ of the form $\Witness{n}{n}{i^\alpha}$, the construction of
    $\PWord_3$ from $\PWord_2$ iteratively compresses adjacent identical letters $i^\alpha
    \in \SetLetters_\Monoid$ of the form $\Witness{y}{n}{i^\alpha} \Witness{n}{y}{i^\alpha}$
    into $\Witness{y}{y}{i^\alpha}$, which again does not influence the left or right
    $1$-witness state of any other letter in the word. Therefore, $\PWord_3$ does not
    contain letters $i^\alpha \in \SetLetters_\Monoid$ of the form
    $\Witness{n}{n}{i^\alpha}$ (is $\CoveringRT^{(1)}$-maximal), nor adjacent letters $i_1,
    i_2 \in \SetLetters_\Monoid$ of the form $\Witness{y}{n}{i_1^{\alpha_1}}
    \Witness{n}{y}{i_2^{\alpha_2}}$ with $i_1 \ne i_2$ (is $\CoveringRT^{(2)}$-maximal).
    Therefore, $\PWord_3$ is maximal w.r.t.\ all three partial orders.
\end{Proof}

\begin{Statement}{Theorem}{thm:p_symbol_magn}
    For any monoid $\Monoid$, $\PSymbol_{\Equiv}$ is a $\PSymbol$-symbol for $\Equiv$ and
    $\PSymbol_{\Equiv}(\P(\Monoid))$ is the set of $\Monoid$-pigmented magnets. Moreover,
    the graded set $\Magn_{1, 1}(\Monoid)$ is isomorphic to the graded set of
    $\Monoid$-pigmented magnets.
\end{Statement}
\begin{Proof}
    By Proposition~\ref{prop:covering_relation_magn_equivalence} and
    Lemmas~\ref{lem:p_symbol_magn_equiv_p_symbol}
    and~\ref{lem:covering_relation_magn_p_symbol}, $\PSymbol_{\Equiv}$ is a
    $\PSymbol$-symbol for $\Equiv$. Therefore, $\PSymbol_{\Equiv}$ is idempotent, which
    implies together with Lemma~\ref{lem:p_symbol_magn_magnets} that
    $\PSymbol_{\Equiv}(\P(\Monoid))$ is the set of $\Monoid$-pigmented magnets. The last
    part of the statement is a direct implication of
    Proposition~\ref{prop:clone_realization_p_symbol} and the fact that $\PSymbol_{\Equiv}$
    is, as we have just shown, a $\PSymbol$-symbol for~$\Equiv$.
\end{Proof}

By Proposition~\ref{prop:clone_realization_p_symbol} and Theorem~\ref{thm:p_symbol_magn},
$\Magn_{1, 1}(\Monoid)$ can be seen as a clone on $\Monoid$-pigmented magnets with
superposition maps satisfying~\eqref{equ:superposition_p_symbol}. For instance, in
$\Magn_{1, 1}\Par{A^*}$, where $A^*$ is the free monoid over $\Bra{a, b}$, we have, up to
isomorphism,
\begin{align}
    1^a 1^b 4^b 3^{ba} 2^b & \Superposition{
        3^b 3^a,
        1^\epsilon 1^{ba} 3^\epsilon 2^\epsilon 2^{ab} 3^{ab},
        1^\epsilon 1^a,
        2^\epsilon 3^a 3^b 1^a
    }
    \\
    & =
    \PSymbol_{\Equiv}\Par{
        \Witness{y}{n}{3^{ab}}
        \Witness{n}{n}{3^{aa}}
        \Witness{n}{n}{3^{bb}}
        \Witness{n}{n}{3^{ba}}
        \Witness{y}{n}{2^b}
        \Witness{n}{n}{3^{ba}}
        \Witness{n}{n}{3^{bb}}
        \Witness{y}{n}{1^{ba}}
        \Witness{n}{n}{1^{ba}}
        \Witness{n}{n}{1^{baa}}
        \Witness{n}{n}{1^b}
        \Witness{n}{y}{1^{bba}}
        \Witness{n}{n}{3^b}
        \Witness{n}{n}{2^b}
        \Witness{n}{y}{2^{bab}}
        \Witness{n}{y}{3^{bab}}
    }
    \notag \\
    & =
    \Witness{y}{n}{3^{ab}}
    \Witness{y}{n}{2^b}
    \Witness{y}{n}{1^{ba}}
    \Witness{n}{y}{1^{bba}}
    \Witness{n}{y}{2^{bab}}
    \Witness{n}{y}{3^{bab}}.
    \notag
\end{align}

By Lemma~\ref{lem:covering_relation_magn_first} and
Proposition~\ref{prop:covering_relation_magn_equivalence}, each $\Equiv$-equivalence class
$[\PWord]_{\Equiv} \in \Magn_{1, 1}$ is equal to $\Han{\First_1(\PWord) \Conc
\First_1^\Reverse(\PWord)}_{\Equiv}$. Moreover, by Theorem~\ref{thm:p_symbol_magn}, the set
of $\Monoid$-pigmented magnets is a system of representatives of $\Magn_{1, 1}$. Therefore,
any $\Monoid$-pigmented magnet $\PWord$ can be represented by the concatenation
$\First_1(\PWord) \Conc \First_1^\Reverse(\PWord)$ of $\Monoid$-pigmented $1$-arrangements.
Since $\First_1(\PWord)$ and $\First_1^\Reverse(\PWord)$ have the same length, $\PWord$ can
be represented by the pair $\Par{\First_1(\PWord), \First_1^\Reverse(\PWord)}$. This
property leads to the fact that when $\Monoid$ is finite, for any $n \geq 0$,
\begin{equation} \label{equ:enumeration_magn}
    \# \Magn_{1, 1}(\Monoid)(n)
    = \sum_{j \in \HanL{n}} \binom{n}{j} j!^2 (\# \Monoid)^{2 j}.
\end{equation}
Let us explain~\eqref{equ:enumeration_magn}. A pair $\Par{\PWord', \PWord''}$ of
$\Monoid$-pigmented $1$-arrangements of arity $n \geq 0$ and length $j \geq 0$, provided
that there exists $\PWord \in \P(\Monoid)$ such that $\PWord' = \First_1(\PWord)$ and
$\PWord'' = \First_1^\Reverse(\PWord)$, is specified by
\begin{enumerate}[label=(S\arabic*)]
    \item \label{item:enumeration_magn_1}
    a subset $J$ of $[n]$, which is the common set of values for the $\Monoid$-pigmented
    letters of $\PWord'$ and~$\PWord''$;
    \item \label{item:enumeration_magn_2}
    two bijective maps $f'$ and $f''$ from $J$ to $[j]$, such that $f'(i)$ (resp.\ $f''(i)$)
    is the position of the $\Monoid$-pigmented letter of value $i$ in $\PWord'$ (resp.\
    $\PWord''$);
    \item \label{item:enumeration_magn_3}
    two words $w'$ and $w''$ of length $j$ on $\Monoid$ such that for any $r \in [j]$,
    $w'(r)$ (resp.\ $w''(r)$) is the pigment of the letter $\PWord'(r)$ (resp.\
    $\PWord''(r)$).
\end{enumerate}
The cardinality of $\Magn_{1, 1}(\Monoid)(n)$ expressed by~\eqref{equ:enumeration_magn}
follows from this specification. Indeed, \ref{item:enumeration_magn_1} gives rise to the
binomial coefficient enumerating all possible subsets $J$ of $[n]$ of cardinality $j$,
\ref{item:enumeration_magn_2} gives rise to the factor $j!^2$ enumerating all possible pairs
$\Par{f', f''}$ of bijections, and~\ref{item:enumeration_magn_3} gives rise to the last term
enumerating all possible pairs $\Par{w', w''}$ of such words.

In particular, the sequence of sizes of $\Magn_{1, 1}(\TrivialMonoid)$ starts by
\begin{equation}
    1, 2, 7, 52, 749, 17686, 614227, 29354312, 1844279257,
\end{equation}
and forms Sequence~\OEIS{A046662} of~\cite{Slo}.

\subsubsection{Presentation}
In order to establish a presentation of $\Magn_{1, 1}(\Monoid)$, we introduce an alternative
description of the clone congruence $\Equiv$ through a new equivalence relation $\Equiv'$.
For this, let us define $\Equiv'$ as the equivalence relation on $\P(\Monoid)$ generated by
\begin{equation} \label{equ:alternative_relation_magn_1}
    \PWord \Conc \QWord \Conc \QWord \Conc \PWord'
    \Equiv'
    \PWord \Conc \QWord \Conc \PWord',
\end{equation}
\begin{equation} \label{equ:alternative_relation_magn_2}
    \PWord \Conc \Par{\alpha_1 \MonoidProductExtension \QWord} \Conc
    \RWord \Conc \Par{\alpha_2 \MonoidProductExtension \QWord} \Conc
    \RWord' \Conc \Par{\alpha_3 \MonoidProductExtension \QWord} \Conc \PWord'
    \Equiv'
    \PWord \Conc \Par{\alpha_1 \MonoidProductExtension \QWord} \Conc
    \RWord \Conc
    \RWord' \Conc \Par{\alpha_3 \MonoidProductExtension \QWord} \Conc \PWord',
\end{equation}
where $\PWord, \PWord', \QWord, \RWord, \RWord' \in \P(\Monoid)$ and $\alpha_1, \alpha_2,
\alpha_3 \in \Monoid$.

\begin{Statement}{Lemma}{lem:alternative_relation_magn}
    For any monoid $\Monoid$, the binary relations $\Equiv$ and $\Equiv'$ on $\P(\Monoid)$
    are equal.
\end{Statement}
\begin{Proof}
    Let $\PWord, \PWord' \in \P(\Monoid)$ such that $\PWord \Equiv' \PWord'$. Since
    $\Equiv'$ is generated by~\eqref{equ:alternative_relation_magn_1}
    and~\eqref{equ:alternative_relation_magn_2}, we have two cases to consider.
    \begin{enumerate}[label={\it (\Roman*)}]
        \item If $\PWord$ and $\PWord'$ decompose as
        \begin{math}
            \PWord = \PWord'' \Conc \QWord \Conc \QWord \Conc \PWord'''
        \end{math}
        and
        \begin{math}
            \PWord' = \PWord'' \Conc \QWord \Conc \PWord'''
        \end{math}
        where $\PWord'', \PWord''', \QWord \in \P(\Monoid)$, then
        \begin{math}
            \First_1(\PWord)
            = \First_1\Par{\PWord'' \Conc \QWord \Conc \PWord'''}
            = \First_1\Par{\PWord'}
        \end{math}
        and
        \begin{math}
            \First_1^\Reverse(\PWord)
            = \First_1^\Reverse\Par{\PWord'' \Conc \QWord \Conc \PWord'''}
            = \First_1^\Reverse\Par{\PWord'}.
        \end{math}
        Therefore, $\PWord \Equiv \PWord'$.
        \item If $\PWord$ and $\PWord'$ decompose as
        \begin{math}
            \PWord
            = \PWord'' \Conc \Par{\alpha_1 \MonoidProductExtension \QWord} \Conc
            \RWord \Conc \Par{\alpha_2 \MonoidProductExtension \QWord} \Conc
            \RWord' \Conc \Par{\alpha_3 \MonoidProductExtension \QWord} \Conc \PWord'''
        \end{math}
        and
        \begin{math}
            \PWord'
            = \PWord'' \Conc \Par{\alpha_1 \MonoidProductExtension \QWord} \Conc
            \RWord \Conc \RWord' \Conc \Par{\alpha_3 \MonoidProductExtension \QWord}
            \Conc \PWord'''
        \end{math}
        where $\PWord'', \PWord''', \QWord, \RWord, \RWord' \in \P(\Monoid)$ and $\alpha_1,
        \alpha_2, \alpha_3 \in \Monoid$, then
        \begin{math}
            \First_1(\PWord)
            = \First_1\Par{
                \PWord'' \Conc \Par{\alpha_1 \MonoidProductExtension \QWord} \Conc
                \RWord \Conc \RWord' \Conc \PWord'''}
            = \First_1\Par{\PWord'}
        \end{math}
        and
        \begin{math}
            \First_1^\Reverse(\PWord)
            = \First_1^\Reverse\Par{
                \PWord'' \Conc \RWord \Conc \RWord' \Conc
                \Par{\alpha_3 \MonoidProductExtension \QWord} \Conc \PWord'''}
            = \First_1^\Reverse\Par{\PWord'}.
        \end{math}
        Therefore, $\PWord \Equiv \PWord'$.
    \end{enumerate}
    This shows that $\PWord \Equiv' \PWord'$ implies $\PWord \Equiv \PWord'$.

    Conversely, let $\PWord, \PWord' \in \P(\Monoid)$ such that $\PWord \Equiv \PWord'$. By
    Proposition~\ref{prop:covering_relation_magn_equivalence}, this is equivalent to the
    fact that $\PWord \EquivCovering \PWord'$. Since $\EquivCovering$ is generated by
    $\Covering$, we have three cases to explore depending whether $\PWord \Covering^{(1)}
    \PWord'$, $\PWord \Covering^{(2)} \PWord'$, or $\PWord \Covering^{(3)} \PWord'$.
    \begin{enumerate}[label={\it (\Roman*)}]
        \item If $\PWord \Covering^{(1)} \PWord'$, then $\PWord$ and $\PWord'$ decompose as
        \begin{math}
            \PWord = \QWord \Conc \Witness{n}{n}{i^\alpha} \Conc \QWord'
        \end{math}
        and
        \begin{math}
            \PWord' = \QWord \Conc \QWord'
        \end{math}
        where $\QWord, \QWord' \in \P(\Monoid)$ and $i^\alpha \in \SetLetters_\Monoid$.
        Since the position $\Length(\QWord) + 1$ of $\PWord$ is neither a left $1$-witness
        nor a right $1$-witness, there is necessarily an occurrence of an
        $\Monoid$-pigmented letter having $i$ as value both in $\QWord$ and in $\QWord'$.
        Hence,
        \begin{math}
            \PWord = \RWord \Conc \Witness{y}{n}{i^{\alpha_1}} \Conc \RWord' \Conc
            \Witness{n}{n}{i^\alpha} \Conc \RWord'' \Conc \Witness{n}{y}{i^{\alpha_2}}
            \Conc \RWord'''
        \end{math}
        and
        \begin{math}
            \PWord' = \RWord \Conc \Witness{y}{n}{i^{\alpha_1}} \Conc \RWord' \Conc \RWord''
            \Conc \Witness{n}{y}{i^{\alpha_2}} \Conc \RWord'''
        \end{math}
        where $\RWord, \RWord', \RWord'', \RWord''' \in \P(\Monoid)$, $\QWord = \RWord \Conc
        i^{\alpha_1} \Conc \RWord'$, $\QWord' = \RWord'' \Conc i^{\alpha_2} \Conc
        \RWord'''$, and $i^{\alpha_1}, i^{\alpha_2} \in \SetLetters_\Monoid$.
        By~\eqref{equ:alternative_relation_magn_2}, we have $\PWord \Equiv' \PWord'$.
        \item If $\PWord \Covering^{(2)} \PWord'$, then $\PWord$ and $\PWord'$ decompose as
        \begin{math}
            \PWord = \QWord \Conc \Witness{u}{n}{i_1^{\alpha_1}}
            \Witness{n}{y}{i_2^{\alpha_2}} \Conc \QWord'
        \end{math}
        and
        \begin{math}
            \PWord' = \QWord \Conc \Witness{n}{y}{i_2^{\alpha_2}}
            \Witness{u}{n}{i_1^{\alpha_1}} \Conc \QWord'
        \end{math}
        where $\QWord, \QWord' \in \P(\Monoid)$, $i_1^{\alpha_1}, i_2^{\alpha_2} \in
        \SetLetters_\Monoid$, and $i_1 \ne i_2$. Since the position $\Length(\QWord) + 1$ of
        $\PWord$ is not a right $1$-witness and the position $\Length(\QWord) + 2$ of
        $\PWord$ is not a left $1$-witness, there is necessarily an occurrence of an
        $\Monoid$-pigmented letter having $i_2$ as value in $\QWord$ and an occurrence of an
        $\Monoid$-pigmented letter having $i_1$ as value in $\QWord'$. Hence,
        \begin{align} \label{equ:alternative_relation_magn}
            \PWord & = \RWord \Conc \Witness{u}{n}{i_2^{\beta_2}} \Conc \RWord'
                \Conc \Witness{u}{n}{i_1^{\alpha_1}} \Witness{n}{y}{i_2^{\alpha_2}} \Conc
                \RWord'' \Conc \Witness{n}{u}{i_1^{\beta_1}} \Conc \RWord'''
            \\
            & \Equiv'
            \RWord \Conc \Witness{u}{n}{i_2^{\beta_2}} \Conc \RWord' \Conc
                \Witness{n}{n}{i_2^{\alpha_2}} \Witness{u}{n}{i_1^{\alpha_1}}
                \Witness{n}{y}{i_2^{\alpha_2}} \Conc \RWord'' \Conc
                \Witness{n}{u}{i_1^{\beta_1}} \Conc \RWord'''
            \notag \\
            & \Equiv'
            \RWord \Conc \Witness{u}{n}{i_2^{\beta_2}} \Conc \RWord' \Conc
                \Witness{n}{n}{i_2^{\alpha_2}} \Witness{u}{n}{i_1^{\alpha_1}}
                \Witness{n}{y}{i_2^{\alpha_2}} \Witness{n}{n}{i_1^{\alpha_1}} \Conc
                \RWord'' \Conc \Witness{n}{u}{i_1^{\beta_1}} \Conc \RWord'''
            \notag \\
            & \Equiv'
            \RWord \Conc \Witness{u}{n}{i_2^{\beta_2}} \Conc \RWord' \Conc
                \Witness{n}{y}{i_2^{\alpha_2}} \Witness{u}{n}{i_1^{\alpha_1}} \Conc
                \RWord'' \Conc \Witness{n}{u}{i_1^{\beta_1}} \Conc \RWord'''
            = \PWord'
            \notag
        \end{align}
        where $\RWord, \RWord', \RWord'', \RWord''' \in \P(\Monoid)$, $\QWord = \RWord \Conc
        i_2^{\beta_2} \Conc \RWord'$, $\QWord' = \RWord'' \Conc i_1^{\beta_1} \Conc
        \RWord'''$, and $i_1^{\beta_1}, i_2^{\beta_2} \in \SetLetters_\Monoid$. The first
        and second $\Equiv'$-equivalences of~\eqref{equ:alternative_relation_magn} are
        consequences of~\eqref{equ:alternative_relation_magn_2} considered from right to
        left and the third $\Equiv'$-equivalence of~\eqref{equ:alternative_relation_magn} is
        a consequence of~\eqref{equ:alternative_relation_magn_1} considered from left to
        right.
        \item If $\PWord \Covering^{(3)} \PWord'$, then $\PWord$ and $\PWord'$ decompose as
        \begin{math}
            \PWord
            = \QWord \Conc \Witness{u}{n}{i^\alpha} \Witness{n}{u}{i^\alpha} \Conc \QWord'
        \end{math}
        and
        \begin{math}
            \PWord' = \QWord \Conc \Witness{u}{u}{i^\alpha} \Conc \QWord'
        \end{math}
        where $\QWord, \QWord' \in \P(\Monoid)$ and $i^\alpha \in \SetLetters_\Monoid$.
        By~\eqref{equ:alternative_relation_magn_1}, we have $\PWord \Equiv' \PWord'$.
    \end{enumerate}
    This shows that $\PWord \Equiv \PWord'$ implies $\PWord \Equiv' \PWord'$ and establishes
    the statement of the lemma.
\end{Proof}

\begin{Statement}{Theorem}{thm:presentation_magn}
    For any monoid $\Monoid$, the clone $\Magn_{1, 1}(\Monoid)$ admits the presentation
    $\Par{\GeneratingSet_\Monoid, \RelationSet_\Monoid'}$ where $\RelationSet_\Monoid'$ is
    the set $\RelationSet_\Monoid$ augmented with the $\GeneratingSet_\Monoid$-equations
    \begin{equation} \label{equ:presentation_magn_1}
        \RightComb_\Monoid\Par{1^\MonoidUnit 1^\MonoidUnit}
        \ \RelationSet_\Monoid' \
        \RightComb_\Monoid\Par{1^\MonoidUnit},
    \end{equation}
    \begin{equation} \label{equ:presentation_magn_2}
        \RightComb_\Monoid\Par{1^{\alpha_1} 2^\MonoidUnit 1^{\alpha_2} 3^\MonoidUnit
            1^{\alpha_3}}
        \ \RelationSet_\Monoid' \
        \RightComb_\Monoid\Par{1^{\alpha_1} 2^\MonoidUnit 3^\MonoidUnit 1^{\alpha_3}},
    \end{equation}
    where $\alpha_1, \alpha_2, \alpha_3 \in \Monoid$ and $\MonoidUnit$ is the unit of
    $\Monoid$.
\end{Statement}
\begin{Proof}
    Let $\Equiv''$ be the clone congruence of $\P(\Monoid)$ generated by
    \begin{equation}
        1^\MonoidUnit 1^\MonoidUnit \Equiv'' 1^\MonoidUnit,
    \end{equation}
    \begin{equation}
        1^{\alpha_1} 2^\MonoidUnit 1^{\alpha_2} 3^\MonoidUnit 1^{\alpha_3}
        \Equiv''
        1^{\alpha_1} 2^\MonoidUnit 3^\MonoidUnit 1^{\alpha_3}
    \end{equation}
    with $\alpha_1, \alpha_2, \alpha_3 \in \Monoid$. Let us show that the clone congruences
    $\Equiv$ and $\Equiv''$ of $\P(\Monoid)$ are equal. This will imply, by the remark
    stated in Section~\ref{subsubsec:remark_presentations}, that $\P(\Monoid) /_{\Equiv''} =
    \P(\Monoid) /_{\Equiv} = \Magn_{1, 1}(\Monoid)$ admits the stated presentation.

    First, since $\First_1\Par{1^\MonoidUnit 1^\MonoidUnit} = 1^\MonoidUnit =
    \First_1\Par{1^\MonoidUnit}$ and $\First_1^\Reverse\Par{1^\MonoidUnit 1^\MonoidUnit} =
    1^\MonoidUnit = \First_1^\Reverse\Par{1^\MonoidUnit}$, we have $1^\MonoidUnit
    1^\MonoidUnit \Equiv 1^\MonoidUnit$. Moreover, since for any $\alpha_1, \alpha_2,
    \alpha_3 \in \Monoid$,
    \begin{math}
        \First_1\Par{
        1^{\alpha_1} 2^\MonoidUnit 1^{\alpha_2} 3^\MonoidUnit 1^{\alpha_3}}
        = 1^{\alpha_1} 2^\MonoidUnit 3^\MonoidUnit
        = \First_1\Par{1^{\alpha_1} 2^\MonoidUnit 3^\MonoidUnit 1^{\alpha_3}}
    \end{math}
    and
    \begin{math}
        \First_1^\Reverse\Par{
        1^{\alpha_1} 2^\MonoidUnit 1^{\alpha_2} 3^\MonoidUnit 1^{\alpha_3}}
        = 2^\MonoidUnit 3^\MonoidUnit 1^{\alpha_3}
        = \First_1^\Reverse\Par{1^{\alpha_1} 2^\MonoidUnit 3^\MonoidUnit 1^{\alpha_3}},
    \end{math}
    we have
    \begin{math}
        1^{\alpha_1} 2^\MonoidUnit 1^{\alpha_2} 3^\MonoidUnit 1^{\alpha_3}
        \Equiv
        1^{\alpha_1} 2^\MonoidUnit 3^\MonoidUnit 1^{\alpha_3}.
    \end{math}
    This shows that $\Equiv''$ is contained in $\Equiv$.

    To prove that $\Equiv$ is contained in $\Equiv''$, let us show that $\Equiv'$ is
    contained in $\Equiv''$. By Lemma~\ref{lem:alternative_relation_magn}, the targeted
    property will follow. For any $\PWord, \PWord', \QWord \in \P(\Monoid)$, we have
    \begin{equation} \label{equ:presentation_magn_A}
        \PWord \Conc \QWord \Conc \QWord \Conc \PWord'
        = 1^\MonoidUnit 2^\MonoidUnit 3^\MonoidUnit \Superposition{
            \PWord, 1^\MonoidUnit 1^\MonoidUnit \Superposition{\QWord}, \PWord'}
        \Equiv'' 1^\MonoidUnit 2^\MonoidUnit 3^\MonoidUnit \Superposition{
            \PWord, 1^\MonoidUnit \Superposition{\QWord}, \PWord'}
        = \PWord \Conc \QWord \Conc \PWord'
    \end{equation}
    so that the first and the last members of~\eqref{equ:presentation_magn_A} are
    $\Equiv''$-equivalent. Moreover, for any $\PWord, \PWord', \QWord, \RWord, \RWord' \in
    \P(\Monoid)$ and $\alpha_1, \alpha_2, \alpha_3 \in \Monoid$, we have
    \begin{align} \label{equ:presentation_magn_B}
        \PWord \Conc \Par{\alpha_1 \MonoidProductExtension \QWord} \Conc \RWord \Conc
            \Par{\alpha_2 \MonoidProductExtension \QWord} \Conc \RWord'
            \Conc \Par{\alpha_3 \MonoidProductExtension \QWord} \Conc \PWord'
        & =
        1^\MonoidUnit 2^\MonoidUnit 3^\MonoidUnit \Superposition{
            \PWord,
            1^{\alpha_1} 2^\MonoidUnit 1^{\alpha_2} 3^\MonoidUnit 1^{\alpha_3}
            \Superposition{
                \QWord, \RWord, \RWord'},
            \PWord'}
        \\
        & \Equiv''
        1^\MonoidUnit 2^\MonoidUnit 3^\MonoidUnit \Superposition{
            \PWord,
            1^{\alpha_1} 2^\MonoidUnit 3^\MonoidUnit 1^{\alpha_3} \Superposition{
                \QWord, \RWord, \RWord'},
            \PWord'}
        \notag
        \\
        & =
        \PWord \Conc \Par{\alpha_1 \MonoidProductExtension \QWord} \Conc \RWord \Conc
            \RWord'
            \Conc \Par{\alpha_3 \MonoidProductExtension \QWord} \Conc \PWord'
        \notag
    \end{align}
    so that the first and the last members of~\eqref{equ:presentation_magn_B} are
    $\Equiv''$-equivalent. Since $\Equiv'$ is the equivalence relation generated
    by~\eqref{equ:alternative_relation_magn_1} and~\eqref{equ:alternative_relation_magn_2},
    the targeted property is shown. This establishes the statement of the theorem.
\end{Proof}

By Theorem~\ref{thm:presentation_magn}, any $\Magn_{1, 1}(\Monoid)$-algebra is, up to term
equivalence, an $\Monoid$-pigmented monoid $\Par{\Algebra, \Product, \Identity,
\Par{\Pigmentation_\alpha}_{\alpha \in \Monoid}}$ where $\Product$ is idempotent, and
$\Product$ and $\Par{\Pigmentation_\alpha}_{\alpha \in \Monoid}$ satisfy, by spelling
out~\eqref{equ:presentation_magn_2} and simplifying it modulo the background theory
$\Equiv_{\RelationSet_\Monoid}$,
\begin{equation}
    \Pigmentation_{\alpha_1} \Par{x_1} \Product x_2 \Product
    \Pigmentation_{\alpha_2} \Par{x_1} \Product x_3 \Product
    \Pigmentation_{\alpha_3} \Par {x_1}
    =
    \Pigmentation_{\alpha_1} \Par{x_1} \Product x_2 \Product x_3 \Product
    \Pigmentation_{\alpha_3} \Par {x_1}
\end{equation} 
for any $x_1, x_2, x_3 \in \Algebra$ and $\alpha_1, \alpha_2, \alpha_3 \in \Monoid$. In
particular, $\Magn_{1, 1}(\TrivialMonoid)$ is a clone realization of the variety of regular
bands equipped with an additional unary operation acting identically.

\subsection{On pigmented stalactites} \label{subsec:stal}
By considering the intersection of the clone congruences $\Equiv_\Sort$ and
$\Equiv_{\First_k}$, $k \geq 0$, we construct a quotient clone $\Stal_k(\Monoid)$ of
$\P(\Monoid)$. A description through new combinatorial objects named $\Monoid$-pigmented
stalactites is introduced and a finitely equationally axiomatizable presentation is
described. These results are based on the introduction of a $\PSymbol$-symbol for the
underlying equivalence relation.

\subsubsection{Clone construction}
For any parameter $k \geq 0$, let $\Equiv_k$ be the clone congruence $\Equiv_\Sort \cap
\Equiv_{\First_k}$ and
\begin{equation}
    \Stal_k(\Monoid) := \P(\Monoid) /_{\Equiv_k}.
\end{equation}
By Propositions~\ref{prop:congruence_sort} and~\ref{prop:congruence_first},
$\Stal_k(\Monoid)$ is a well-defined clone, and $\WInc(\Monoid)$ and $\Arra_k(\Monoid)$ are
both isomorphic to quotients of $\Stal_k(\Monoid)$. Since for any $0 \leq k \leq k'$,
$\Equiv_{k'}$ is a refinement of $\Equiv_k$, $\Stal_k(\Monoid)$ is isomorphic to a quotient
of $\Stal_{k'}(\Monoid)$. Moreover, since $\Equiv_0$ and $\Equiv_\Sort$ are the same
equivalence relations, $\Stal_0(\Monoid)$ is identical to $\WInc(\Monoid)$. Besides, the
clone
\begin{math}
    \Stal_k^\Reverse(\Monoid)
    := \Stal_k(\Monoid)^\Reverse
    = \P(\Monoid) /_{\Equiv_{\Sort} \cap \Equiv_{\First_k}^\Reverse}
\end{math}
is by
Proposition~\ref{prop:reversed_congruence} isomorphic to $\Stal_k(\Monoid)$.

\subsubsection{Equivalence relation}
By definition, for any $\PWord, \PWord' \in \P(\Monoid)$, $\PWord \Equiv_k \PWord'$ holds if
and only
\begin{math}
    \Par{\First_k\Par{\PWord}, \Sort_{\Leq}\Par{\PWord}}
    = \Par{\First_k\Par{\PWord'}, \Sort_{\Leq}\Par{\PWord'}}
\end{math}
where $\Leq$ is any total order relation on $\Monoid$.

In order to obtain properties about the clone $\Stal_k(\Monoid)$, $k \geq 0$, we introduce
an alternative equivalence relation $\EquivCovering_k$ for which it appears that it is
equal to $\Equiv_k$. Let $\Covering^{(1)}_k$ and $\Covering^{(2)}_k$ be the two binary
relations on $\P(\Monoid)$ satisfying
\begin{equation} \label{equ:covering_relation_stal_1}
    \PWord \Conc \Witness{n}{u}{i_1^{\alpha_1}} \Witness{y}{u}{i_2^{\alpha_2}} \Conc \PWord'
    \enspace \Covering^{(1)}_k \enspace
    \PWord \Conc \Witness{y}{u}{i_2^{\alpha_2}} \Witness{n}{u}{i_1^{\alpha_1}} \Conc
        \PWord'
    \qquad
    \text{where }
    i_1 \ne i_2,
\end{equation}
\begin{equation} \label{equ:covering_relation_stal_2}
    \PWord \Conc \Witness{n}{u}{i_1^{\alpha_1}} \Witness{n}{u}{i_2^{\alpha_2}} \Conc \PWord'
    \enspace \Covering^{(2)}_k \enspace
    \PWord \Conc \Witness{n}{u}{i_2^{\alpha_2}} \Witness{n}{u}{i_1^{\alpha_1}} \Conc
        \PWord'
    \qquad
    \text{where }
    i_1^{\alpha_1} \ne i_2^{\alpha_2}, \text{ and } i_2^{\alpha_2} \PLeq i_1^{\alpha_1},
\end{equation}
where $\PWord, \PWord' \in \P(\Monoid)$ and $i_1^{\alpha_1}, i_2^{\alpha_2} \in
\SetLetters_\Monoid$. Note that these definitions depend on $k$ because the properties of
being a left $k$-witness of the shown pigmented letters
in~\eqref{equ:covering_relation_stal_1} and~\eqref{equ:covering_relation_stal_2} depend
themselves on $k$. Let $\CoveringRT^{(j)}_k$, $j \in [2]$, be the reflexive and transitive
closure of $\Covering^{(j)}_k$, $\Covering_k$ be the union $\Covering^{(1)}_k \cup
\Covering^{(2)}_k$, and $\EquivCovering_k$ be the reflexive, symmetric, and transitive
closure of~$\Covering_k$.

Let $\RemoveFirst_k : \P(\Monoid) \to \P(\Monoid)$ be the map sending any $\PWord \in
\P(\Monoid)$ to the $\Monoid$-pigmented word defined as the subword of $\PWord$ consisting
of the letters whose positions are not left $k$-witnesses. For instance, in $\P\Par{A^*}$,
where $A^*$ is the free monoid over $\Bra{a, b}$, we have
\begin{equation}
    \RemoveFirst_1\Par{
        \Witness{y}{u}{1^\epsilon} \Witness{y}{u}{3^{ab}} \Witness{n}{u}{1^b}
        \Witness{n}{u}{3^b} \Witness{n}{u}{1^{aa}} \Witness{n}{u}{3^\epsilon}
        \Witness{y}{u}{2^{aa}} \Witness{n}{u}{3^{bba}}}
    = 1^b 3^b 1^{aa} 3^\epsilon 3^{bba},
\end{equation}
\begin{equation}
    \RemoveFirst_2\Par{
        \Witness{y}{u}{1^\epsilon} \Witness{y}{u}{3^{ab}} \Witness{y}{u}{1^b}
        \Witness{y}{u}{3^b} \Witness{n}{u}{1^{aa}} \Witness{n}{u}{3^\epsilon}
        \Witness{y}{u}{2^{aa}} \Witness{n}{u}{3^{bba}}}
    = 1^{aa} 3^\epsilon  3^{bba}.
\end{equation}

\begin{Statement}{Lemma}{lem:covering_relation_stal}
    For any monoid $\Monoid$ endowed with a total order relation $\Leq$, any $k \geq 0$, and
    any $\PWord \in \P(\Monoid)$,
    \begin{equation}
        \PWord
        \CoveringRT^{(1)}_k
        \First_k(\PWord) \Conc \RemoveFirst_k(\PWord)
        \CoveringRT^{(2)}_k
        \First_k(\PWord) \Conc \Sort_{\Leq} \Par{\RemoveFirst_k(\PWord)}.
    \end{equation}
\end{Statement}
\begin{Proof}
    The binary relation $\Covering^{(1)}_k$ acts on an $\Monoid$-pigmented word by
    transposing two of its positions $j$ and $j + 1$ such that the position $j$ is not a
    left $k$-witness while the position $j + 1$ is. Therefore, iterated applications of
    $\Covering^{(1)}_k$ on $\PWord$ move the $\Monoid$-pigmented letters whose positions are
    left $k$-witnesses to the left of the word while preserving the relative order among such
    letters and while preserving the relative order of $\Monoid$-pigmented letters whose
    positions are not left $k$-witnesses. Thus, this process builds the $\Monoid$-pigmented
    word $\First_k(\PWord) \Conc \QWord$, where $\QWord$ is the subword of $\PWord$ formed
    by the $\Monoid$-pigmented letters of $\PWord$ such that their positions in $\PWord$ are
    not left $k$-witnesses. Hence, $\QWord$ is the $\Monoid$-pigmented word
    $\RemoveFirst_k(\PWord)$. This shows that
    \begin{math}
        \PWord \CoveringRT^{(1)}_k \First_k(\PWord) \Conc \RemoveFirst_k(\PWord).
    \end{math}

    Now, observe that $\Covering^{(2)}_k$ acts on an $\Monoid$-pigmented word by transposing
    two of its positions $j$ and $j + 1$ which are not left $k$-witnesses and such that the
    $\Monoid$-pigmented letter at position $j$ is greater than the $\Monoid$-pigmented
    letter at position $j + 1$ w.r.t.\ the order relation $\PLeq$. Therefore, iterated
    applications of $\Covering^{(2)}_k$ on $\First_k(\PWord) \Conc \RemoveFirst_k(\PWord)$
    sort w.r.t.\ $\PLeq$ the suffix of this word made of letters whose positions are not
    left $k$-witnesses. Thus, this process builds the $\Monoid$-pigmented word
    $\First_k(\PWord) \Conc \Sort_{\Leq}\Par{\RemoveFirst_k(\PWord)}$. This shows that
    \begin{math}
        \First_k(\PWord) \Conc \RemoveFirst_k(\PWord)
        \CoveringRT^{(2)}_k
        \First_k(\PWord) \Conc \Sort_{\Leq} \Par{\RemoveFirst_k(\PWord)}
    \end{math}
    as expected.
\end{Proof}

\begin{Statement}{Proposition}{prop:covering_relation_stal_equivalence}
    For any monoid $\Monoid$ and any $k \geq 0$, the binary relations $\Equiv_k$ and
    $\EquivCovering_k$ on $\P(\Monoid)$ are equal.
\end{Statement}
\begin{Proof}
    Let $\Leq$ be any total order relation on $\Monoid$.

    Let us show that $\PWord \Covering_k \PWord'$ for $\PWord, \PWord' \in \P(\Monoid)$
    entails $\PWord \Equiv_k \PWord'$. Since the equivalence relation $\EquivCovering_k$ is
    generated by $\Covering_k$, this will show that $\EquivCovering_k$ is contained in
    $\Equiv_k$. We are going to consider the following two cases depending whether $\PWord
    \Covering^{(1)}_k \PWord'$ or $\PWord \Covering^{(2)}_k \PWord'$.
    \begin{enumerate}[label={\it (\Roman*)}]
        \item Assume that $\PWord \Covering^{(1)}_k \PWord'$. Hence, $\PWord$ and $\PWord'$
        decompose as $\PWord = \QWord \Conc \Witness{n}{u}{i_1^{\alpha_1}}
        \Witness{y}{u}{i_2^{\alpha_2}} \Conc \RWord$ and $\PWord' = \QWord \Conc
        \Witness{y}{u}{i_2^{\alpha_2}} \Witness{n}{u}{i_1^{\alpha_1}} \Conc \RWord$ where
        $\QWord, \RWord \in \P(\Monoid)$, $i_1^{\alpha_1}, i_2^{\alpha_2} \in
        \SetLetters_\Monoid$, and $i_1 \ne i_2$. Since
        \begin{equation}
            \First_k \Par{\PWord}
            = \First_k \Par{
                \QWord \Conc \Witness{y}{u}{i_2^{\alpha_2}} \Conc \RWord
            }
            = \First_k \Par{\PWord'}
        \end{equation}
        and
        \begin{equation}
            \Sort_{\Leq}\Par{\PWord}
            = \Sort_{\Leq}\Par{
                \QWord \Conc
                \Witness{y}{u}{i_2^{\alpha_2}}
                \Witness{n}{u}{i_1^{\alpha_1}}
                \Conc \RWord
            }
            = \Sort_{\Leq} \Par{\PWord'},
        \end{equation}
        by definition of $\Equiv_k$, $\PWord \Equiv_k \PWord'$.
        \item Assume that $\PWord \Covering^{(2)}_k \PWord'$. Hence, $\PWord$ and $\PWord'$
        decompose as $\PWord = \QWord \Conc \Witness{n}{u}{i_1^{\alpha_1}}
        \Witness{n}{u}{i_2^{\alpha_2}} \Conc \RWord$ and $\PWord' = \QWord \Conc
        \Witness{n}{u}{i_2^{\alpha_2}} \Witness{n}{u}{i_1^{\alpha_1}} \Conc \RWord$ where
        $\QWord, \RWord \in \P(\Monoid)$, $i_1^{\alpha_1}, i_2^{\alpha_2} \in
        \SetLetters_\Monoid$, $i_1^{\alpha_1} \ne i_2^{\alpha_2}$, and $i_2^{\alpha_2} \PLeq
        i_1^{\alpha_1}$. Since
        \begin{equation}
            \First_k \Par{\PWord}
            = \First_k \Par{
                \QWord \Conc \RWord
            }
            =
            \First_k \Par{\PWord'}
        \end{equation}
        and
        \begin{equation}
            \Sort_{\Leq}\Par{\PWord}
            = \Sort_{\Leq} \Par{
                \QWord \Conc
                \Witness{n}{u}{i_2^{\alpha_2}}
                \Witness{n}{u}{i_1^{\alpha_1}} \Conc
                \RWord
            }
            =
            \Sort_{\Leq}\Par{\PWord'},
        \end{equation}
        by definition of $\Equiv_k$, $\PWord \Equiv_k \PWord'$.
    \end{enumerate}

    Conversely, for any $\PWord, \PWord' \in \P(\Monoid)$ such that $\PWord \Equiv_k
    \PWord'$, we have $\First_k\Par{\PWord} = \First_k\Par{\PWord'}$ and
    $\Sort_{\Leq}\Par{\PWord} = \Sort_{\Leq}\Par{\PWord'}$. This implies that $\PWord$ and
    $\PWord'$ are built on the same multiset of $\Monoid$-pigmented letters and thus, since
    $\First_k(\PWord) = \First_k\Par{\PWord'}$, the $\Monoid$-pigmented words
    $\RemoveFirst_k(\PWord)$ and $\RemoveFirst_k\Par{\PWord'}$ are also built on the same
    multiset of $\Monoid$-pigmented letters. Therefore, by
    Lemma~\ref{lem:covering_relation_stal},
    \begin{math}
        \PWord \EquivCovering_k
        \First_k(\PWord) \Conc \Sort_{\Leq}\Par{\RemoveFirst_k(\PWord)}
        = \First_k\Par{\PWord'} \Conc \Sort_{\Leq}\Par{\RemoveFirst_k\Par{\PWord'}}
        \EquivCovering_k \PWord'.
    \end{math}
    For this reason, we have $\PWord \EquivCovering_k \PWord'$, showing that $\Equiv_k$ is
    contained in~$\EquivCovering_k$.
\end{Proof}

\subsubsection{$\PSymbol$-symbol algorithm}
With the aim of describing $\Stal_k(\Monoid)$, we propose now a $\PSymbol$-symbol for
$\Equiv_k$. Let $\PSymbol_{\Equiv_k} : \P(\Monoid) \to \P(\Monoid)$ be the map defined for
any $\PWord \in \P(\Monoid)$ by
\begin{equation}
    \PSymbol_{\Equiv_k}(\PWord)
    := \First_k(\PWord) \Conc \Sort_{\Leq} \Par{\RemoveFirst_k(\PWord)}.
\end{equation}
For instance, in $\P\Par{A^*}$, where $A^*$ is the free monoid over $\Bra{a, b}$, we have
\begin{equation}
    \PSymbol_{\Equiv_1} \Par{
        \Witness{y}{u}{3^a} \Witness{y}{u}{2^\epsilon} \Witness{y}{u}{1^a}
        \Witness{n}{u}{1^b} \Witness{n}{u}{1^{ba}} \Witness{n}{u}{2^\epsilon}
        \Witness{n}{u}{1^{ba}} \Witness{n}{u}{1^\epsilon} \Witness{n}{u}{2^a}
        \Witness{y}{u}{4^a} \Witness{n}{u}{4^b}}
    =
    \Witness{y}{u}{3^a} \Witness{y}{u}{2^\epsilon} \Witness{y}{u}{1^a}
    \Witness{y}{u}{4^a} \Witness{n}{u}{1^\epsilon} \Witness{n}{u}{1^b}
    \Witness{n}{u}{1^{ba}} \Witness{n}{u}{1^{ba}} \Witness{n}{u}{2^\epsilon}
    \Witness{n}{u}{2^a} \Witness{n}{u}{4^b}
\end{equation}
and
\begin{equation}
    \PSymbol_{\Equiv_2} \Par{
        \Witness{y}{u}{3^a} \Witness{y}{u}{2^\epsilon} \Witness{y}{u}{1^a}
        \Witness{y}{u}{1^b} \Witness{n}{u}{1^{ba}} \Witness{y}{u}{2^\epsilon}
        \Witness{n}{u}{1^{ba}} \Witness{n}{u}{1^\epsilon} \Witness{n}{u}{2^a}
        \Witness{y}{u}{4^a} \Witness{y}{u}{4^b}}
    =
    \Witness{y}{u}{3^a} \Witness{y}{u}{2^\epsilon} \Witness{y}{u}{1^a} \Witness{y}{u}{1^b}
    \Witness{y}{u}{2^\epsilon} \Witness{y}{u}{4^a} \Witness{y}{u}{4^b}
    \Witness{n}{u}{1^\epsilon} \Witness{n}{u}{1^{ba}} \Witness{n}{u}{1^{ba}}
    \Witness{n}{u}{2^a}.
\end{equation}

\begin{Statement}{Lemma}{lem:p_symbol_stal_equiv_p_symbol}
    For any monoid $\Monoid$, any $k \geq 0$, and any $\PWord \in \P(\Monoid)$, $\PWord
    \EquivCovering_k \PSymbol_{\Equiv_k}(\PWord)$.
\end{Statement}
\begin{Proof}
    By definition of $\EquivCovering_k$, $\CoveringRT^{(1)}_k$, and $\CoveringRT^{(2)}_k$,
    it follows that both $\CoveringRT^{(1)}_k$ and $\CoveringRT^{(2)}_k$ are contained in
    $\EquivCovering_k$. The statement of the lemma is now a direct consequence of
    Lemma~\ref{lem:covering_relation_stal} and of the definition of the
    map~$\PSymbol_{\Equiv_k}$.
\end{Proof}

\begin{Statement}{Lemma}{lem:covering_relation_stal_p_symbol}
    For any monoid $\Monoid$, any $k \geq 0$, and any $\PWord, \PWord' \in \P(\Monoid)$,
    $\PWord \EquivCovering_k \PWord'$ implies $\PSymbol_{\Equiv_k}(\PWord) =
    \PSymbol_{\Equiv_k}\Par{\PWord'}$.
\end{Statement}
\begin{Proof}
    Assume that $\PWord \EquivCovering_k \PWord'$. By
    Proposition~\ref{prop:covering_relation_stal_equivalence}, this implies $\PWord \Equiv_k
    \PWord'$. Now, by definition of the map $\PSymbol_{\Equiv_k}$, this entails
    $\PSymbol_{\Equiv_k}(\PWord) = \PSymbol_{\Equiv_k}\Par{\PWord'}$.
\end{Proof}

By Proposition~\ref{prop:covering_relation_stal_equivalence} and
Lemmas~\ref{lem:p_symbol_stal_equiv_p_symbol} and~\ref{lem:covering_relation_stal_p_symbol},
$\PSymbol_{\Equiv_k}$ is a $\PSymbol$-symbol for~$\Equiv_k$.

\subsubsection{Description}
An \Def{$\Monoid$-pigmented $k$-stalactite} (or simply \Def{pigmented $k$-stalactite} when
the context is clear) of arity $n \geq 0$ is an $\Monoid$-pigmented word $\PWord$ of arity
$n$ which is in the image of $\PSymbol_{\Equiv_k}$. For instance, in $\P\Par{A^*}$, where
$A^*$ is the free monoid over $\Bra{a, b}$,
\begin{equation}
    \Witness{y}{u}{3^b} \Witness{y}{u}{2^a} \Witness{y}{u}{2^a} \Witness{n}{u}{2^b}
    \Witness{y}{u}{3^b} \Witness{y}{u}{1^a} \Witness{n}{u}{3^b} \Witness{y}{u}{1^b}
    \Witness{n}{u}{3^a}
\end{equation}
is not an $A^*$-pigmented $2$-stalactite. In contrast,
\begin{equation}
    \Witness{y}{u}{2^b} \Witness{y}{u}{2^a} \Witness{y}{u}{1^a} \Witness{y}{u}{3^a}
    \Witness{y}{u}{4^b} \Witness{y}{u}{1^a} \Witness{y}{u}{3^b} \Witness{n}{u}{1^a}
    \Witness{n}{u}{2^b} \Witness{n}{u}{2^b} \Witness{n}{u}{3^a} \Witness{n}{u}{3^b}
\end{equation}
is an $A^*$-pigmented $2$-stalactite but not an $A^*$-pigmented $1$-stalactite.

\begin{Statement}{Theorem}{thm:p_symbol_stal}
    For any monoid $\Monoid$ and any $k \geq 0$, $\PSymbol_{\Equiv_k}$ is a
    $\PSymbol$-symbol for $\Equiv_k$ and $\PSymbol_{\Equiv_k}(\P(\Monoid))$ is the set of
    $\Monoid$-pigmented $k$-stalactites. Moreover, the graded set $\Stal_k(\Monoid)$ is
    isomorphic to the graded set of $\Monoid$-pigmented $k$-stalactites.
\end{Statement}
\begin{Proof}
    By Proposition~\ref{prop:covering_relation_stal_equivalence} and
    Lemmas~\ref{lem:p_symbol_stal_equiv_p_symbol}
    and~\ref{lem:covering_relation_stal_p_symbol}, $\PSymbol_{\Equiv_k}$ is a
    $\PSymbol$-symbol for $\Equiv_k$. Moreover, the set of $\Monoid$-pigmented
    $k$-stalactites is defined as the set $\PSymbol_{\Equiv_k}(\P(\Monoid))$. The last part
    of the statement is a direct implication of
    Proposition~\ref{prop:clone_realization_p_symbol} and the fact that
    $\PSymbol_{\Equiv_k}$ is, as we have just shown, a $\PSymbol$-symbol for~$\Equiv_k$.
\end{Proof}

By Proposition~\ref{prop:clone_realization_p_symbol} and Theorem~\ref{thm:p_symbol_stal},
$\Stal_k(\Monoid)$ can be seen as a clone on $\Monoid$-pigmented $k$-stalactites with
superposition maps satisfying~\eqref{equ:superposition_p_symbol}. For instance, in
$\Stal_1\Par{A^*}$, where $A^*$ is the free monoid over $\Bra{a, b}$, we have, up to
isomorphism,
\begin{align}
    4^{ab} 1^a 2^{ab} 3^a 3^\epsilon & \Superposition{
        2^{ba} 3^b,
        3^{ba} 1^b 1^b 3^\epsilon,
        2^\epsilon 3^{ab} 2^{ba} 3^b,
        2^a}
    \\
    & =
    \PSymbol_{\Equiv_1}\Par{
        \Witness{y}{u}{2^{aba}}
        \Witness{n}{u}{2^{aba}}
        \Witness{y}{u}{3^{ab}}
        \Witness{n}{u}{3^{abba}}
        \Witness{y}{u}{1^{abb}}
        \Witness{n}{u}{1^{abb}}
        \Witness{n}{u}{3^{ab}}
        \Witness{n}{u}{2^a}
        \Witness{n}{u}{3^{aab}}
        \Witness{n}{u}{2^{aba}}
        \Witness{n}{u}{3^{ab}}
        \Witness{n}{u}{2^\epsilon}
        \Witness{n}{u}{3^{ab}}
        \Witness{n}{u}{2^{ba}}
        \Witness{n}{u}{3^b}
    }
    \notag \\
    & =
    \Witness{y}{u}{2^{aba}}
    \Witness{y}{u}{3^{ab}}
    \Witness{y}{u}{1^{abb}}
    \Witness{n}{u}{1^{abb}}
    \Witness{n}{u}{2^\epsilon}
    \Witness{n}{u}{2^a}
    \Witness{n}{u}{2^{aba}}
    \Witness{n}{u}{2^{aba}}
    \Witness{n}{u}{2^{ba}}
    \Witness{n}{u}{3^{aab}}
    \Witness{n}{u}{3^{ab}}
    \Witness{n}{u}{3^{ab}}
    \Witness{n}{u}{3^{ab}}
    \Witness{n}{u}{3^{abba}}
    \Witness{n}{u}{3^b}
    \notag
\end{align}
and in $\Stal_2\Par{A^*}$, we have, up to isomorphism,
\begin{align}
    3^a 2^a 1^b 3^{ba} 3^\epsilon & \Superposition{
        2^a 1^{ab},
        3^b 3^\epsilon 2^{ab},
        1^{ba} 3^b,
        1^a 1^{ab}
    }
    \\
    & =
    \PSymbol_{\Equiv_2}\Par{
        \Witness{y}{u}{1^{aba}}
        \Witness{y}{u}{3^{ab}}
        \Witness{y}{u}{3^{ab}}
        \Witness{n}{u}{3^a}
        \Witness{y}{u}{2^{aab}}
        \Witness{y}{u}{2^{ba}}
        \Witness{y}{u}{1^{bab}}
        \Witness{n}{u}{1^{baba}}
        \Witness{n}{u}{3^{bab}}
        \Witness{n}{u}{1^{ba}}
        \Witness{n}{u}{3^b}
    }
    \notag \\
    & =
    \Witness{y}{u}{1^{aba}}
    \Witness{y}{u}{3^{ab}}
    \Witness{y}{u}{3^{ab}}
    \Witness{y}{u}{2^{aab}}
    \Witness{y}{u}{2^{ba}}
    \Witness{y}{u}{1^{bab}}
    \Witness{n}{u}{1^{ba}}
    \Witness{n}{u}{1^{baba}}
    \Witness{n}{u}{3^a}
    \Witness{n}{u}{3^b}
    \Witness{n}{u}{3^{bab}}.
    \notag
\end{align}

\subsubsection{Presentation}
In order to establish a presentation of $\Stal_k(\Monoid)$, we introduce an alternative
description of the clone congruence $\Equiv_k$ through a new equivalence relation
$\Equiv_k'$. For this, let us define $\Equiv_k'$ as the equivalence relation on
$\P(\Monoid)$ generated by
\begin{multline} \label{equ:alternative_relation_stal}
    \PWord \Conc
    \Par{\alpha_1 \MonoidProductExtension \QWord} \Conc \RWord_1 \Conc
    \Par{\alpha_2 \MonoidProductExtension \QWord} \Conc \RWord_2 \Conc
    \ \dots \ \Conc
    \Par{\alpha_k \MonoidProductExtension \QWord} \Conc
    \RWord_k \Conc
    \Par{\alpha_{k + 1} \MonoidProductExtension \QWord}
    \Conc \PWord'
    \\
    \Equiv_k'
    \PWord \Conc
    \Par{\alpha_1 \MonoidProductExtension \QWord} \Conc \RWord_1 \Conc
    \Par{\alpha_2 \MonoidProductExtension \QWord} \Conc \RWord_2 \Conc
    \ \dots \ \Conc
    \Par{\alpha_k \MonoidProductExtension \QWord} \Conc
    \Par{\alpha_{k + 1} \MonoidProductExtension \QWord}
    \Conc \RWord_k \Conc
    \PWord',
\end{multline}
where $\PWord, \PWord', \QWord, \RWord_1, \RWord_2, \dots, \RWord_k \in
\P(\Monoid)$ and $\alpha_1, \alpha_2, \dots, \alpha_k, \alpha_{k + 1} \in \Monoid$.

\begin{Statement}{Lemma}{lem:alternative_relation_stal}
    For any monoid $\Monoid$ and any $k \geq 0$, the binary relations $\Equiv_k$ and
    $\Equiv_k'$ on $\P(\Monoid)$ are equal.
\end{Statement}
\begin{Proof}
    Let $\Leq$ be any total order relation of $\Monoid$. Let us show that the left-hand side
    and the right-hand side of Equation~\eqref{equ:alternative_relation_stal} are
    $\Equiv_k$-equivalent. Observe first that the images by $\First_k$ of the left-hand side
    and the right-hand side of~\eqref{equ:alternative_relation_stal} are, by using the
    notations introduced in this equation, both equal to the image by $\First_k$ of the
    $\Monoid$-pigmented word
    \begin{math}
        \PWord \Conc
        \Par{\alpha_1 \MonoidProductExtension \QWord} \Conc \RWord_1 \Conc
        \Par{\alpha_2 \MonoidProductExtension \QWord} \Conc \RWord_2 \Conc
        \ \dots \ \Conc
        \Par{\alpha_k \MonoidProductExtension \QWord} \Conc
        \RWord_k \Conc
        \PWord'.
    \end{math}
    Moreover, since the $\Monoid$-pigmented words of the left-hand side and right-hand side
    of~\eqref{equ:alternative_relation_stal} differ only by moving some $\Monoid$-pigmented
    letters, their images by $\Sort_{\Leq}$ are equal. Therefore, since $\Equiv_k'$ is
    generated by~\eqref{equ:alternative_relation_stal}, this shows that $\Equiv_k'$ is
    contained in $\Equiv_k$.

    Conversely, let $\PWord, \PWord' \in \P(\Monoid)$ such that $\PWord \Equiv_k \PWord'$.
    By Proposition~\ref{prop:covering_relation_stal_equivalence}, this is equivalent to the
    fact that $\PWord \EquivCovering_k \PWord'$. Since $\EquivCovering_k$ is generated by
    $\Covering_k$, we have two cases to explore depending whether $\PWord \Covering^{(1)}_k
    \PWord'$ or $\PWord \Covering^{(2)}_k \PWord'$. These two cases are treated uniformly as
    follows. For any $j \in [2]$, if $\PWord \Covering^{(j)}_k \PWord'$, then $\PWord$ and
    $\PWord'$ decompose as
    \begin{math}
        \PWord
        =
        \QWord \Conc \Witness{n}{u}{i_1^{\alpha_1}}
        \Witness{u}{u}{i_2^{\alpha_2}} \Conc \QWord'
    \end{math}
    and
    \begin{math}
        \PWord'
        =
        \QWord \Conc \Witness{u}{u}{i_2^{\alpha_2}}
        \Witness{n}{u}{i_1^{\alpha_1}} \Conc \QWord'
    \end{math}
    where $\QWord, \QWord' \in \P(\Monoid)$, $i_1^{\alpha_1}, i_2^{\alpha_2} \in
    \SetLetters_\Monoid$, and $i_1^{\alpha_1} \ne i_2^{\alpha_2}$. Since the position
    $\Length(\QWord) + 1$ of $\PWord$ is not a left $k$-witness, there are necessarily at
    least $k$ occurrences of $\Monoid$-pigmented letters having $i_1$ as value in $\QWord$.
    Hence,
    \begin{align} \label{equ:alternative_relation_stal_1}
        \PWord
        & =
        \RWord_1 \Conc \Witness{y}{u}{i_1^{\beta_1}} \Conc \RWord_2 \Conc
        \Witness{y}{u}{i_1^{\beta_2}}
        \Conc \ \dots \ \Conc
        \RWord_k \Conc \Witness{y}{u}{i_1^{\beta_k}} \Conc \RWord_{k + 1}
        \Conc \Witness{n}{u}{i_1^{\alpha_1}}
        \Witness{u}{u}{i_2^{\alpha_2}} \Conc \QWord'
        \\
        & \Equiv_k'
        \RWord_1 \Conc \Witness{y}{u}{i_1^{\beta_1}} \Conc \RWord_2 \Conc
        \Witness{y}{u}{i_1^{\beta_2}}
        \Conc \ \dots \ \Conc
        \RWord_k \Conc \Witness{y}{u}{i_1^{\beta_k}} \Conc
        \Witness{n}{u}{i_1^{\alpha_1}} \Conc
        \RWord_{k + 1}
        \Conc
        \Witness{u}{u}{i_2^{\alpha_2}} \Conc \QWord'
        \notag
        \\
        & \Equiv_k'
        \RWord_1 \Conc \Witness{y}{u}{i_1^{\beta_1}} \Conc \RWord_2 \Conc
        \Witness{y}{u}{i_1^{\beta_2}}
        \Conc \ \dots \ \Conc
        \RWord_k \Conc \Witness{y}{u}{i_1^{\beta_k}} \Conc
        \RWord_{k + 1}
        \Conc
        \Witness{u}{u}{i_2^{\alpha_2}} \Conc
        \Witness{n}{u}{i_1^{\alpha_1}} \Conc
        \QWord'
        = \PWord'.
        \notag
    \end{align}
    where $\RWord_1, \RWord_2, \dots, \RWord_k, \RWord_{k + 1} \in \P(\Monoid)$ and
    $\beta_1, \beta_2, \dots, \beta_k \in \Monoid$. The first $\Equiv_k'$-equivalence
    of~\eqref{equ:alternative_relation_stal_1} is a consequence
    of~\eqref{equ:alternative_relation_stal} considered from left to right and the second
    $\Equiv_k'$-equivalence of~\eqref{equ:alternative_relation_stal_1} is a consequence
    of~\eqref{equ:alternative_relation_stal} considered from right to left. This shows that
    $\PWord \Equiv_k \PWord'$ implies $\PWord \Equiv_k' \PWord'$ and establishes the
    statement of the lemma.
\end{Proof}

\begin{Statement}{Theorem}{thm:presentation_stal}
    For any monoid $\Monoid$ and any $k \geq 0$, the clone $\Stal_k(\Monoid)$ admits the
    presentation $\Par{\GeneratingSet_\Monoid, \RelationSet_\Monoid'}$ where
    $\RelationSet_\Monoid'$ is the set $\RelationSet_\Monoid$ augmented with the
    $\GeneratingSet_\Monoid$-equation
    \begin{equation} \label{equ:presentation_stal}
        \RightComb_\Monoid\Par{
            1^{\alpha_1} 2^\MonoidUnit 1^{\alpha_2} 3^\MonoidUnit
            \dots 1^{\alpha_k} (k + 1)^\MonoidUnit 1^{\alpha_{k + 1}}
        }
        \ \RelationSet_\Monoid' \
        \RightComb_\Monoid\Par{
            1^{\alpha_1} 2^\MonoidUnit 1^{\alpha_2} 3^\MonoidUnit
            \dots 1^{\alpha_k} 1^{\alpha_{k + 1}} (k + 1)^\MonoidUnit
        },
    \end{equation}
    where $\alpha_1, \alpha_2, \dots, \alpha_k, \alpha_{k + 1} \in \Monoid$ and
    $\MonoidUnit$ is the unit of $\Monoid$.
\end{Statement}
\begin{Proof}
    Let $\Equiv''_k$ be the clone congruence of $\P(\Monoid)$ generated by
    \begin{equation} \label{equ:presentation_stal_1}
        1^{\alpha_1} 2^\MonoidUnit 1^{\alpha_2} 3^\MonoidUnit
        \dots 1^{\alpha_k} (k + 1)^\MonoidUnit 1^{\alpha_{k + 1}}
        \Equiv''_k
        1^{\alpha_1} 2^\MonoidUnit 1^{\alpha_2} 3^\MonoidUnit
        \dots 1^{\alpha_k} 1^{\alpha_{k + 1}} (k + 1)^\MonoidUnit
    \end{equation}
    with $\alpha_1, \alpha_2, \dots, \alpha_k, \alpha_{k + 1} \in \Monoid$. Let us show that
    the clone congruences $\Equiv_k$ and $\Equiv''_k$ of $\P(\Monoid)$ are equal. This will
    imply, by the remark stated in Section~\ref{subsubsec:remark_presentations}, that
    $\P(\Monoid) /_{\Equiv''_k} = \P(\Monoid) /_{\Equiv_k} = \Stal_k(\Monoid)$ admits the
    stated presentation.

    First, since the images by the map $\First_k$ of the left-hand side and the right-hand
    side of~\eqref{equ:presentation_stal_1} are both equal to
    \begin{math}
        1^{\alpha_1} 2^\MonoidUnit 1^{\alpha_2} 3^\MonoidUnit
        \dots 1^{\alpha_k} (k + 1)^\MonoidUnit
    \end{math}
    and since the images by $\Sort_{\Leq}$, where $\Leq$ is any total order relation on
    $\Monoid$, of the left-hand side and the right-hand side
    of~\eqref{equ:presentation_stal_1} are the same, we have
    \begin{equation}
        1^{\alpha_1} 2^\MonoidUnit 1^{\alpha_2} 3^\MonoidUnit
        \dots 1^{\alpha_k} (k + 1)^\MonoidUnit 1^{\alpha_{k + 1}}
        \Equiv_k
        1^{\alpha_1} 2^\MonoidUnit 1^{\alpha_2} 3^\MonoidUnit
        \dots 1^{\alpha_k} 1^{\alpha_{k + 1}} (k + 1)^\MonoidUnit.
    \end{equation}
    This shows that $\Equiv''_k$ is contained in $\Equiv_k$.

    To prove that $\Equiv_k$ is contained in $\Equiv''_k$, let us show that $\Equiv_k'$ is
    contained in $\Equiv''_k$. By Lemma~\ref{lem:alternative_relation_stal}, the targeted
    property will follow. For any $\PWord, \PWord', \QWord, \RWord_1, \RWord_2, \dots,
    \RWord_k \in \P(\Monoid)$ and $\alpha_1, \alpha_2, \dots, \alpha_k, \alpha_{k + 1} \in
    \Monoid$, we have
    \begin{align} \label{equ:presentation_stal_2}
        \PWord \Conc
        \Par{\alpha_1 \MonoidProductExtension \QWord} \Conc \RWord_1 \Conc
        &
        \Par{\alpha_2 \MonoidProductExtension \QWord} \Conc \RWord_2 \Conc
        \ \dots \ \Conc
        \Par{\alpha_k \MonoidProductExtension \QWord} \Conc
        \RWord_k \Conc
        \Par{\alpha_{k + 1} \MonoidProductExtension \QWord} \Conc
        \PWord'
        \\
        & =
        1^\MonoidUnit 2^\MonoidUnit 3^\MonoidUnit
        \Han{
            \PWord,
            1^{\alpha_1} 2^\MonoidUnit 1^{\alpha_2} 3^\MonoidUnit
            \dots 1^{\alpha_k} (k + 1)^\MonoidUnit 1^{\alpha_{k + 1}}
            \Han{
                \QWord, \RWord_1, \RWord_2, \dots, \RWord_k
            },
            \PWord'
        }
        \notag
        \\
        & \Equiv''_k
        1^\MonoidUnit 2^\MonoidUnit 3^\MonoidUnit
        \Han{
            \PWord,
            1^{\alpha_1} 2^\MonoidUnit 1^{\alpha_2} 3^\MonoidUnit
            \dots 1^{\alpha_k} 1^{\alpha_{k + 1}} (k + 1)^\MonoidUnit
            \Han{
                \QWord, \RWord_1, \RWord_2, \dots, \RWord_k
            },
            \PWord'
        }
        \notag
        \\
        & =
        \PWord \Conc
        \Par{\alpha_1 \MonoidProductExtension \QWord} \Conc \RWord_1 \Conc
        \Par{\alpha_2 \MonoidProductExtension \QWord} \Conc \RWord_2 \Conc
        \ \dots \ \Conc
        \Par{\alpha_k \MonoidProductExtension \QWord}
        \Conc
        \Par{\alpha_{k + 1} \MonoidProductExtension \QWord} \Conc
        \RWord_k \Conc
        \PWord'
        \notag
    \end{align}
    so that the first and last members of~\eqref{equ:presentation_stal_2} are
    $\Equiv''_k$-equivalent. Since $\Equiv_k'$ is the equivalence relation generated
    by~\eqref{equ:alternative_relation_stal}, the targeted property is shown. This
    establishes the statement of the theorem.
\end{Proof}

By Theorem~\ref{thm:presentation_stal}, any $\Stal_k(\Monoid)$-algebra is, up to term
equivalence, an $\Monoid$-pigmented monoid $\Par{\Algebra, \Product, \Identity,
\Par{\Pigmentation_\alpha}_{\alpha \in \Monoid}}$ where $\Product$ and
$\Par{\Pigmentation_\alpha}_{\alpha \in \Monoid}$, satisfy, by spelling
out~\eqref{equ:presentation_stal} and simplifying it modulo the background theory
$\Equiv_{\RelationSet_{\Monoid}}$,
\begin{multline}
    \Pigmentation_{\alpha_1} \Par{x_1} \Product x_2 \Product
    \Pigmentation_{\alpha_2} \Par{x_1} \Product x_3 \Product
    \dots \Product
    \Pigmentation_{\alpha_k} \Par{x_1} \Product x_{k + 1} \Product
    \Pigmentation_{\alpha_{k + 1}} \Par{x_1}
    \\ =
    \Pigmentation_{\alpha_1} \Par{x_1} \Product x_2 \Product
    \Pigmentation_{\alpha_2} \Par{x_1} \Product x_3 \Product
    \dots \Product
    \Pigmentation_{\alpha_k} \Par{x_1} \Product
    \Pigmentation_{\alpha_{k + 1}} \Par{x_1} \Product
    x_{k + 1}
\end{multline}
for any $x_1, x_2, x_3, \dots, x_{k + 1} \in \Algebra$ and $\alpha_1, \alpha_2, \dots,
\alpha_k, \alpha_{k + 1} \in \Monoid$.

As a side remark, the equivalence relation $\Equiv_1$, as a monoid congruence, has been
introduced in~\cite{HNT08} under the name of the ``stalactic congruence''. As a monoid
congruence, $\Equiv_k$, $k \geq 0$, is therefore a generalization of the previous one.

\subsection{On pigmented pillars} \label{subsec:pill}
By considering the intersection of the clone congruences $\Equiv_\Sort$,
$\Equiv_{\First_k}$, $k \geq 0$, and their reversions $\Equiv_{\First_{k'}}^\Reverse$, $k'
\geq 0$, we construct a quotient clone $\Pill_{k, k'}(\Monoid)$ of $\P(\Monoid)$. This clone
is studied in detail for the case $k = 1 = k'$. A description through new combinatorial
objects named $\Monoid$-pigmented pillars is introduced and a finitely equationally
axiomatizable presentation is described. These results are based on the introduction of a
$\PSymbol$-symbol for the underlying equivalence relation.

\subsubsection{Clone construction}
For any parameters $k, k' \geq 0$, let $\Equiv_{k, k'}$ be the clone congruence
$\Equiv_{\First_k} \cap \Equiv_\Sort \cap \Equiv^\Reverse_{\First_{k'}}$ and
\begin{equation}
    \Pill_{k, k'}(\Monoid) := \P(\Monoid) /_{\Equiv_{k, k'}}.
\end{equation}
By Propositions~\ref{prop:congruence_sort}, \ref{prop:congruence_first},
and~\ref{prop:reversed_congruence}, $\Pill_{k, k'}(\Monoid)$ is a well-defined clone, and
$\Stal_k(\Monoid)$, $\Magn_{k, k'}(\Monoid)$, and $\Stal_{k'}^\Reverse(\Monoid)$ are
isomorphic to quotients of $\Pill_{k, k'}$. Since for any $0 \leq k \leq k''$ and $0 \leq k'
\leq k'''$, $\Equiv_{k'', k'''}$ is a refinement of $\Equiv_{k, k'}$, $\Pill_{k,
k'}(\Monoid)$ is isomorphic to a quotient of $\Pill_{k'', k'''}(\Monoid)$. Moreover, since
$\Equiv_{0, 0}$ and $\Equiv_{\Sort}$ are the same equivalence relations, $\Pill_{0,
0}(\Monoid)$ is identical to $\WInc(\Monoid)$. Besides, the clone
\begin{math}
    \Pill_{k, k'}^\Reverse(\Monoid)
    := \Pill_{k, k'}(\Monoid)^\Reverse
    = \P(\Monoid) /_{\Equiv_{k, k'}^\Reverse}
\end{math}
is by Proposition~\ref{prop:reversed_congruence} isomorphic to $\Pill_{k, k'}(\Monoid)$.
Since the reversion operation on congruences is involutive, $\Equiv_{k, k'}^\Reverse =
\Equiv_{k', k}$, the clones $\Pill_{k, k'}^\Reverse(\Monoid)$ and $\Pill_{k', k}(\Monoid)$
are identical and $\Pill_{k, k'}(\Monoid)$ and $\Pill_{k', k}(\Monoid)$ are isomorphic.

\subsubsection{Equivalence relation}
To lighten the notation, we denote by $\Equiv$ the equivalence relation $\Equiv_{1, 1}$ on
$\P(\Monoid)$. By definition, for any $\PWord, \PWord' \in \P(\Monoid)$, $\PWord \Equiv
\PWord'$ holds if and only if
\begin{math}
    \Par{\First_1\Par{\PWord}, \Sort_{\Leq}\Par{\PWord}, \First_1^\Reverse\Par{\PWord}}
    = \Par{\First_1\Par{\PWord'}, \Sort_{\Leq}\Par{\PWord'}, \First_1^\Reverse\Par{\PWord'}}
\end{math}
where $\Leq$ is any total order relation on~$\Monoid$.

In order to obtain properties about the clone $\Pill_{1, 1}(\Monoid)$, we introduce an
alternative equivalence relation $\EquivCovering$ for which we will show that it is equal to
$\Equiv$. Let $\Covering^{(1)}$, $\Covering^{(2)}$, and $\Covering^{(3)}$ be the three
binary relations on $\P(\Monoid)$ satisfying
\begin{equation} \label{equ:covering_relation_pill_1}
    \PWord \Conc \Witness{u}{n}{i^{\alpha_1}} \Conc \QWord \Conc
        \Witness{n}{n}{i^{\alpha_2}} \Conc \PWord'
    \enspace \Covering^{(1)} \enspace
    \PWord \Conc \Witness{u}{n}{i^{\alpha_1}} \Conc
        \Witness{n}{n}{i^{\alpha_2}} \Conc \QWord \Conc \PWord'
    \qquad
    \text{ where }
    \QWord \ne \PWordEmpty,
    \text{ and } i \notin \QWord,
\end{equation}
\begin{equation} \label{equ:covering_relation_pill_2}
    \PWord \Conc \Witness{n}{n}{i^{\alpha_1}} \Witness{n}{n}{i^{\alpha_2}} \Conc \PWord'
    \enspace \Covering^{(2)} \enspace
    \PWord \Conc \Witness{n}{n}{i^{\alpha_2}} \Witness{n}{n}{i^{\alpha_1}} \Conc \PWord'
    \qquad
    \text{ where }
    \alpha_1 \ne \alpha_2
    \text{ and } \alpha_2 \Leq \alpha_1,
\end{equation}
\begin{equation} \label{equ:covering_relation_pill_3}
    \PWord \Conc \Witness{n}{y}{i_1^{\alpha_1}} \Witness{u}{n}{i_2^{\alpha_2}}
        \Conc \PWord'
    \enspace \Covering^{(3)} \enspace
    \PWord \Conc \Witness{u}{n}{i_2^{\alpha_2}} \Witness{n}{y}{i_1^{\alpha_1}}
        \Conc \PWord'
    \qquad
    \text{ where }
    i_1 \ne i_2,
\end{equation}
where $\PWord, \PWord', \QWord \in \P(\Monoid)$, $i^{\alpha_1}, i^{\alpha_2},
i_1^{\alpha_1}, i_2^{\alpha_2} \in \SetLetters_\Monoid$, and where the notation $i \notin
\RWord$ means that the $\Monoid$-pigmented word $\RWord$ has no occurrence of any
$\Monoid$-pigmented letter having $i$ as value. Let $\CoveringRT^{(j)}$, $j \in [3]$, be the
reflexive and transitive closure of $\Covering^{(j)}$, $\Covering$ be the union
$\Covering^{(1)} \cup \Covering^{(2)} \cup \Covering^{(3)}$, and $\EquivCovering$ be the
reflexive, symmetric, and transitive closure of~$\Covering$.

\begin{Statement}{Lemma}{lem:covering_relation_algorithm_pill}
    For any monoid $\Monoid$ endowed with a total order relation $\Leq$ and any $\PWord \in
    \P(\Monoid)$, we have
    \begin{math}
        \PWord \CoveringRT^{(1)} \PWord^{(1)} \CoveringRT^{(2)} \PWord^{(1, 2)}
        \CoveringRT^{(3)} \PWord^{(1, 2, 3)}
    \end{math}
    where
    \begin{enumerate}
        \item the $\Monoid$-pigmented word $\PWord^{(1)}$ is obtained from $\PWord$ by the
        following process. For any value $i$ between $1$ and the maximal value appearing in
        $\PWord$, extract the subword of $\PWord$ consisting of $\Monoid$-pigmented letters
        that have $i$ as their value and whose positions are neither left nor right
        $1$-witnesses, and place it to the right of the letter of value $i$ whose position
        is a left $1$-witness;
        \item the $\Monoid$-pigmented word $\PWord^{(1, 2)}$ is obtained from $\PWord^{(1)}$
        by sorting w.r.t.\ the total order relation~$\PLeq$ each factor consisting of
        letters having the same value and whose positions are neither left nor right
        $1$-witnesses;
        \item the $\Monoid$-pigmented word $\PWord^{(1, 2, 3)}$ is obtained from
        $\PWord^{(1, 2)}$ by the following process. Consider, from right to left, each
        $\Monoid$-pigmented letter whose position is a right $1$-witness but not a left
        $1$-witness. For each such letter, except the rightmost one, extract it and insert
        it to the right within the word, placing it immediately to the left of the first
        $\Monoid$-pigmented letter encountered to its right whose position is a right
        $1$-witness.
    \end{enumerate}
\end{Statement}
\begin{Proof}
    First, observe that applying the relation $\Covering^{(1)}$ to an $\Monoid$-pigmented
    word consists of moving an $\Monoid$-pigmented letter $i^{\alpha_2}$ whose position is
    neither a left nor a right $1$-witness to the right of the nearest letter of the same
    value, provided that the position of this letter is not a right $1$-witness. This last
    condition about not being a right $1$-witness obviously always holds. Iterating this
    operation moves each subword, consisting of $\Monoid$-pigmented letters of the same
    value whose positions are neither left nor right $1$-witnesses, to the right of the
    letter having $i$ as value and whose position is a left $1$-witness. Such a letter
    exists because, as the position of $i^{\alpha_2}$ is never a left $1$-witness during all
    this process, there is necessarily a letter at the left of $i^{\alpha_2}$ having $i$ as
    value and whose position is a left $1$-witness. Consequently, the $\Monoid$-pigmented
    word $\PWord^{(1)}$ is eventually obtained from $\PWord$ through this process. Hence, we
    conclude that $\PWord \CoveringRT^{(1)} \PWord^{(1)}$, as expected.

    Besides, the application of relation $\Covering^{(2)}$ on an $\Monoid$-pigmented word
    has the effect of transposing two adjacent $\Monoid$-pigmented letters having the same
    value $i$ and whose positions are neither left nor right $1$-witnesses, in such a way
    that they become sorted w.r.t.\ the order relation $\PLeq$. The iteration of such an
    operation will sort, w.r.t.\ $\PLeq$, each factor consisting of letters having $i$ as
    value and whose positions are neither left nor right $1$-witnesses. Therefore, the
    $\Monoid$-pigmented word $\PWord^{(1,2)}$ is eventually obtained through this process
    from $\PWord^{(1)}$. Thus, we have $\PWord^{(1)} \CoveringRT^{(2)} \PWord^{(1,2)}$, as
    expected.

    Finally, the application of relation $\Covering^{(3)}$ on an $\Monoid$-pigmented word
    consists of transposing an $\Monoid$-pigmented letter $i_1^{\alpha_1}$, whose position
    is a right $1$-witness but not a left $1$-witness, with the neighboring letter
    $i_2^{\alpha_2}$ to its right, provided that $i_1 \ne i_2$ and the position of
    $i_2^{\alpha_2}$ is not a right $1$-witness. The iteration of such an operation will
    move this letter $i_1^{\alpha_1}$ to the right of the word, placing it on the left of
    the nearest $\Monoid$-pigmented letter whose position is a right $1$-witness. Therefore,
    the $\Monoid$-pigmented word $\PWord^{(1, 2, 3)}$ is eventually obtained through this
    process from $\PWord^{(1, 2)}$. Thus, we have $\PWord^{(1, 2)} \CoveringRT^{(3)}
    \PWord^{(1, 2, 3)}$ as expected.
\end{Proof}

For instance, in $\P\Par{A^*}$, where $A^*$ is the free monoid over $\{a, b\}$, by setting
\begin{equation}
    \PWord :=
    \Witness{y}{n}{2^{ab}} \Witness{n}{n}{2^a} \Witness{y}{n}{4^b} \Witness{n}{n}{4^b}
    \Witness{n}{n}{2^\epsilon} \Witness{n}{n}{4^{ab}} \Witness{n}{y}{4^\epsilon}
    \Witness{y}{n}{3^a} \Witness{n}{n}{3^a} \Witness{n}{n}{3^{ba}} \Witness{n}{y}{2^{ab}}
    \Witness{y}{y}{5^b} \Witness{n}{y}{3^{ab}},
\end{equation}
we have
\begin{equation}
    \PWord^{(1)} =
    \Witness{y}{n}{2^{ab}} \Witness{n}{n}{2^a} \Witness{n}{n}{2^\epsilon}
    \Witness{y}{n}{4^b} \Witness{n}{n}{4^b} \Witness{n}{n}{4^{ab}}
    \Witness{n}{y}{4^\epsilon} \Witness{y}{n}{3^a} \Witness{n}{n}{3^a}
    \Witness{n}{n}{3^{ba}} \Witness{n}{y}{2^{ab}} \Witness{y}{y}{5^b}
    \Witness{n}{y}{3^{ab}},
\end{equation}
\begin{equation}
    \PWord^{(1, 2)} =
    \Witness{y}{n}{2^{ab}} \Witness{n}{n}{2^\epsilon} \Witness{n}{n}{2^a}
    \Witness{y}{n}{4^b} \Witness{n}{n}{4^{ab}} \Witness{n}{n}{4^b}
    \Witness{n}{y}{4^\epsilon} \Witness{y}{n}{3^a} \Witness{n}{n}{3^a}
    \Witness{n}{n}{3^{ba}} \Witness{n}{y}{2^{ab}} \Witness{y}{y}{5^b}
    \Witness{n}{y}{3^{ab}},
\end{equation}
and
\begin{equation}
    \PWord^{(1, 2, 3)} =
    \Witness{y}{n}{2^{ab}} \Witness{n}{n}{2^\epsilon} \Witness{n}{n}{2^a}
    \Witness{y}{n}{4^b} \Witness{n}{n}{4^{ab}} \Witness{n}{n}{4^b} \Witness{y}{n}{3^a}
    \Witness{n}{n}{3^a} \Witness{n}{n}{3^{ba}} \Witness{n}{y}{4^\epsilon}
    \Witness{n}{y}{2^{ab}} \Witness{y}{y}{5^b} \Witness{n}{y}{3^{ab}}.
\end{equation}
Moreover, by setting
\begin{equation}
    \PWord :=
    \Witness{y}{n}{1^{c}} \Witness{n}{n}{1^{c}} \Witness{y}{n}{2^{a}} \Witness{n}{y}{1^{a}}
    \Witness{y}{n}{4^{b}} \Witness{y}{n}{6^{c}} \Witness{n}{n}{6^{b}} \Witness{y}{n}{3^{b}}
    \Witness{n}{y}{2^{b}} \Witness{n}{n}{4^{c}} \Witness{y}{y}{5^{a}} \Witness{n}{n}{4^{b}}
    \Witness{n}{n}{6^{a}} \Witness{n}{n}{3^{a}} \Witness{n}{n}{3^{b}} \Witness{n}{y}{3^{c}}
    \Witness{n}{y}{4^{a}} \Witness{n}{y}{6^{a}},
\end{equation}
we have
\begin{equation}
    \PWord^{(1)} =
    \Witness{y}{n}{1^{c}} \Witness{n}{n}{1^{c}} \Witness{y}{n}{2^{a}} \Witness{n}{y}{1^{a}}
    \Witness{y}{n}{4^{b}} \Witness{n}{n}{4^{c}} \Witness{n}{n}{4^{b}} \Witness{y}{n}{6^{c}}
    \Witness{n}{n}{6^{b}} \Witness{n}{n}{6^{a}} \Witness{y}{n}{3^{b}} \Witness{n}{n}{3^{a}}
    \Witness{n}{n}{3^{b}} \Witness{n}{y}{2^{b}} \Witness{y}{y}{5^{a}} \Witness{n}{y}{3^{c}}
    \Witness{n}{y}{4^{a}} \Witness{n}{y}{6^{a}},
\end{equation}
\begin{equation}
    \PWord^{(1, 2)} =
    \Witness{y}{n}{1^{c}} \Witness{n}{n}{1^{c}} \Witness{y}{n}{2^{a}} \Witness{n}{y}{1^{a}}
    \Witness{y}{n}{4^{b}} \Witness{n}{n}{4^{b}} \Witness{n}{n}{4^{c}} \Witness{y}{n}{6^{c}}
    \Witness{n}{n}{6^{a}} \Witness{n}{n}{6^{b}} \Witness{y}{n}{3^{b}} \Witness{n}{n}{3^{a}}
    \Witness{n}{n}{3^{b}} \Witness{n}{y}{2^{b}} \Witness{y}{y}{5^{a}} \Witness{n}{y}{3^{c}}
    \Witness{n}{y}{4^{a}} \Witness{n}{y}{6^{a}},
\end{equation}
and
\begin{equation}
    \PWord^{(1, 2, 3)} =
    \Witness{y}{n}{1^{c}} \Witness{n}{n}{1^{c}} \Witness{y}{n}{2^{a}} \Witness{y}{n}{4^{b}}
    \Witness{n}{n}{4^{b}} \Witness{n}{n}{4^{c}} \Witness{y}{n}{6^{c}} \Witness{n}{n}{6^{a}}
    \Witness{n}{n}{6^{b}} \Witness{y}{n}{3^{b}} \Witness{n}{n}{3^{a}} \Witness{n}{n}{3^{b}}
    \Witness{n}{y}{1^{a}} \Witness{n}{y}{2^{b}} \Witness{y}{y}{5^{a}} \Witness{n}{y}{3^{c}}
    \Witness{n}{y}{4^{a}} \Witness{n}{y}{6^{a}}.
\end{equation}

\begin{Statement}{Proposition}{prop:covering_relation_pill_equivalence}
    For any monoid $\Monoid$, the binary relations $\Equiv$ and $\EquivCovering$ on
    $\P(\Monoid)$ are equal.
\end{Statement}
\begin{Proof}
    Let $\Leq$ be any total order relation on $\Monoid$.

    Let us show that $\PWord \Covering \PWord'$ for $\PWord, \PWord' \in \P(\Monoid)$
    entails $\PWord \Equiv \PWord'$. Since the equivalence relation $\EquivCovering$ is
    generated by $\Covering$, this will show that $\EquivCovering$ is contained in $\Equiv$.
    We are going to consider the following three cases depending whether $\PWord
    \Covering^{(1)} \PWord'$, $\PWord \Covering^{(2)} \PWord'$, or $\PWord \Covering^{(3)}
    \PWord'$.
    \begin{enumerate}[label={\it (\Roman*)}]
        \item Assume that $\PWord \Covering^{(1)} \PWord'$. Hence, $\PWord$ and $\PWord'$
        decompose as
        \begin{math}
            \PWord = \RWord \Conc \Witness{u}{n}{i^{\alpha_1}} \Conc \QWord \Conc
            \Witness{n}{n}{i^{\alpha_2}} \Conc \RWord'
        \end{math}
        and
        \begin{math}
            \PWord' = \RWord \Conc \Witness{u}{n}{i^{\alpha_1}} \Conc
            \Witness{n}{n}{i^{\alpha_2}} \Conc \QWord \Conc \RWord'
        \end{math}
        where $\QWord, \RWord, \RWord' \in \P(\Monoid)$, $i^{\alpha_1}, i^{\alpha_2} \in
        \SetLetters_\Monoid$, $\QWord \ne \PWordEmpty$, and $i \notin \QWord$. Since
        \begin{equation}
            \First_1(\PWord)
            = \First_1\Par{
                \RWord \Conc \Witness{u}{n}{i^{\alpha_1}} \Conc \QWord \Conc
                \RWord'
            }
            = \First_1\Par{\PWord'},
        \end{equation}
        \begin{equation}
            \Sort_{\Leq}(\PWord)
            =
            \Sort_{\Leq}\Par{
                \RWord \Conc \Witness{u}{n}{i^{\alpha_1}} \Conc
                \Witness{n}{n}{i^{\alpha_2}} \Conc \QWord \Conc \RWord'
            }
            = \Sort_{\Leq}\Par{\PWord'},
        \end{equation}
        and
        \begin{equation}
            \First_1^\Reverse(\PWord)
            = \First_1^\Reverse\Par{\RWord \Conc \QWord \Conc \RWord'}
            = \First_1^\Reverse\Par{\PWord'},
        \end{equation}
        by definition of $\Equiv$, $\PWord \Equiv \PWord'$.
        \item Assume that $\PWord \Covering^{(2)} \PWord'$. Hence, $\PWord$ and $\PWord'$
        decompose as
        \begin{math}
            \PWord =
            \QWord \Conc \Witness{n}{n}{i^{\alpha_1}} \Witness{n}{n}{i^{\alpha_2}} \Conc
            \QWord'
        \end{math}
        and
        \begin{math}
            \PWord' =
            \QWord \Conc \Witness{n}{n}{i^{\alpha_2}} \Witness{n}{n}{i^{\alpha_1}} \Conc
            \QWord'
        \end{math}
        where $\QWord, \QWord' \in \P(\Monoid)$, $i^{\alpha_1}, i^{\alpha_2} \in
        \SetLetters_\Monoid$, $\alpha_1 \ne \alpha_2$, and $\alpha_2 \Leq \alpha_1$. Since
        \begin{equation}
            \First_1(\PWord) = \First_1\Par{\QWord \Conc \QWord'} = \First_1\Par{\PWord'},
        \end{equation}
        \begin{equation}
            \Sort_{\Leq}(\PWord)
            =
            \Sort_{\Leq}\Par{
                \QWord \Conc \Witness{n}{n}{i^{\alpha_2}} \Conc
                \Witness{n}{n}{i^{\alpha_1}} \Conc \QWord'
            }
            = \Sort_{\Leq}\Par{\PWord'},
        \end{equation}
        and
        \begin{equation}
            \First_1^\Reverse(\PWord) =
            \First_1\Par{\QWord \Conc \QWord'} =
            \First_1^\Reverse\Par{\PWord'},
        \end{equation}
        by definition of $\Equiv$, $\PWord \Equiv \PWord'$.
        \item Assume that $\PWord \Covering^{(3)} \PWord'$. Hence, $\PWord$ and $\PWord'$
        decompose as
        \begin{math}
            \PWord =
            \QWord \Conc \Witness{n}{y}{i_1^{\alpha_1}} \Witness{u}{n}{i_2^{\alpha_2}} \Conc
            \QWord'
        \end{math}
        and
        \begin{math}
            \PWord' =
            \QWord \Conc \Witness{u}{n}{i_2^{\alpha_2}} \Witness{n}{y}{i_1^{\alpha_1}} \Conc
            \QWord'
        \end{math}
        where $\QWord, \QWord' \in \P(\Monoid)$, $i_1^{\alpha_1}, i_2^{\alpha_2} \in
        \SetLetters_\Monoid$, and $i_1 \ne i_2$. Since
        \begin{equation}
            \First_1(\PWord)
            = \First_1\Par{\QWord \Conc \Witness{u}{n}{i_2^{\alpha_2}} \Conc \QWord'}
            = \First_1\Par{\PWord'},
        \end{equation}
        \begin{equation}
            \Sort_{\Leq}(\PWord)
            =
            \Sort_{\Leq}\Par{
                \QWord \Conc \Witness{u}{n}{i_2^{\alpha_2}} \Witness{n}{y}{i_1^{\alpha_1}}
                \Conc \QWord'
            }
            = \Sort_{\Leq}\Par{\PWord'},
        \end{equation}
        and
        \begin{equation}
            \First_1^\Reverse(\PWord)
            = \First_1^\Reverse\Par{
                \QWord \Conc \Witness{n}{y}{i_1^{\alpha_1}} \Conc \QWord'}
            = \First_1^\Reverse\Par{\PWord'},
        \end{equation}
        by definition of $\Equiv$, $\PWord \Equiv \PWord'$.
    \end{enumerate}

    Conversely, assume that $\PWord$ and $\PWord'$ are two $\Monoid$-pigmented words such
    that $\PWord \Equiv \PWord'$, Hence, by definition of $\Equiv$, $\First_1\Par{\PWord} =
    \First_1\Par{\PWord'}$, $\Sort_{\Leq}\Par{\PWord} = \Sort_{\Leq}\Par{\PWord'}$, and
    $\First_1^\Reverse\Par{\PWord} = \First_1^\Reverse\Par{\PWord'}$. By
    Lemma~\ref{lem:covering_relation_algorithm_pill} and by using the notations defined in
    its statement, the three previous properties imply that $\PWord^{(1, 2, 3)} =
    {\PWord'}^{(1, 2, 3)}$. Indeed, the fact that $\First_1\Par{\PWord} =
    \First_1\Par{\PWord'}$ and $\Sort_{\Leq}\Par{\PWord} = \Sort_{\Leq}\Par{\PWord'}$ lead
    to the fact that both $\PWord^{(1, 2)} = {\PWord'}^{(1, 2)}$ contains the same factors
    of letters with identical value and whose positions are not right $1$-witnesses.
    Moreover, the fact that $\First_1^\Reverse\Par{\PWord} = \First_1^\Reverse\Par{\PWord'}$
    leads to the fact that the two $\Monoid$-pigmented words $\PWord^{(1, 2, 3)}$ and
    ${\PWord'}^{(1, 2, 3)}$, obtained respectively from $\PWord^{(1, 2)}$ and
    ${\PWord'}^{(1, 2)}$, are such that each letter whose position is a right $1$-witnesses
    appears at the same location in both words. By definition of the equivalence relation
    $\Equiv$, this implies that $\PWord \EquivCovering \PWord'$.
\end{Proof}

\subsubsection{$\PSymbol$-symbol algorithm}
With the aim of describing of $\Pill_{1, 1}(\Monoid)$, we propose now a $\PSymbol$-symbol
for~$\Equiv$. Let $\PSymbol_{\Equiv} : \P(\Monoid) \to \P(\Monoid)$ be the map defined for
any $\PWord \in \P(\Monoid)$ by $\PSymbol_{\Equiv}(\PWord) := \PWord^{(1, 2, 3)}$ where
$\PWord^{(1, 2, 3)}$ is the $\Monoid$-pigmented word built from $\PWord$ as described in the
statement of Lemma~\ref{lem:covering_relation_algorithm_pill}.

\begin{Statement}{Lemma}{lem:p_symbol_pill_equiv_p_symbol}
    For any monoid $\Monoid$ and any $\PWord \in \P(\Monoid)$, $\PWord \EquivCovering
    \PSymbol_{\Equiv}(\PWord)$.
\end{Statement}
\begin{Proof}
    By definition of $\EquivCovering$, $\CoveringRT^{(1)}$, $\CoveringRT^{(2)}$, and
    $\CoveringRT^{(3)}$, it follows that both $\CoveringRT^{(1)}$, $\CoveringRT^{(2)}$, and
    $\CoveringRT^{(3)}$ are contained in $\EquivCovering$. The statement of the lemma is now
    a direct consequence of Lemma~\ref{lem:covering_relation_algorithm_pill} and of the
    definition of the map~$\PSymbol_{\Equiv}$.
\end{Proof}

\begin{Statement}{Lemma}{lem:covering_relation_pill_p_symbol}
    For any monoid $\Monoid$ and any $\PWord, \PWord' \in \P(\Monoid)$, $\PWord
    \EquivCovering \PWord'$ implies $\PSymbol_{\Equiv}(\PWord) =
    \PSymbol_{\Equiv}\Par{\PWord'}$.
\end{Statement}
\begin{Proof}
    It is straightforward to see that for any $j \in [3]$, if $\PWord \Covering^{(j)}
    \PWord'$, then $\PSymbol_{\Equiv}(\PWord) = \PSymbol_{\Equiv}\Par{\PWord'}$. Indeed,
    from the description of the computation steps of these two $\PSymbol$-symbols provided
    by Lemma~\ref{lem:covering_relation_algorithm_pill}, any change performed by the
    relation $\Covering^{(j)}$ on an $\Monoid$-pigmented word does not influence the
    final result. Since the equivalence relation $\EquivCovering$ is generated by the union
    of $\Covering^{(1)}$, $\Covering^{(2)}$, and $\Covering^{(3)}$, the statement of the
    lemma follows.
\end{Proof}

By Proposition~\ref{prop:covering_relation_pill_equivalence} and
Lemmas~\ref{lem:p_symbol_pill_equiv_p_symbol} and~\ref{lem:covering_relation_pill_p_symbol},
$\PSymbol_{\Equiv}$ is a $\PSymbol$-symbol for $\Equiv$.

\subsubsection{Description}
An \Def{$\Monoid$-pigmented pillar} (or simply \Def{pigmented pillar} when the context is
clear) of arity $n \geq 0$ is an $\Monoid$-pigmented word $\PWord$ of arity $n$ which is in
the image of $\PSymbol_{\Equiv}$. For instance, in $\P\Par{A^*}$, where $A^*$ is the free
monoid over $\Bra{a, b}$,
\begin{equation}
    \Witness{y}{n}{2^b} \Witness{y}{n}{1^a} \Witness{n}{n}{2^{ba}}
    \Witness{y}{n}{4^\epsilon} \Witness{n}{y}{1^{ba}} \Witness{y}{y}{3^a}
    \Witness{n}{y}{2^a} \Witness{n}{y}{4^b}
    \quad \text{and} \quad
    \Witness{y}{y}{1^\epsilon} \Witness{y}{n}{2^b} \Witness{y}{n}{6^b} \Witness{n}{n}{6^a}
    \Witness{n}{n}{2^{ba}} \Witness{y}{n}{5^{ab}} \Witness{n}{y}{2^{ab}}
    \Witness{n}{y}{6^\epsilon} \Witness{n}{n}{5^a} \Witness{n}{y}{5^{ab}}
\end{equation}
are not $A^*$-pigmented pillars. In contrast,
\begin{equation}
    \Witness{y}{n}{1^\epsilon} \Witness{y}{n}{2^\epsilon} \Witness{n}{n}{2^\epsilon}
    \Witness{y}{n}{4^a} \Witness{n}{n}{4^b} \Witness{y}{y}{3^{ab}} \Witness{y}{n}{2^a}
    \Witness{n}{y}{4^\epsilon} \Witness{n}{y}{1^\epsilon}
    \quad \text{and} \quad
    \Witness{y}{n}{4^b} \Witness{n}{n}{4^b} \Witness{y}{n}{5^{ba}} \Witness{y}{n}{3^{ba}}
    \Witness{n}{n}{3^{ab}} \Witness{y}{y}{1^{ba}} \Witness{n}{y}{5^\epsilon}
    \Witness{n}{y}{3^b} \Witness{y}{y}{6^b} \Witness{n}{y}{4^b}
\end{equation}
are $A^*$-pigmented pillars.

\begin{Statement}{Theorem}{thm:p_symbol_pill}
    For any monoid $\Monoid$, $\PSymbol_{\Equiv}$ is a $\PSymbol$-symbol for $\Equiv$ and
    $\PSymbol_{\Equiv}(\P(\Monoid))$ is the set of $\Monoid$-pigmented pillars.
    Moreover, the graded set $\Pill_{1, 1}(\Monoid)$ is isomorphic to the graded set of
    $\Monoid$-pigmented pillars.
\end{Statement}
\begin{Proof}
    By Proposition~\ref{prop:covering_relation_pill_equivalence} and
    Lemmas~\ref{lem:p_symbol_pill_equiv_p_symbol}
    and~\ref{lem:covering_relation_pill_p_symbol}, $\PSymbol_{\Equiv}$ is a
    $\PSymbol$-symbol for $\Equiv$. Moreover, the set of $\Monoid$-pigmented pillars is
    defined as the set $\PSymbol_{\Equiv}(\P(\Monoid))$. The last part of the statement is a
    direct implication of Proposition~\ref{prop:clone_realization_p_symbol} and the fact
    that $\PSymbol_{\Equiv}$ is, as we have just shown, a $\PSymbol$-symbol for~$\Equiv$.
\end{Proof}

By Proposition~\ref{prop:clone_realization_p_symbol} and Theorem~\ref{thm:p_symbol_pill},
$\Pill_{1, 1}(\Monoid)$ can be seen as a clone on $\Monoid$-pigmented pillars with
superposition maps satisfying~\eqref{equ:superposition_p_symbol}. For instance, in
$\Pill_{1, 1}\Par{A^*}$, where $A^*$ is the free monoid over $\Bra{a, b}$, we have, up to
isomorphism,
\begin{align}
    3^\epsilon 2^{ab} 1^b 1^a 4^a
    & \Superposition{
        2^{ba} 2^{ba} 1^{ab} 1^\epsilon,
        2^a 3^a,
        1^{ba},
        3^{ba} 3^a 1^{ab} 2^{ab} 1^b}
    \\
    & = \PSymbol_{\Equiv}\Par{
        \Witness{y}{n}{1^{ba}} \Witness{y}{n}{2^{aba}} \Witness{y}{n}{3^{aba}}
        \Witness{n}{n}{2^{bba}} \Witness{n}{n}{2^{bba}} \Witness{n}{n}{1^{bab}}
        \Witness{n}{n}{1^b} \Witness{n}{n}{2^{aba}} \Witness{n}{n}{2^{aba}}
        \Witness{n}{n}{1^{aab}} \Witness{n}{n}{1^a} \Witness{n}{n}{3^{aba}}
        \Witness{n}{y}{3^{aa}} \Witness{n}{n}{1^{aab}} \Witness{n}{y}{2^{aab}}
        \Witness{n}{y}{1^{ab}}
    }
    \notag \\
    & =
    \Witness{y}{n}{1^{ba}} \Witness{n}{n}{1^a} \Witness{n}{n}{1^{aab}}
    \Witness{n}{n}{1^{aab}} \Witness{n}{n}{1^b} \Witness{n}{n}{1^{bab}}
    \Witness{y}{n}{2^{aba}} \Witness{n}{n}{2^{aba}} \Witness{n}{n}{2^{aba}}
    \Witness{n}{n}{2^{bba}} \Witness{n}{n}{2^{bba}} \Witness{y}{n}{3^{aba}}
    \Witness{n}{n}{3^{aba}} \Witness{n}{y}{3^{aa}} \Witness{n}{y}{2^{aab}}
    \Witness{n}{y}{1^{ab}}
    \notag
\end{align}

\subsubsection{Presentation}
In order to establish a presentation of $\Pill_{1, 1}(\Monoid)$, we introduce an alternative
description of the clone congruence $\Equiv$ through a new equivalence relation $\Equiv'$.
For this, let us define $\Equiv'$ as the equivalence relation on $\P(\Monoid)$ generated by
\begin{equation} \label{equ:alternative_relation_pill_1}
    \PWord \Conc \Par{\alpha_1 \MonoidProductExtension \QWord} \Conc
    \RWord \Conc \RWord' \Conc
    \Par{\alpha_2 \MonoidProductExtension \QWord} \Conc
    \RWord'' \Conc \Par{\alpha_3 \MonoidProductExtension \QWord} \Conc \PWord'
    \Equiv'
    \PWord \Conc \Par{\alpha_1 \MonoidProductExtension \QWord} \Conc
    \RWord \Conc \Par{\alpha_2 \MonoidProductExtension \QWord} \Conc
    \RWord' \Conc \RWord'' \Conc
    \Par{\alpha_3 \MonoidProductExtension \QWord} \Conc \PWord',
\end{equation}
\begin{equation} \label{equ:alternative_relation_pill_2}
    \PWord \Conc \Par{\alpha_1 \MonoidProductExtension \QWord_1} \Conc
    \RWord \Conc
    \Par{\beta_1 \MonoidProductExtension \QWord_2} \Conc
    \Par{\alpha_2 \MonoidProductExtension \QWord_1} \Conc
    \RWord' \Conc
    \Par{\beta_2 \MonoidProductExtension \QWord_2} \Conc
    \PWord'
    \Equiv'
    \PWord \Conc \Par{\alpha_1 \MonoidProductExtension \QWord_1} \Conc
    \RWord \Conc
    \Par{\alpha_2 \MonoidProductExtension \QWord_1} \Conc
    \Par{\beta_1 \MonoidProductExtension \QWord_2} \Conc
    \RWord' \Conc
    \Par{\beta_2 \MonoidProductExtension \QWord_2} \Conc
    \PWord',
\end{equation}
where $\PWord, \PWord', \QWord, \QWord_1, \QWord_2, \RWord, \RWord', \RWord'' \in
\P(\Monoid)$ and $\alpha_1, \alpha_2, \alpha_3, \beta_1, \beta_2 \in \Monoid$.

\begin{Statement}{Lemma}{lem:alternative_relation_pill}
    For any monoid $\Monoid$, the binary relations $\Equiv$ and $\Equiv'$ on $\P(\Monoid)$
    are equal.
\end{Statement}
\begin{Proof}
    Let $\Leq$ be any total order relation on $\Monoid$. Let $\PWord, \PWord' \in
    \P(\Monoid)$ such that $\PWord \Equiv' \PWord'$. Since $\Equiv'$ is generated
    by~\eqref{equ:alternative_relation_pill_1} and~\eqref{equ:alternative_relation_pill_2},
    we have two cases to consider.
    \begin{enumerate}[label={\it (\Roman*)}]
        \item
        If $\PWord$ and $\PWord'$ decompose as
        \begin{math}
            \PWord =
            \PWord'' \Conc \Par{\alpha_1 \MonoidProductExtension \QWord} \Conc
            \RWord \Conc \RWord' \Conc
            \Par{\alpha_2 \MonoidProductExtension \QWord} \Conc
            \RWord'' \Conc \Par{\alpha_3 \MonoidProductExtension \QWord} \Conc \PWord'''
        \end{math}
        and
        \begin{math}
            \PWord' =
            \PWord'' \Conc \Par{\alpha_1 \MonoidProductExtension \QWord} \Conc
            \RWord \Conc \Par{\alpha_2 \MonoidProductExtension \QWord} \Conc
            \RWord' \Conc \RWord'' \Conc
            \Par{\alpha_3 \MonoidProductExtension \QWord} \Conc \PWord'''
        \end{math}
        where $\PWord'', \PWord''', \QWord, \RWord, \RWord', \RWord'' \in \P(\Monoid)$ and
        $\alpha_1, \alpha_2, \alpha_3, \in \Monoid$, then
        \begin{math}
            \First_1(\PWord)
            =
            \First_1 \Par{
                \PWord'' \Conc \Par{\alpha_1 \MonoidProductExtension \QWord} \Conc
                \RWord \Conc \RWord' \Conc
                \RWord'' \Conc \PWord'''
            }
            =
            \First_1\Par{\PWord'},
        \end{math}
        \begin{math}
            \Sort_{\Leq}\Par{\PWord} = \Sort_{\Leq}\Par{\PWord'},
        \end{math}
        and
        \begin{math}
            \First_1^{\Reverse}(\PWord)
            =
            \First_1 ^{\Reverse}\Par{
                \PWord'' \Conc
                \RWord \Conc \RWord' \Conc \RWord'' \Conc
                \Par{\alpha_3 \MonoidProductExtension \QWord} \Conc
                \PWord'''
            }
            =
            \First_1^{\Reverse}\Par{\PWord'}.
        \end{math}
        Therefore, $\PWord \Equiv \PWord'$.
        \item If $\PWord$ and $\PWord'$ decompose as
        \begin{math}
            \PWord =
            \PWord'' \Conc \Par{\alpha_1 \MonoidProductExtension \QWord_1} \Conc
            \RWord \Conc
            \Par{\beta_1 \MonoidProductExtension \QWord_2} \Conc
            \Par{\alpha_2 \MonoidProductExtension \QWord_1} \Conc
            \RWord' \Conc
            \Par{\beta_2 \MonoidProductExtension \QWord_2} \Conc
            \PWord'''
        \end{math}
        and
        \begin{math}
            \PWord' =
            \PWord'' \Conc \Par{\alpha_1 \MonoidProductExtension \QWord_1} \Conc
            \RWord \Conc
            \Par{\alpha_2 \MonoidProductExtension \QWord_1} \Conc
            \Par{\beta_1 \MonoidProductExtension \QWord_2} \Conc
            \RWord' \Conc
            \Par{\beta_2 \MonoidProductExtension \QWord_2} \Conc
            \PWord'''
        \end{math}
        where $\PWord'', \PWord''', \QWord_1, \QWord_2, \RWord, \RWord' \in \P(\Monoid)$ and
        $\alpha_1, \alpha_2, \beta_1, \beta_2 \in \Monoid$, then
        \begin{math}
            \First_1 \Par{\PWord} =
            \First_1 \Par{
                \PWord'' \Conc \Par{\alpha_1 \MonoidProductExtension \QWord_1} \Conc
                \RWord \Conc
                \Par{\beta_1 \MonoidProductExtension \QWord_2} \Conc
                \RWord' \Conc
                \PWord'''
            }
            =
            \First_1 \Par{\PWord'},
        \end{math}
        \begin{math}
            \Sort_{\Leq}\Par{\PWord} = \Sort_{\Leq}\Par{\PWord'},
        \end{math}
        and
        \begin{math}
            \First_1^{\Reverse} \Par{\PWord} =
            \First_1^{\Reverse} \Par{
                \PWord'' \Conc
                \RWord \Conc
                \Par{\alpha_2 \MonoidProductExtension \QWord_1} \Conc
                \RWord' \Conc
                \Par{\beta_2 \MonoidProductExtension \QWord_2} \Conc
                \PWord'''
            }
            =
            \First_1^{\Reverse} \Par{\PWord'}.
        \end{math}
        Therefore, $\PWord \Equiv \PWord'$.
    \end{enumerate}
    This shows that $\PWord \Equiv' \PWord'$ implies $\PWord \Equiv \PWord'$.

    Conversely, let $\PWord, \PWord' \in \P(\Monoid)$ such that $\PWord \Equiv \PWord'$. By
    Proposition~\ref{prop:covering_relation_pill_equivalence}, this is equivalent to the
    fact that $\PWord \EquivCovering \PWord'$. Since $\EquivCovering$ is generated by
    $\Covering$, we have three cases to explore depending whether $\PWord \Covering^{(1)}
    \PWord'$, $\PWord \Covering^{(2)} \PWord'$, or $\PWord \Covering^{(3)} \PWord'$,
    \begin{enumerate}[label={\it (\Roman*)}]
        \item If $\PWord \Covering^{(1)} \PWord'$, then $\PWord$ and $\PWord'$ decompose as
        \begin{math}
            \PWord =
            \PWord'' \Conc \Witness{u}{n}{i^{\alpha_1}} \Conc \QWord \Conc
            \Witness{n}{n}{i^{\alpha_2}} \Conc \PWord'''
        \end{math}
        and
        \begin{math}
            \PWord' =
            \PWord'' \Conc \Witness{u}{n}{i^{\alpha_1}} \Conc
            \Witness{n}{n}{i^{\alpha_2}} \Conc \QWord \Conc \PWord'''
        \end{math}
        where $\PWord'', \PWord''', \QWord \in \P(\Monoid)$, $i^{\alpha_1}, i^{\alpha_2} \in
        \SetLetters_\Monoid$, $\QWord \ne \PWordEmpty$, and $i \notin \QWord$. Since the
        position $\Length\Par{\PWord''} + \Length\Par{\QWord} + 2$ of $\PWord$ is not a
        right $1$-witness, there is necessarily an occurrence of an $\Monoid$-pigmented
        letter having $i$ as value in $\PWord'''$. Hence,
        \begin{math}
            \PWord =
            \PWord'' \Conc \Witness{u}{n}{i^{\alpha_1}} \Conc \QWord \Conc
            \Witness{n}{n}{i^{\alpha_2}} \Conc
            \PWord'''' \Conc i^{\alpha_3} \Conc \PWord'''''
        \end{math}
        and
        \begin{math}
            \PWord' =
            \PWord'' \Conc \Witness{u}{n}{i^{\alpha_1}} \Conc
            \Witness{n}{n}{i^{\alpha_2}} \Conc \QWord \Conc
            \PWord'''' \Conc i^{\alpha_3} \Conc \PWord'''''
        \end{math}
        where $\PWord'''', \PWord''''' \in \P(\Monoid)$ and $\alpha_3 \in \Monoid$.
        By~\eqref{equ:alternative_relation_pill_1}, we have $\PWord \Equiv' \PWord'$.
        \item If $\PWord \Covering^{(2)} \PWord'$, then $\PWord$ and $\PWord'$ decompose as
        \begin{math}
            \PWord =
            \PWord'' \Conc \Witness{n}{n}{i^{\alpha_1}} \Witness{n}{n}{i^{\alpha_2}} \Conc
            \PWord'''
        \end{math}
        and
        \begin{math}
            \PWord' =
            \PWord'' \Conc
            \Witness{n}{n}{i^{\alpha_2}} \Witness{n}{n}{i^{\alpha_1}} \Conc \PWord'''
        \end{math}
        where $\PWord'', \PWord''' \in \P(\Monoid)$, $i^{\alpha_1}, i^{\alpha_2} \in
        \SetLetters_\Monoid$, $\alpha_1 \ne \alpha_2$, and $\alpha_2 \Leq \alpha_1$. Since
        the position $\Length\Par{\PWord''} + 1$ of $\PWord$ is not a left $1$-witness,
        there is necessarily an occurrence of an $\Monoid$-pigmented letter having $i$ as
        value in $\PWord''$. Similarly, since the position $\Length\Par{\PWord''} + 2$ of
        $\PWord$ is not a right $1$-witness, there is necessarily an occurrence of an
        $\Monoid$-pigmented letter having $i$ as value in $\PWord'''$. Hence,
        \begin{math}
            \PWord =
            \QWord \Conc i^{\beta_1} \Conc \QWord'
            \Conc \Witness{n}{n}{i^{\alpha_1}} \Witness{n}{n}{i^{\alpha_2}} \Conc
            \QWord'' \Conc i^{\beta_2} \Conc \QWord'''
        \end{math}
        and
        \begin{math}
            \PWord' =
            \QWord \Conc i^{\beta_1} \Conc \QWord'
            \Conc \Witness{n}{n}{i^{\alpha_2}} \Witness{n}{n}{i^{\alpha_1}} \Conc
            \QWord'' \Conc i^{\beta_2} \Conc \QWord'''
        \end{math}
        where $\QWord, \QWord', \QWord'', \QWord''' \in \P(\Monoid)$ and $\beta_1, \beta_2
        \in \Monoid$. By~\eqref{equ:alternative_relation_pill_2}, we have $\PWord \Equiv'
        \PWord'$.
        \item If $\PWord \Covering^{(3)} \PWord'$, then $\PWord$ and $\PWord'$ decompose as
        \begin{math}
            \PWord =
            \PWord'' \Conc \Witness{n}{y}{i_1^{\alpha_1}} \Witness{u}{n}{i_2^{\alpha_2}}
            \Conc \PWord'''
        \end{math}
        and
        \begin{math}
            \PWord' =
            \PWord'' \Conc \Witness{u}{n}{i_2^{\alpha_2}} \Witness{n}{y}{i_1^{\alpha_1}}
            \Conc \PWord'''
        \end{math}
        where $\PWord'', \PWord''' \in \P(\Monoid)$, $i_1^{\alpha_1}, i_2^{\alpha_2} \in
        \SetLetters_\Monoid$, and $i_1 \ne i_2$. Since the position $\Length\Par{\PWord''} +
        1$ of $\PWord$ is not a left $1$-witness, there is necessarily an occurrence of an
        $\Monoid$-pigmented letter having $i_1$ as value in $\PWord''$. Similarly, since the
        position $\Length\Par{\PWord''} + 2$ of $\PWord$ is not a right $1$-witness, there
        is necessarily an occurrence of an $\Monoid$-pigmented letter having $i_2$ as value
        in $\PWord'''$. Hence,
        \begin{math}
            \PWord =
            \QWord \Conc i_1^{\beta_1} \Conc \QWord'
            \Conc \Witness{n}{y}{i_1^{\alpha_1}} \Witness{u}{n}{i_2^{\alpha_2}} \Conc
            \QWord'' \Conc i_2^{\beta_2} \Conc \QWord'''
        \end{math}
        and
        \begin{math}
            \PWord' =
            \QWord \Conc i_1^{\beta_1} \Conc \QWord'
            \Conc \Witness{u}{n}{i_2^{\alpha_2}} \Witness{n}{y}{i_1^{\alpha_1}} \Conc
            \QWord'' \Conc i_2^{\beta_2} \Conc \QWord'''
        \end{math}
        where $\QWord, \QWord', \QWord'', \QWord''' \in \P(\Monoid)$ and $\beta_1, \beta_2
        \in \Monoid$. By~\eqref{equ:alternative_relation_pill_2}, we have $\PWord \Equiv'
        \PWord'$.
    \end{enumerate}
    This shows that $\PWord \Equiv \PWord'$ implies $\PWord \Equiv' \PWord'$ and establishes
    the statement of the lemma.
\end{Proof}

\begin{Statement}{Theorem}{thm:presentation_pill}
    For any monoid $\Monoid$, the clone $\Pill_{1, 1}(\Monoid)$ admits the presentation
    $\Par{\GeneratingSet_\Monoid, \RelationSet_\Monoid'}$ where $\RelationSet_\Monoid'$ is
    the set $\RelationSet_\Monoid$ augmented with the $\GeneratingSet_\Monoid$-equations
    \begin{equation} \label{equ:presentation_pill_1}
        \RightComb_\Monoid\Par{
            1^{\alpha_1} 2^\MonoidUnit 3^\MonoidUnit 1^{\alpha_2} 4^\MonoidUnit 1^{\alpha_3}
        }
        \ \RelationSet_\Monoid' \
        \RightComb_\Monoid\Par{
            1^{\alpha_1} 2^\MonoidUnit 1^{\alpha_2} 3^\MonoidUnit 4^\MonoidUnit 1^{\alpha_3}
        },
    \end{equation}
    \begin{equation} \label{equ:presentation_pill_2}
        \RightComb_\Monoid\Par{1^{\alpha_1} 2^\MonoidUnit 3^{\beta_1} 1^{\alpha_2}
            4^\MonoidUnit 3^{\beta_2}}
        \ \RelationSet_\Monoid' \
        \RightComb_\Monoid\Par{1^{\alpha_1} 2^\MonoidUnit 1^{\alpha_2} 3^{\beta_1}
            4^\MonoidUnit 3^{\beta_2}},
    \end{equation}
    where $\alpha_1, \alpha_2, \alpha_3, \beta_1, \beta_2 \in \Monoid$ and $\MonoidUnit$ is
    the unit of $\Monoid$.
\end{Statement}
\begin{Proof}
    Let $\Equiv''$ be the clone congruence of $\P(\Monoid)$ generated by
    \begin{equation} \label{equ:presentation_pill_3}
        1^{\alpha_1} 2^\MonoidUnit 3^\MonoidUnit 1^{\alpha_2} 4^\MonoidUnit 1^{\alpha_3}
        \Equiv''
        1^{\alpha_1} 2^\MonoidUnit 1^{\alpha_2} 3^\MonoidUnit 4^\MonoidUnit 1^{\alpha_3},
    \end{equation}
    \begin{equation} \label{equ:presentation_pill_4}
        1^{\alpha_1} 2^\MonoidUnit 3^{\beta_1} 1^{\alpha_2} 4^\MonoidUnit 3^{\beta_2}
        \Equiv''
        1^{\alpha_1} 2^\MonoidUnit 1^{\alpha_2} 3^{\beta_1} 4^\MonoidUnit 3^{\beta_2},
    \end{equation}
    with $\alpha_1, \alpha_2, \alpha_3, \beta_1, \beta_2 \in \Monoid$. Let us show that the
    congruence $\Equiv$ and $\Equiv''$ of $\P(\Monoid)$ are equal. This will imply, by the
    remark stated in Section~\ref{subsubsec:remark_presentations}, that $\P(\Monoid)
    /_{\Equiv''} = \P(\Monoid) /_{\Equiv} = \Pill_{1, 1}(\Monoid)$ admits the stated
    presentation.

    First, since the images by the map $\First_1$ of the left-hand side and the right-hand
    side of~\eqref{equ:presentation_pill_3} are both equal to
    \begin{math}
        1^{\alpha_1} 2^\MonoidUnit 3^\MonoidUnit 4^\MonoidUnit,
    \end{math}
    the images by the map $\First_1^{\Reverse}$ of the left-hand side and the right-hand
    side of~\eqref{equ:presentation_pill_3} are both equal to
    \begin{math}
        2^\MonoidUnit 3^\MonoidUnit 4^\MonoidUnit 1^{\alpha_3},
    \end{math}
    and the image by $\Sort_{\Leq}$, where $\Leq$ is any total order relation on $\Monoid$,
    of the left-hand side and the right-hand side of~\eqref{equ:presentation_pill_3} are the
    same, we have
    \begin{equation}
        1^{\alpha_1} 2^\MonoidUnit 3^\MonoidUnit 1^{\alpha_2} 4^\MonoidUnit 1^{\alpha_3}
        \Equiv
        1^{\alpha_1} 2^\MonoidUnit 1^{\alpha_2} 3^\MonoidUnit 4^\MonoidUnit 1^{\alpha_3}.
    \end{equation}
    Moreover, since the images by the map $\First_1$ of the left-hand side and the right-hand
    side of~\eqref{equ:presentation_pill_4} are both equal to
    \begin{math}
        1^{\alpha_1} 2^\MonoidUnit 3^{\beta_1} 4^\MonoidUnit,
    \end{math}
    the images by the map $\First_1^{\Reverse}$ of the left-hand side and the right-hand
    side of~\eqref{equ:presentation_pill_4} are both equal to
    \begin{math}
        2^\MonoidUnit 1^{\alpha_2} 4^\MonoidUnit 3^{\beta_2},
    \end{math}
    and the image by $\Sort_{\Leq}$ of the left-hand side and the right-hand side
    of~\eqref{equ:presentation_pill_4} are the same, we have
    \begin{equation}
        1^{\alpha_1} 2^\MonoidUnit 3^{\beta_1} 1^{\alpha_2} 4^\MonoidUnit 3^{\beta_2}
        \Equiv
        1^{\alpha_1} 2^\MonoidUnit 1^{\alpha_2} 3^{\beta_1} 4^\MonoidUnit 3^{\beta_2}.
    \end{equation}
    This shows that $\Equiv''$ is contained in $\Equiv$.

    To prove that $\Equiv$ is contained in $\Equiv''$, let us show that $\Equiv'$ is
    contained in $\Equiv''$. By Lemma~\ref{lem:alternative_relation_pill}, the targeted
    property will follow. For any $\PWord, \PWord', \QWord, \RWord, \RWord', \RWord'' \in
    \P(\Monoid)$ and $\alpha_1, \alpha_2, \alpha_3 \in \Monoid$, we have
    \begin{align} \label{equ:presentation_pill_5}
        \PWord \Conc \Par{\alpha_1 \MonoidProductExtension \QWord} \Conc
        \RWord \Conc \RWord' \Conc
        \Par{\alpha_2 \MonoidProductExtension \QWord} \Conc
        \RWord'' \Conc \Par{\alpha_3 \MonoidProductExtension \QWord} \Conc \PWord'
        & =
        1^\MonoidUnit 2^\MonoidUnit 3^\MonoidUnit
        \Han{
            \PWord,
            1^{\alpha_1} 2^\MonoidUnit 3^\MonoidUnit 1^{\alpha_2} 4^\MonoidUnit 1^{\alpha_3}
            \Han{
                \QWord, \RWord, \RWord', \RWord''
            },
            \PWord'
        }
        \\
        & \Equiv''
        1^\MonoidUnit 2^\MonoidUnit 3^\MonoidUnit
        \Han{
            \PWord,
            1^{\alpha_1} 2^\MonoidUnit 1^{\alpha_2} 3^\MonoidUnit 4^\MonoidUnit 1^{\alpha_3}
            \Han{
                \QWord, \RWord, \RWord', \RWord''
            },
            \PWord'
        }
        \notag
        \\
        & =
        \PWord \Conc \Par{\alpha_1 \MonoidProductExtension \QWord} \Conc
        \RWord \Conc \Par{\alpha_2 \MonoidProductExtension \QWord} \Conc
        \RWord' \Conc \RWord'' \Conc
        \Par{\alpha_3 \MonoidProductExtension \QWord} \Conc \PWord'
        \notag
    \end{align}
    so that the first and last members of~\eqref{equ:presentation_pill_5} are
    $\Equiv''$-equivalent. Moreover, for any $\PWord, \PWord', \QWord_1, \QWord_2, \RWord,
    \RWord' \in \P(\Monoid)$ and $\alpha_1, \alpha_2, \beta_1, \beta_2 \in \Monoid$, we have
    \begin{align} \label{equ:presentation_pill_6}
        \PWord \Conc \Par{\alpha_1 \MonoidProductExtension \QWord_1} \Conc
        \RWord \Conc
        \Par{\beta_1 \MonoidProductExtension \QWord_2} \Conc
        \Par{\alpha_2 \MonoidProductExtension \QWord_1}
        & \Conc
        \RWord' \Conc
        \Par{\beta_2 \MonoidProductExtension \QWord_2} \Conc
        \PWord'
        \\
        & =
        1^\MonoidUnit 2^\MonoidUnit 3^\MonoidUnit
        \Han{
            \PWord,
            1^{\alpha_1} 2^\MonoidUnit 3^{\beta_1} 1^{\alpha_2} 4^\MonoidUnit 3^{\beta_2}
            \Han{
                \QWord_1, \RWord, \QWord_2, \RWord'
            },
            \PWord'
        }
        \notag
        \\
        & \Equiv''
        1^\MonoidUnit 2^\MonoidUnit 3^\MonoidUnit
        \Han{
            \PWord,
            1^{\alpha_1} 2^\MonoidUnit 1^{\alpha_2} 3^{\beta_1} 4^\MonoidUnit 3^{\beta_2}
            \Han{
                \QWord_1, \RWord, \QWord_2, \RWord'
            },
            \PWord'
        }
        \notag
        \\
        & =
        \PWord \Conc \Par{\alpha_1 \MonoidProductExtension \QWord_1} \Conc
        \RWord \Conc
        \Par{\alpha_2 \MonoidProductExtension \QWord_1} \Conc
        \Par{\beta_1 \MonoidProductExtension \QWord_2} \Conc
        \RWord' \Conc
        \Par{\beta_2 \MonoidProductExtension \QWord_2} \Conc
        \PWord'
        \notag
    \end{align}
    so that the first and last members of~\eqref{equ:presentation_pill_6} are
    $\Equiv''$-equivalent. Since $\Equiv'$ is the equivalence relation generated
    by~\eqref{equ:alternative_relation_pill_1} and~\eqref{equ:alternative_relation_pill_2},
    the targeted property is shown. This establishes the statement of the theorem.
\end{Proof}

By Theorem~\ref{thm:presentation_pill}, any $\Pill_{1, 1}(\Monoid)$-algebra is, up to term
equivalence, an $\Monoid$-pigmented monoid $\Par{\Algebra, \Product, \Identity,
\Par{\Pigmentation_\alpha}_{\alpha \in \Monoid}}$ where $\Product$ and
$\Par{\Pigmentation_\alpha}_{\alpha \in \Monoid}$ satisfy, by spelling
out~\eqref{equ:presentation_pill_1} and~\eqref{equ:presentation_pill_2} and simplifying them
modulo the background theory $\Equiv_{\RelationSet_{\Monoid}}$,
\begin{equation}
    \Pigmentation_{\alpha_1}\Par{x_1} \Product x_2 \Product x_3 \Product
    \Pigmentation_{\alpha_2}\Par{x_1} \Product x_4 \Product
    \Pigmentation_{\alpha_3}\Par{x_1}
    =
    \Pigmentation_{\alpha_1}\Par{x_1} \Product x_2 \Product
    \Pigmentation_{\alpha_2}\Par{x_1} \Product x_3 \Product x_4 \Product
    \Pigmentation_{\alpha_3}\Par{x_1},
\end{equation}
\begin{equation}
    \Pigmentation_{\alpha_1}\Par{x_1} \Product x_2 \Product \Pigmentation_{\beta_1}\Par{x_3}
    \Product \Pigmentation_{\alpha_2}\Par{x_1} \Product x_4 \Product
    \Pigmentation_{\beta_2}\Par{x_3}
    =
    \Pigmentation_{\alpha_1}\Par{x_1} \Product x_2 \Product
    \Pigmentation_{\alpha_2}\Par{x_1} \Product \Pigmentation_{\beta_1}\Par{x_3} \Product x_4
    \Product \Pigmentation_{\beta_2}\Par{x_3},
\end{equation}
for any $x_1, x_2, x_3, x_4 \in \Algebra$ and $\alpha_1, \alpha_2, \alpha_3, \beta_1,
\beta_2 \in \Monoid$.

\section{Open questions and future work}
We have introduced the construction $\P$ producing clones from monoids and studied a
selection of quotient clones of $\P(\Monoid)$. This has resulted in a novel hierarchy of
clone realizations of varieties of monoids. Here follow some open questions and future areas
of investigation raised by this work.

\subsection*{Variations around the variety of pigmented monoids}
As shown by Theorem~\ref{thm:construction_p}, $\P(\Monoid)$ is a clone realization of the
variety of $\Monoid$-pigmented monoids. This variety stems from the six
equations~\eqref{equ:pigmented_monoids_1}, \eqref{equ:pigmented_monoids_2},
\eqref{equ:pigmented_monoids_3}, \eqref{equ:pigmented_monoids_4},
\eqref{equ:pigmented_monoids_5}, and~\eqref{equ:pigmented_monoids_6}. A compelling question
to consider involves the alternative varieties resulting from the omission of some of these
equations, and proposing in this way variations of the construction $\P$ in order to
describe the corresponding clone realizations. There are therefore $2^6 - 1 = 63$ such
alternative varieties and among these, $2^3 - 1 = 7$ seem particularly interesting to study
because these equations are naturally paired as outlined at the end of
Section~\ref{subsubsec:varieties_pigmented_monoids}. Indeed, \eqref{equ:pigmented_monoids_1}
pairs with~\eqref{equ:pigmented_monoids_2}, \eqref{equ:pigmented_monoids_3}
with~\eqref{equ:pigmented_monoids_4}, and \eqref{equ:pigmented_monoids_5}
with~\eqref{equ:pigmented_monoids_6}. In particular, in~\cite{Gir18} (see
also~\cite{Gir17,Gir20}), the variety that arises by omitting the pair consisting of
Equations~\eqref{equ:pigmented_monoids_3} and~\eqref{equ:pigmented_monoids_4} (except for
a few details) has been studied via operads and involves configurations of noncrossing and
decorated diagonals in polygons. Such objects recur very frequently in
combinatorics~\cite{CP92,FN99,DRS10,PR14} and considering clone structures on these objects
could give an original point of view and lead to new questions and results in this domain.

\subsection*{Linearization of the construction and equations}
The clones examined in this work are defined within the category of sets. It is of course
possible to extend the construction $\P$ in order to see the produced clones as clones on
the $\K$-linear span of the set of $\Monoid$-pigmented words where $\K$ is any field of zero
characteristic. This type of extension opens a myriad of new questions. Among these, the
broad question of describing the nontrivial equations satisfied by certain linear
combinations of terms of the variety of $\Monoid$-pigmented monoids is worth considering.
When translated into the language of clones, this equates to describe the presentations of
certain subclones of the linearization of $\P(\Monoid)$ which are generated by some linear
combinations of $\Monoid$-pigmented words. More specifically, this question can be posed,
given $\alpha_1, \alpha_2 \in \Monoid$, for the commutator $1^{\alpha_1} 2^{\alpha_2} -
2^{\alpha_2} 1^{\alpha_1}$ and for the anti-commutator $1^{\alpha_1} 2^{\alpha_2} +
2^{\alpha_2} 1^{\alpha_1}$ in the linearization of $\P(\Monoid)$, as well as in the
linearizations of some of its quotients constructed in
Sections~\ref{subsec:clone_congruences} and~\ref{sec:hierarchy}. Similar questions have been
explored for different varieties of algebras: for instance for the anti-commutator of
associative algebras~\cite{Gle70}, for the commutator and anti-commutator of bicommutative
algebras~\cite{DI18}, and for the anti-commutator of pre-Lie algebras~\cite{BL11}.

\subsection*{Finitely generated subclones}
In the present work, the clone $\P(\Monoid)$ is studied along with some of its quotients. A
potential next step in this research involves paying attention to subclones of $\P(\Monoid)$
and to some of its quotients generated by some finite sets of elements. This approach has
been considered in~\cite{Gir15} where a construction $\T$ from monoids to operads has been
introduced and numerous operads on combinatorial objects have been discovered (on several
sorts of words, trees, and paths). Recall, as explained in
Section~\ref{subsubsec:varieties_pigmented_monoids}, that the construction $\P$ can be seen
as a generalization of the construction $\T$ at the level of clones. In this way, we could
expect to develop a hierarchy of clones based on a large collection of sorts of
combinatorial objects. As consequences, mainly by describing presentations of such derived
clones, it may sometimes be feasible to establish a convergent rewrite system on the terms
of the underlying variety. This could lead to new methods for the enumeration of the
involved combinatorial objects and for their ---exhaustive or random--- generation
(see~\cite{Gir19} and~\cite{Gir20b} in the context of operads rather than clones).

\subsection*{Plactic-like monoids and other constructions}
As briefly highlighted in Section~\ref{subsec:p_symbols}, many monoids hold a distinctive
role in algebraic combinatorics. Examples include the plactic monoid~\cite{LS81,Lot02}, the
hypoplactic monoid~\cite{KT97}, the sylvester monoid~\cite{HNT05}, the Bell
monoid~\cite{Rey07}, the Baxter monoid~\cite{Gir12}, the $k$-recoil monoid~\cite{NRT11}, and
the stalactic monoid~\cite{HNT08}. These monoids can be defined through congruences of free
monoids on a totally ordered alphabet. The main observation here is that these monoids
intervene in a crucial way to construct Hopf algebras generalizing the prototypical one of
symmetric functions~\cite{GKLLRT94} (also refer to the previously cited works
and~\cite[Chap.~5]{Gir11} for a comprehensive description and properties of this
construction). A key component here is formed by $\PSymbol$-symbols, which ---akin to the
present work--- are maps sending words to some combinatorial objects encoding the
equivalence classes. In the context of the present work, we are interested in clone
congruences of $\P(\Monoid)$, which are in particular also monoid congruences on words of
integers. As a matter of fact, most of the previously cited congruences do not define clone
congruences of $\P(\Monoid)$. Nevertheless, instead of trying to use already existing
monoids to propose new clone congruences of $\P(\Monoid)$ (which is a possible direction for
future work that deserves to be explored), we can proceed in the opposite direction. This
consists in trying to build Hopf algebras in the same manner by considering the clone
congruences and monoids at the heart of the constructions of $\Arra_k(\Monoid)$, $\Magn_{k,
k'}(\Monoid)$, $\Stal_k(\Monoid)$, and $\Pill_{k, k'}(\Monoid)$.

\subsection*{General case for pigmented magnets and pigmented pillars}
The final question we ask here concerns the clones $\Magn_{k, k'}(\Monoid)$ and $\Pill_{k,
k'}(\Monoid)$. These clones are well understood in the case $k = 1 = k'$. Indeed, both
descriptions and presentations are furnished for each clone in this case. The question here
consists in establishing generalizations of these results working for any nonnegative
integers $k$ and~$k'$.

\MakeReferences

\end{document}